 \documentclass[a4paper]{article}
 \usepackage{etoolbox}
 \newbool{PREPRINT}
 \booltrue{PREPRINT}

\ifbool{PREPRINT}{ 
\usepackage[colorlinks,bookmarksopen,bookmarksnumbered,citecolor=red,urlcolor=red]{hyperref}
\usepackage[affil-it]{authblk}

\usepackage{fullpage}
}{
\usepackage{fullpage}
} 

\usepackage{amsmath,amssymb}
\usepackage{stmaryrd}
\usepackage{graphicx}
\usepackage{algorithm,algorithmic}
\usepackage[usenames]{color}
\usepackage[colorlinks,bookmarksopen,bookmarksnumbered,citecolor=red,urlcolor=red]{hyperref}
\usepackage{graphicx}
\usepackage{multirow}
\usepackage[dvipsnames]{xcolor}
\newcommand{\figdir}{./}

\newcommand{\order}{\mathcal{O}}
\renewcommand{\vec}[1]{\boldsymbol{#1}}
\newcommand{\vecn}[1]{\boldsymbol{{#1}}}
\title{Matrix-free multigrid block-preconditioners for higher order Discontinuous Galerkin discretisations}

\ifbool{PREPRINT}{
\author[1]{Peter~Bastian}
\author[$\dagger$,*,2]{Eike~Hermann~M\"uller}
\author[1]{Steffen~M\"{u}thing}
\author[1]{Marian~Piatkowski}
\affil[1]{Interdisciplinary Centre for Scientific Computing, Heidelberg University, Im Neuenheimer Feld 205, 69120 Heidelberg, Germany}
\affil[2]{Department of Mathematical Sciences, University of Bath, Bath BA2 7AY, United Kingdom}
\affil[$\dagger$]{Email: \texttt{e.mueller@bath.ac.uk}}
\affil[*]{Corresponding author}

}{ 
\author[heidelberg]{Peter~Bastian}
\author[bath]{Eike~Hermann~M\"{u}ller\corref{cor1}\fnref{fn1}}
\author[heidelberg]{Steffen~M\"{u}thing}
\author[heidelberg]{Marian~Piatkowski}
\ead{e.mueller@bath.ac.uk}
\cortext[cor1]{Corresponding author}
\fntext[fn1]{email: \texttt{e.mueller@bath.ac.uk}}
\address[bath]{Department of Mathematical Sciences, University of Bath, Claverton Down, Bath BA2 7AY, United Kingdom}
\address[heidelberg]{Interdisciplinary Centre for Scientific Computing, Heidelberg University, Im Neuenheimer Feld 205, 69120 Heidelberg, Germany}
} 
\date{\today}
\begin{document}
\ifbool{PREPRINT}{ 
\maketitle
}{} 
\begin{abstract}
Efficient and suitably preconditioned iterative solvers for elliptic
partial differential equations (PDEs) of the convection-diffusion type
are used in all fields of science and engineering, including for example
computational fluid dynamics, nuclear reactor simulations and combustion
models. To achieve optimal
performance, solvers have to exhibit high arithmetic intensity and need to
exploit every form of parallelism available in modern manycore CPUs. This
includes both distributed- or shared memory parallelisation between processors
and vectorisation on individual cores.

The computationally most expensive components of the solver are the repeated
applications of the linear operator and the preconditioner. For discretisations
based on higher-order Discontinuous Galerkin methods,
sum-factorisation results in a dramatic reduction of the
computational complexity of the operator application while, at the same time,
the matrix-free implementation can run at a
significant fraction of the theoretical peak floating point performance.
Multigrid methods for high order methods often rely on block-smoothers to reduce
high-frequency error components within one grid cell.
Traditionally,
this requires the assembly and expensive dense matrix solve in each
grid cell, which counteracts any improvements achieved in the fast
matrix-free operator application. To overcome this issue, we present a
new matrix-free implementation of block-smoothers. Inverting the block
matrices iteratively avoids storage and factorisation of the matrix
and makes it is possible to harness the full power of the CPU. We
implemented a hybrid multigrid algorithm with matrix-free
block-smoothers in the high order Discontinuous Galerkin (DG) space combined
with a low order coarse grid correction using algebraic multigrid where
only low order components are explicitly assembled.
The effectiveness of this approach is demonstrated by solving a set of
representative elliptic PDEs of increasing complexity, including a
convection dominated problem and the stationary SPE10 benchmark.
\end{abstract}

\ifbool{PREPRINT}{ 
%
%
\textbf{keywords}:
\newcommand{\sep}{,}
}{ 
\begin{keyword}
} 
multigrid\sep\ elliptic PDE\sep\ Discontinuous Galerkin\sep\ matrix-free methods\sep\ preconditioners\sep\ DUNE
\ifbool{PREPRINT}{}{ 
\end{keyword}
} 
\maketitle
\section{Introduction}
Second order elliptic PDEs of the convection-diffusion-reaction form
\begin{equation}
  -\nabla\cdot\left(K\nabla u\right) + \nabla\cdot(\vec{b} u) + c u = f\label{eqn:elliptic_eqn}
\end{equation}
with spatially varying coefficients
play an important role in many areas of science and engineering. The
Discontinuous Galerkin (DG) method
\cite{Wheeler1978,Arnold1982,Oden1998,Baumann1999,Riviere1999,Arnold2002,Riviere2008,DiPietro2011}
is a popular discretisation scheme for various reasons: it allows
local mass conservation and (when used with an appropriate basis) leads to
diagonal mass matrices which avoid global solves in explicit time integrators.
It also has several computational advantages as will be discussed below.
Problems of the
form in \eqref{eqn:elliptic_eqn} arise for example in
the simulation of subsurface flow phenomena
\cite{Cockburn2002,Bastian2003,Bastian2014II,Siefert2014}. Another important
application area is fluid dynamics. When the incompressible Navier
Stokes equations are solved with Chorin's projection method
\cite{Chorin1968,Temam1969}, a Poisson equation has to be solved for
the pressure correction and the momentum equation needs to be solved for the
tentative velocity, see \cite{2016arXiv161200657P} for a DG-based version.
Implicit DG solvers for the fully compressible atmospheric
equations of motion are described in \cite{Restelli2009,Dedner2015}.

The fast numerical solution of equation \eqref{eqn:elliptic_eqn}
requires algorithmically optimal algorithms, such as multigrid methods
(see \cite{Trottenberg2000} for an overview), which have a numerical
complexity that typically grows linearly with the number of
unknowns. Equally important is the efficient implementation of those
methods on modern computer hardware. While low order discretisation
methods such as finite volume approximations and linear conforming
finite element methods are very popular, from a computational point of
view one of their main drawbacks is the fact that their performance is
limited by main memory bandwidth. Modern computer
architectures such as the Intel Haswell/Broadwell CPUs, Xeon Phi
co-processors and graphics processing units (GPUs) may only achieve their
peak performance floating point performance if data loaded into the CPU is reused,
i.e. the arithmetic intensity of the algorithm
(ratio of floating point operations executed per byte loaded, in short: flop-to-byte-ratio)
is sufficiently high.
For example the 16 core Haswell chip
used for the computations in this work has a theoretical peak
performance of $243.2\cdot 10^9$ fused multiply-add (FMA) operations per
second and a main memory bandwidth (as measured with the STREAM
benchmark \cite{McCalpin1995}) of $40-50\operatorname{GB/s}$. Hence
about 40 FMA operations can be carried out in the time it takes to
read one double precision number from memory. Low order discretisation
methods require frequent sparse matrix-vector products. If the matrix
is stored in the standard compressed sparse row storage (CSR) format,
this multiplication has a very low arithmetic intensity of about 1 FMA
operation per matrix entry loaded from memory (assuming perfect caching of the vectors).
%
Hence traditional solvers which first assemble the system matrix and
then invert it with sparse matrix-vector operations will only use a
fraction of the system's peak performance.

Matrix-free methods, on the other hand, recompute each matrix-element
on-the-fly instead of storing an assembled matrix. In addition to
eliminating setup costs for the matrix assembly, this raises the
computational intensity and allows utilisation of the full performance
of modern manycore chip architectures. An additional advantage are the
reduced storage requirements: since the matrix does not have to be
kept in memory, this allows the solution of much larger
problems and has been quantified by estimating the memory requirements
for different implementations in this paper.
While high arithmetic intensity is a desirable property of a
good implementation, the quantity which really needs to be minimised
is the total computation time. A naive matrix-free implementation can
result in an increase in the runtime since a large number of floating
point operations (FLOPs) is executed for each unknown. Careful matrix-free
implementations of higher-order finite element discretisations avoid
this problem. If $p$ is the polynomial degree of the basis function in
each $d$-dimensional grid cell, a naive implementation requires
$\order(p^{2d})$ FLOPs per cell. Sum factorisation, which
exploits the tensor-product structure of the local basis
functions, reduces this to $\order(d\cdot p^{d+1})$
\cite{Orszag1980,Karniadakis2005,Vos2010,Kronbichler2012}, making the matrix-free implementation particularly
competitive if the order $p$ is sufficiently high.
Numerical results indicate that sum factorisation reduces the runtime even for 
low polynomial degrees $p>1$ \cite{Kronbichler2012,Muething2017}.
 Another advantage of the DG
method is that all data associated with one grid cell is stored
contiguously in memory. This avoids
additional gather/scatter operations which are necessary if data is
also held on other topological entities such as vertices and edges as
in conforming high order finite elements.

The efficient matrix-free implementation of higher-order DG methods
was described by some of us in a previous paper \cite{Muething2017} which
concentrates on the efficient application of the DG operator. Here we
extend this approach to iterative solvers for elliptic PDEs. While matrix-vector
products are one key ingredient in Krylov-subspace solvers, efficient
preconditioners are equally important and often form the computational
bottleneck. Many preconditioners such as algebraic multigrid (AMG)
\cite{Brandt1984,Stueben2001} or incomplete factorisations (ILU)
\cite{Saad2003} require the assembly of the system matrix, which can
offset the advantage of the matrix-free sparse matrix-vector product.

In this work we describe the matrix-free implementation of block-smoothers,
which are the key ingredient of multigrid methods such as those
described in \cite{Bastian2012,Siefert2014} . The idea behind the multigrid
algorithm in \cite{Bastian2012} is that the block-smoother eliminates
high-frequency errors within one grid cell whereas an AMG solver in a
lower order subspace of the full DG space reduces long range error
components. Due to the much smaller size of the coarse level system, the cost
for one AMG V-cycle on the low order subspace (which uses a traditional
matrix-based implementation) is negligible in comparison to the block smoother
on the high order system. In the simplest case the smoother is a block-Jacobi
iteration applied to the system $A\vecn{u}=\vecn{f}$ which arises from high-order DG discretisation
of \eqref{eqn:elliptic_eqn}. In each cell $T$ of
the grid this requires the following operations:\\
\begin{center}
\begin{minipage}{\linewidth}
\begin{center}
\begin{algorithmic}[1]
  \STATE{Calculate local defect $\vecn{d}_T = \vecn{f}_T - \sum_{T', \partial T'\cap \partial T\ne\emptyset} A_{T,T'}\vecn{u}^{(k)}_{T'}$}
  \STATE{Calculate correction $\vecn{\delta u}_T = A_{T,T}^{-1}\vecn{d}_T$}
  \STATE{Update solution $\vecn{u}^{(k+1)}_T=\vecn{u}^{(k)}_T+\vecn{\delta u}_T$}
\end{algorithmic}
\end{center}
\end{minipage}
\end{center}
Here $\vecn{u}_T$ is the solution on cell $T$ and the matrix $A_{T,T'}$
describes the coupling between cells $T$ and $T'$ which share at least one common face.
In addition to dense matrix-vector products to calculate the local defect $\vecn{d}_T$, this
requires the solution of a small dense system
\begin{equation}
A_{T,T} \vecn{\delta u}_T=\vecn{d}_T \label{eqn:small_dense_system}
\end{equation}
with $n_T=(p+1)^d$ unknowns in the second step.  Matrix based methods
use an LU- or Cholesky- factorisation of the matrix $A_{T,T}$ to solve
\eqref{eqn:small_dense_system} in each grid cell. The cost of the
back-substitution step, which has to be carried out for every
block-Jacobi iteration, is $\order(n_T^2)=\order(p^{2d})$. There is an
additional setup cost of $\order(n_T^3)=\order(p^{3d})$ for the initial
factorisation of $A_{T,T}$, which may be relevant if the number of preconditioner
applications is small. In contrast, the matrix-free method which
we present here uses an iterative method to solve this equation and if
the implementation uses sum-factorisation, this reduces the complexity
to $\order(d\cdot n_{\operatorname{iter}}(\epsilon)\cdot p^{d+1})$. An important factor
in this expression is $n_{\operatorname{iter}}(\epsilon)$, the number of
iterations which are required to solve the system in
\eqref{eqn:small_dense_system} to accuracy $\epsilon$. Solving
\eqref{eqn:small_dense_system} exactly, i.e. iterating until the
error is reduced to machine precision, will result in a relatively
large $n_{\operatorname{iter}}$ and might not lead to a reduction of
the total solution time for the full system. However, since the
inverse is only required in the preconditioner, an approximate
inversion with a small number of iterations might be sufficient and
can increase the computational efficiency of the overall method. This
is one of the key optimisations which will be explored below. The
overall efficiency of the implementation also depends on how computationally
expensive the evaluation of the functions $K$, $\vec{b}$ and $c$ in
\eqref{eqn:elliptic_eqn} is. At first sight this might put a
matrix-free preconditioner at a disadvantage, since those functions
have to be re-computed at every quadrature point in each
iteration. However, it is often sufficient to use a simpler functional
dependency in the preconditioner, for example by approximating $K$,
$\vec{b}$ and $c$ by their average value on each cell.
To demonstrate the efficiency of this approach we solve a set of
linear problems of increasing complexity. This includes advection dominated problems which require more advanced SOR-type smoothers and the stationary SPE10 benchmark. The latter is a hard problem since the permeability field $K(\vec{x})$ varies over several magnitudes and has high-contrast jumps between grid cells.


A similar fully matrix-free approach for the solution of elliptic
problems is described in \cite{Kronbichler2016}. The results there are
obtained with the \texttt{deal.II} \cite{dealII84} implementation
discussed in \cite{Kronbichler2012}. As for our code, the matrix-free
operator application in \cite{Kronbichler2012} is highly optimised and
utilises a significant fraction the peak floating point
performance. In contrast to our implementation, the authors of
\cite{Kronbichler2016} use a
(matrix-free) Chebyshev-iteration as a smoother on the DG space (a
Chebyshev smoother is also used for the AMG solver on the lower order
subspace). One drawback of this approach is
that it requires estimates of the largest eigenvalues,
which will result in additional setup costs. While the authors of
\cite{Kronbichler2016} solve a Laplace equation on complex geometries,
it is not clear how the same Chebyshev smoother can be applied to
non-trivial convection-dominated problems or problems with 
large, high contrast jumps of the coefficients.

Another approach for preconditioning higher order finite element
discretisations was originally developed for spectral methods in
\cite{Orszag1980} and used for example in
\cite{Fischer1997,Brown2010}. In \cite{Brown2010}, a preconditioner is
constructed by approximating the higher-order finite element operator
with a low-order discretisation on the finer grid which is obtained by
subdividing each grid cell into elements defined by the nodal
points. The resulting sparse system is then solved directly with an
AMG preconditioner. A similar approach is used in \cite{Fischer1997},
where a two-level Schwarz preconditioner is constructed for a spectral
element discretisation of the pressure correction equation in the
Navier Stokes system. Instead of using an AMG method for the sparse
low-order system, the author employs a two-level Schwarz
preconditioner. As in our method, the coarse level is obtained from
the linear finite element discretisation on the original grid and this
is combined with an iterative solution of the low order system on the
finer nodal grid. Finally, instead of re-discretising the equation on
a finer grid, the authors of \cite{Austin2012} construct a sparse
approximation with algebraic methods, which requires additional setup
costs. Compared to our approach, the drawback of all those methods is
that the fine grid problem is obtained from a lowest-order
discretisation (or an algebraic representation leading to a sparse
matrix), and will therefore result in a memory-bound implementation
which can not use the full capacity of modern manycore processors.
\paragraph{Structure}
This paper is organised as follows: in Section \ref{sec:Methods} we
describe the higher order DG discretisation and our hybrid multigrid
solver algorithm which is used throughout this work. We introduce the
principles of the matrix-free smoother implementation and analyse its
computational complexity as a function of the order of the DG
discretisation. The implementation of those algorithms 
based on the DUNE software framework \cite{Bastian2008b,Bastian2014}
is described in Section \ref{sec:Implementation}, where we
also briefly review the principles of sum-factorised operator
application. Numerical results and comparisons between matrix-based
and matrix-free solvers are presented in Section \ref{sec:Results}, which also 
contains an in-depth discussion of the memory requirements of different
implementations. We conclude in Section \ref{sec:Conclusion}. Some more
technical details are relegated to the appendices. In
\ref{sec:lower_order_discretisation} we discuss the re-discretisation
of the operator in the low-order subspace which is used in the coarse
level correction of the multigrid algorithm. A detailed breakdown of
setup costs for the multigrid algorithm can be found in
\ref{sec:setup_costs}.
\section{Methods}\label{sec:Methods}
We consider linear second order elliptic problems of the form
\begin{equation}
  \nabla\cdot\left(\vec{b} u - K \nabla u\right) + cu = f \qquad \text{in} \; \Omega\label{eqn:linearproblem}
\end{equation}
with boundary conditions
\begin{xalignat}{2}
  u &= g \quad\text{on}\; \Gamma_D\subset \partial\Omega,&
  \left(\vec{b}u - K\nabla u\right)\cdot\vec{\nu} &= j \quad\text{on}\; \Gamma_N=\partial\Omega\backslash \Gamma_D
\label{eqn:boundary_conditions}
\end{xalignat}
in the $d$-dimensional domain $\Omega$.
Here the unit outer normal vector on the boundary $\partial\Omega$ is
denoted by \(\vec{\nu} = \vec{\nu}(\vec{x})\). $K=K(\vec{x})\in\mathbb{R}^{d\times d}$ is the diffusion tensor, $\vec{b}=\vec{b}(\vec{x})\in\mathbb{R}^d$ is
the advection vector and $c=c(\vec{x})$ is a scalar reaction coefficient. 
All those terms can depend on the position $\vec{x}\in\Omega$; the fields $K$, $\vec{b}$ and $c$ are not necessarily smooth
and can have large, high contrast jumps. A typical example for a
problem with discontinuous coefficients is the SPE10 benchmark
\cite{Christie2001}. The problem in \eqref{eqn:linearproblem}
could also appear in the linearisation of a non-linear problem such as the 
Navier Stokes equations \cite{2016arXiv161200657P}. In this case it needs to
be solved at every iteration of a Newton- or
Picard- iteration. The reaction term often arises from the implicit
time discretisation of an instationary problem.

\subsection{Higher-order DG discretisation}
To discretise \eqref{eqn:linearproblem} and
\eqref{eqn:boundary_conditions} we use the weighted symmetric interior
penalty (WSIPG) method derived in \cite{Bastian2012} based on
\cite{Wheeler1978,Arnold1982,Arnold2002,Epshteyn2007,Hartmann2008,Ern2009}
together with an upwinding flux for the advective part.  The full
method is described in \cite{Muething2017} and the most important
properties are recalled here. We work on a grid $\mathcal{T}_h$ consisting of
a set of cuboid
cells $T$ with smallest size $h=\min_{T\in\mathcal{T}_h}
\{\operatorname{diam}(T)\}$. Each cell $T$ is the image of a
diffeomorphic map \(\mu_T : \hat T\rightarrow T\) where \(\hat T\) is the
reference cube in \(d\) dimensions. Given $\mathcal{T}_h$, the Discontinuous Galerkin
space $V_h^p\subset L^2(\Omega)$ is the subspace of piecewise polynomials of
degree $p$
\begin{equation}
  V_h^p = \left\{v\in L^2(\Omega) : \forall T\in \mathcal{T}_h, v|_T = q\circ\mu_T^{-1}
  \,\text{with}\, q\in\mathbb{Q}_p^d\right\}.
\end{equation}
Here $\mathbb{Q}_p^d = \{P : P(x_1,\dots,x_d)=\sum_{0\le\alpha_1,\dots,\alpha_d\le p}
c_{\alpha_1,\dots,\alpha_d} x_1^{\alpha_1}\cdot \dots\cdot x_d^{\alpha_d}\}$ is the set of
polynomials of at most degree $p$. It is worth stressing that, in contract to conforming discretisations, functions in $V_h^p$ are not necessarily continuous between cells.
We also define the $L^2$ inner product of two (possibly vector-valued) functions on any domain $Q$ as
\begin{equation}
  (v,w)_Q := \int_Q v\cdot w \;dx.\label{eqn:inner_product}
\end{equation}
For a function $u \in V_h^p$ the discrete weak form of \eqref{eqn:linearproblem} and \eqref{eqn:boundary_conditions} is given by
\begin{equation}
a_h(u,v) = \ell_h(v) \qquad \forall v\in V_h^p
\label{eqn:weak_form}
\end{equation}
with the bilinear form
\begin{equation}
  a_h(u,v) = a_h^{\textup{v}}(u,v) + a_h^{\textup{if}}(u,v) +
  a_h^{\textup{bf}}(u,v)
\label{eqn:bilinear_form}
\end{equation}
which has been split into a volume (v), interface (if) and boundary (b) term defined by
\begin{equation}
\begin{aligned}
  a_h^{\textup{v}}(u,v) &= \sum_{T\in\mathcal{T}_h} \big[\left(K\nabla u - \vec{b}u,\nabla v\right)_T +  \left(c u, v\right)_T\big] \\
  a_h^{\textup{if}}(u,v) &= \sum_{F\in\mathcal{F}_h^{i}}
  \big[\left( \Phi(u^-,u^+,\vec{\nu}_F\cdot \vec{b}),\llbracket v \rrbracket\right)_F
    - \left(\vec{\nu}_F\cdot \{K\nabla u\}_\omega,\llbracket v \rrbracket\right)_F-
  \left(\vec{\nu}_F\cdot\{K \nabla v\}_\omega,\llbracket u \rrbracket\right)_F
 +
 \gamma_{F}\left(\llbracket u \rrbracket,\llbracket v \rrbracket  \right)_F\big]\\
  a_h^{\textup{bf}}(u,v) &=
  \sum_{F\in\mathcal{F}_h^{D}} \big[\left(
  \Phi(u,0,\vec{\nu}_F\cdot \vec{b}),v\right)_F
 -\left(\vec{\nu}_F\cdot (K\nabla u),v\right)_F
  -\left(\vec{\nu}_F\cdot K \nabla v,u\right)_F
  + \gamma_{F} \left(u,v\right)_F\big] .
\end{aligned}
\label{eqn:bilinear_form_split}
\end{equation}
The right hand side functional is
\begin{equation}
\begin{split}
\ell_h(v) &=  \sum_{T\in\mathcal{T}_h} \left( f,v\right)_T
- \sum_{F\in\mathcal{F}_h^{N}} \left(j,v\right)_F
- \sum_{F\in\mathcal{F}_h^{D}} \big[\left(\Phi(0,g,\vec{\nu}_F\cdot \vec{b}), v\right)_F
+ \left(\vec{\nu}_F\cdot (K \nabla v),g\right)_F
-\gamma_{F} \left(g,v \right)_F\big].
\end{split}
\label{eqn:rhs}
\end{equation}
In those expressions the set of interior faces is denoted as
$\mathcal{F}_h^i$. The sets of faces on the Neumann- and Dirichlet boundary are
$\mathcal{F}_h^N$ and $\mathcal{F}_h^D$ respectively. With each face
$F\in\mathcal{F}_h=\mathcal{F}_h^i\cup\mathcal{F}_h^N \cup \mathcal{F}_h^D$ we
associate an oriented unit normal vector $\vec{\nu}_F$. For any point $\vec{x}\in F \in
\mathcal{F}_h^i$ on an internal face we define the jump
\begin{equation}
\llbracket v \rrbracket (\vec{x}) = v^-(\vec{x})-v^+(\vec{x})
\end{equation}
and the weighted average
\begin{equation}
\{ v \}_\omega (\vec{x}) = \omega^-(\vec{x}) v^-(\vec{x}) + \omega^+(\vec{x}) v^+(\vec{x})\quad
\text{for weights}\quad \omega^-(\vec{x}) + \omega^+(\vec{x}) = 1, \omega^\pm(\vec{x}) \geq 0
\end{equation}
where $v^{\pm}=\lim_{\epsilon\downarrow 0}v(\vec{x}\pm \epsilon\vec{\nu}_F)$ is the value of the field
on one side of the facet.  In the expressions $\ell_h(v)$ and
$a_h^{\textup{bf}}(u,v)$ we implicitly assume that $u$ and $v$
are evaluated on the interior side of the boundary facets. To account
for strongly varying diffusion fields, we choose the weights
introduced in \cite{Ern2009}:
\begin{equation}
\omega^{\pm}(\vec{x}) = \frac{\delta_{K\nu}^{\mp}(\vec{x})}{\delta_{K\nu}^-(\vec{x}) +
  \delta_{K\nu}^+(\vec{x})},\quad \text{with}\quad \delta_{K\nu}^{\pm}(\vec{x})
=\nu^T(\vec{x})K^\pm(\vec{x}) \nu(\vec{x}).
\end{equation}
Defining the harmonic average of two numbers as $\langle a,b\rangle=2ab/(a+b)$,
the penalty factor $\gamma_{F}$ is chosen based on combination of choices
in \cite{Epshteyn2007,Hartmann2008,Ern2009} as in \cite{Bastian2012}
\begin{equation}
\gamma_{F} = \alpha \cdot p (p+d-1)
\cdot \begin{cases}  \frac{\langle \delta_{K\nu_F}^- , \delta_{K\nu_F}^+ \rangle|F|}{\min(|T^-(F)|,|T^+(F)|)} & \text{for $F\in\mathcal{F}_h^i$}\\[1ex]
  \frac{\delta_{K\nu_F}^-|F|}{|T^-(F)|} & \text{for $F\in\mathcal{F}_h^D$}.\\
\end{cases}
\end{equation}
In this expression $\alpha$ is a free parameter chosen to be $\alpha=1.25$ in
all the computations reported below.  To discretise the advection
terms, we use the upwind flux on the face $F$ which is given by
\begin{equation*}
\Phi(u^-,u^+,b_{\vec{\nu}}) = \left\{\begin{array}{ll}
b_{\vec{\nu}} u^- & \text{if}\quad b_{\vec{\nu}}\geq 0\\
b_{\vec{\nu}} u^+ & \text{otherwise}
\end{array}\right. .
\end{equation*}
While the formulation in \eqref{eqn:bilinear_form}, \eqref{eqn:rhs} is
valid on any grid, here we only consider grids based on hexahedral
elements. To use sum-factorisation techniques, we also assume that the
basis functions on the reference element as well as the quadrature rules
have tensor-product structure. Similar techniques
can be used on simplicial elements \cite{Kirby2011,Kirby2012}. A basis
$\Psi_h = \{\psi_1^p, \hdots,\psi_{N_h}^p\}$ is chosen for $V_h^p$, i.e. every
function $u\in V_h^p$ can be written as
\begin{equation}
  u(x) = \sum_{i\in I} u_i \psi_i^p(x),\qquad I = \{1,\dots,N_h\}\subset \mathbb{N}.
\end{equation}
Here $N_h$ is the total number of degrees of freedom and $\vecn{u}=(u_1,u_2,\dots,u_{N_h})^T\in \mathbb{R}^n$ is the vector of unknowns. With this basis the weak formulation in
\eqref{eqn:weak_form} is equivalent to a matrix equation
\begin{equation}
  A \vecn{u} = \vecn{f}\qquad \text{with}\quad A_{ij}=a_h(\psi^p_j,\psi_i^p)\quad\text{and}\quad \vecn{f}\in\mathbb{R}^n\;\;\text{with}\;\;f_i := \ell_h(\psi_i^p).\label{eqn:linear_equation}
\end{equation}
It is this equation which we will solve in the following. Since we
consider a discontinuous function space, the basis can be chosen such
that the support of any basis function $\psi_i^p\in \Psi_h$ is restricted to a
single element $T\in \mathcal{T}_h$ and this implies that the matrix $A$
is block-structured. This property is important for the algorithms
developed in the following. Throughout this work we use a basis which is
constructed from the tensor-product of one-dimensional Gauss-Lobatto
basis functions (i.e. Lagrange polynomials with nodes at the Gauss-Lobatto points) on each grid cell. We also assume that the basis is nodal,
i.e. there is a set of points $\vec{\zeta}_j$, $j\in I$ such that
$\psi_i^p(\vec{\zeta}_j)=\delta_{ij}$.
\paragraph{Block notation of linear systems}
To efficiently deal with block-structured matrices we introduce the
following notation. For any finite index set $I\subset\mathbb{N}$ we define
the vector space $\mathbb{R}^I$ to be isomorphic to $\mathbb{R}^{|I|}$
with components indexed by $i\in I$. Thus any vector of unknowns $\vecn{u} \in
\mathbb{R}^I$ can be interpreted as a mapping $\vecn{u} : I \to \mathbb{R}$ and
$\vecn{u}(i) = u_i$. In the same way, for any two finite index sets $I,
J\subset\mathbb{N}$ we write $A \in \mathbb{R}^{I\times J}$ with the interpretation
$A : I \times J \to \mathbb{R}$ and $A(i,j) = A_{i,j}$. Finally, for any
subset $I'\subseteq I$ we define the restriction matrix $R_{I,I'} :
\mathbb{R}^I \to \mathbb{R}^{I'}$ as $ (R_{I,I'} \vec{u})_i = u_i \ \forall i\in I'$.
Now, let $I_T = \{ i\in\mathbb{N} : \textup{supp} \;\psi_i^p \subseteq T
\}$ denote the subset of indices which are associated with basis
functions that are nonzero on element $T$. Then
\begin{align}
\bigcup_{T\in \mathcal{T}_h} I_T &= I = \{1,\ldots,n\} \quad \text{and} \quad
I_T \cap I_{T'} = \emptyset, \ \forall T\neq T',
\end{align}
where the last property follows from the fact that we consider a
discontinuous function space and therefore $\{I_T\}$ forms a disjoint
partitioning of the index set $I$. Using this partitioning and
imposing an ordering $\mathcal{T}_h=\{T_1,\ldots,T_m\}$ on the mesh
elements, the linear system in \eqref{eqn:linear_equation} can be
written in block form
\begin{equation}
  \left(
    \begin{array}{ccc}
      A_{T_1,T_1} & \hdots & A_{T_1, T_m} \\
      \vdots & \ddots & \vdots\\
      A_{T_m,T_1} & \hdots & A_{T_m,T_m}
    \end{array}
  \right)
  \left(
    \begin{array}{c}
      \vec{u}_{T_1}\\\vdots\\ \vec{u}_{T_m}
    \end{array}
  \right) =
  \left(
    \begin{array}{c}
      \vec{f}_{T_1} \\ \vdots \\ \vec{f}_{T_m}
    \end{array}
\right)
\label{eq:block_ls}
\end{equation}
where $A_{T_\rho,T_\sigma} = R_{I,I_{T_\rho}} A R_{I,I_{T_\sigma}}^T$,
$\vec{u}_{T_\rho} = R_{I,I_{T_\rho}} \vec{u}$ and $\vec{f}_{T_\rho} = R_{I,I_{T_\rho}} \vec{f}$.
For each cell $T$ and function $u\in V_h^p$ we also define the function $u_{\cap T}\in V_h^p$ as
\begin{equation}
  u_{\cap T}(x)
=  \begin{cases}
  u(x) & \text{for $x \in T$}\\
  0 & \text{otherwise}
\end{cases}
\qquad
\Leftrightarrow
\qquad
u_{\cap T} = \sum_{i\in I_{T}} \left(u_{T}\right)_i \psi_i^p = \sum_{i\in I_{T}} u_i \psi_i^p.\label{eqn:ucupT}
\end{equation}
Our choice of Gauss-Lobatto basis functions with non-overlapping
support leads to the so called \emph{minimal stencil}; the matrices
$A_{T_\rho,T_\sigma}$ are only nonzero if $\rho=\sigma$ or if the elements $T_\rho$ and
$T_\sigma$ share a face, and hence the matrix in \eqref{eq:block_ls} is
block-sparse. In general and for the basis considered in this paper,
the $(p+1)^d\times (p+1)^d$ matrix $A_{T_\rho,T_\sigma}$ is dense.

This partitioning of \(A\) and the minimal stencil property form are key
for the matrix-free preconditioner implementation discussed below.
\subsection{Hybrid multigrid preconditioner}
\label{sec:hybr-mult-prec}
To solve the linear system in \eqref{eqn:linear_equation} we use the
hybrid multigrid preconditioner described in
\cite{Bastian2012}. However, in contrast to \cite{Bastian2012}, the
computationally most expensive components, namely the operator- and
preconditioner application in the DG space $V_h^p$, are implemented in
a matrix free way. As the numerical results in Section
\ref{sec:Results} confirm, this leads to significant speedups at
higher polynomial orders $p$ and has benefits even for relatively
modest orders of $p=2,3$.

The key idea in \cite{Bastian2012} is to combine a block-smoother in
the high-order DG space $V_h^p$ with an AMG correction in a low-order
subspace $\hat{V}_h\subset V_h$. In this work we will consider low order
subspaces $\hat{V}_h$ spanned by the piecewise constant $P_0$ and
conforming piecewise linear elements $Q_1$. One V-cycle of the
multigrid algorithm is shown in Algorithm
\ref{alg:hybrid_multigrid}. At high polynomial orders $p$ the
computationally most expensive steps are the operator application in
the defect (or residual) calculation and the pre-/post-smoothing on
the DG level. Usually the multigrid algorithm is used as a
preconditioner of a Krylov subspace method, which requires additional
operator applications.
\begin{algorithm}
\caption{Hybrid multigrid V-cycle. Input: Right hand side $\vecn{f}$, initial guess $\vecn{u}_0$ (often set to zero); Output:
  $\vecn{u}$, the approximate solution to $A\vecn{u}=\vecn{f}$}
\label{alg:hybrid_multigrid}
\begin{center}
\begin{algorithmic}[1]
\STATE{$\vecn{u}\mapsfrom S(f,\vecn{u}_0;n_{\operatorname{pre}},\omega)$}\qquad\COMMENT{\textit{Presmooth with $n_{\operatorname{pre}}$ iterations in DG space}}
\STATE{$\vecn{r}\mapsfrom \vecn{f}-A\vecn{u}$}\qquad\COMMENT{\textit{Calculate residual}}
\STATE{$\hat{\vecn{r}} \mapsfrom R(\vecn{r})$}\qquad\COMMENT{\textit{Restrict to subspace}}
\STATE{$\hat{\vecn{u}}\mapsfrom \operatorname{AMG}(\hat{\vecn{r}})\approx \hat{A}^{-1}\hat{\vecn{r}}$}\qquad\COMMENT{\textit{AMG V-cycle on subspace fields}}
\STATE{$\vecn{u}\mapsfrom \vecn{u} + P(\hat{\vecn{u}})$}\qquad\COMMENT{\textit{Prolongate and add to DG solution}}
\STATE{$\vecn{u}\mapsfrom S(f,\vecn{u};n_{\operatorname{post}},\omega)$}\qquad\COMMENT{\textit{Postsmooth with $n_{\operatorname{post}}$ iterations in DG space}}
\end{algorithmic}
\end{center}
\end{algorithm}
To solve the problem in the lower-order subspace we use the
aggregation based AMG solver described in
\cite{Blatt2010,Blatt2012}. In the following we discuss the efficient
implementation of the individual components in Algorithm
\ref{alg:hybrid_multigrid}:
\begin{enumerate}
\item Matrix-free operator application in the DG space to calculate the residual (line 2)
\item Matrix-free pre-/post-smoothing in the DG space (lines 1 and 6)
\item Transfer of grid functions between between DG space and low order subspace (lines 3 and 5)
\item Assembly of the low order subspace matrix which is used in the AMG V-cycle (line 4)
\end{enumerate}
The matrix-free operator application based on sum factorisation techniques has been previously described in \cite{Muething2017} (and is briefly reviewed in Section \ref{sec:operator_application} below). Here we will concentrate on the other components.
\subsubsection{Matrix-free smoothers}
\label{sec:matr-free-prec}
In the context of multigrid methods, iterative improvement of the
solution with a stationary method is referred to as \textit{smoothing}
since a small number of iterations with a simple method reduces
high-frequency error components.  For a linear equation $A\vecn{u}=\vecn{f}$, this
process can be summarised as follows: given a matrix \(W\approx A\),
calculate the defect \(\vecn{d} = \vecn{f} - A \vecn{u}^{(k-1)}\) for some guess
$\vecn{u}^{(k-1)}$ and solve for \(\vecn{\delta u}\) such that
\begin{equation}
  W \vecn{\delta u} = \vecn{d}.
\end{equation}
This correction $\vecn{\delta u}$ can then be used to obtain an improved solution
\begin{equation}
  \vecn{u}^{(k)} = \vecn{u}^{(k-1)}+\vecn{\delta u} = \vecn{u}^{(k-1)} + W^{-1} (\vecn{f}-A\vecn{u}^{(k-1)}).
\end{equation}
Starting from some initial guess $\vecn{u}^{(0)}$, the iterates
$\vecn{u}^{(0)},\vecn{u}^{(1)},\dots,\vecn{u}^{(k)}$ form an increasingly better
approximation to the exact solution $\vecn{u}$ of $A\vecn{u}=\vecn{f}$ for a convergent method. Choosing a
block-smoother guarantees that any high-frequency oscillations
inside one cell are treated exactly, which is particularly important
for higher polynomial orders. Lower-frequency errors with variations
between the grid cells are reduced by the AMG preconditioner in the
low-order subspace (line 4 in Alg. \ref{alg:hybrid_multigrid}).  Following the
block partitioning in
\eqref{eq:block_ls}, let \(A = D + L + U\) be decomposed into the
strictly lower block triangular part \(L\), the block diagonal \(D\)
and the strictly upper block triangular part \(U\) with
\begin{xalignat}{3}
  D_{T_\rho,T_\sigma}&=
  \begin{cases}
    D_{T_\rho} = A_{T_\rho,T_\rho} & \text{for $\rho=\sigma$}\\
    0 & \text{otherwise,}
  \end{cases}
&
  L_{T_\rho,T_\sigma}&=
  \begin{cases}
    A_{T_\rho,T_\sigma} & \text{for $\rho>\sigma$}\\
    0 & \text{otherwise,}
  \end{cases}
&
  U_{T_\rho,T_\sigma}&=
  \begin{cases}
    A_{T_\rho,T_\sigma} & \text{for $\rho<\sigma$}\\
    0 & \text{otherwise.}
  \end{cases}
\end{xalignat}
In this work we consider the following choices for the block preconditioner
$W$: (1) block-Jacobi $W=\omega^{-1} D$, (2) successive block
over-relaxation (block-SOR) $W=\omega^{-1}D+L$ and (3) symmetric block-SOR
(block-SSOR) $W=(\omega(2-\omega))^{-1}(D+\omega L)D^{-1}(D+\omega U)$. The over-relaxation
parameter $\omega$ can be tuned to accelerate convergence. The explicit
form of the block-Jacobi and block-SOR algorithms are written down in
Algorithms \ref{alg:block_jacobi} and \ref{alg:block_SOR}. Note that
for block-SOR, the field $\vecn{u}^{(k)}$ can be obtained by transforming
$\vecn{u}^{(k-1)}$ in-place, which reduces storage requirements. The
corresponding backward-SOR
sweep is obtained by reversing the direction of the grid traversal. It
can then be shown that one application of the block-SSOR smoother is
equivalent to a forward-block-SOR sweep following by a
backward-block-SOR iteration.
\begin{algorithm}
\caption{Block-Jacobi smoother. Input: Right hand side $\vecn{f}$, initial
  guess $\vecn{u}^{(0)}$, number of smoothing steps
  $n_{\operatorname{smooth}}$. Output: Approximate solution
  $\vecn{u}^{(n_{\text{smooth}})}$}
\label{alg:block_jacobi}
\begin{center}
\begin{algorithmic}[1]
\FOR{$k=1,\dots,n_{\operatorname{smooth}}$}
\STATE{Calculate defect $\vecn{d} = \vecn{f} - A\vecn{u}^{(k-1)}$}
\STATE{Solve $D_{T} \vecn{\delta u}_{T}=\vecn{d}_{T}$ for $\vecn{\delta u}_{T}$ on all grid elements $T$}
\STATE{Update solution $\vecn{u}^{(k)} = \vecn{u}^{(k-1)} + \omega\vecn{\delta u}$}
\ENDFOR
\end{algorithmic}
\end{center}
\end{algorithm}
\begin{algorithm}
\caption{Successive block-over-relaxation, forward Block-SOR
  smoother. Input: Right hand side $f$, initial guess $\vecn{u}^{(0)}$,
  number of smoothing steps $n_{\operatorname{smooth}}$. Output:
  Approximate solution $\vecn{u}^{(n_{\text{smooth}})}$}
\label{alg:block_SOR}
\begin{center}
\begin{algorithmic}[1]
\FOR{$k=1,\dots,n_{\operatorname{smooth}}$}
\FOR{all grid elements $T_\rho$ with $\rho=1,\dots,m$}
\STATE{Calculate local defect\begin{equation}\vecn{d}_{T_\rho} = \vecn{f}_{T_\rho} - \sum_{\substack{\sigma<\rho\\\partial T_{\sigma}\cap \partial T_\rho\ne\emptyset}} A_{T_\rho,T_{\sigma}} \vecn{u}^{(k)}_{T_\sigma}-\sum_{\substack{\sigma>\rho\\\partial T_{\sigma}\cap \partial T_\rho\ne\emptyset}} A_{T_\rho,T_{\sigma}} \vecn{u}^{(k-1)}_{T_{\sigma}}\end{equation}}
\STATE{Solve $D_{T_\rho} \vecn{\delta u}_{T_\rho}=\vecn{d}_{T_\rho}$ for $\vecn{\delta u}_{T_\rho}$}
\STATE{Update solution $\vecn{u}^{(k)}_{T_\rho} = (1-\omega) \vecn{u}^{(k-1)}_{T_\rho} + \omega\vecn{\delta u}_{T_\rho}$}
\ENDFOR
\ENDFOR
\end{algorithmic}
\end{center}
\end{algorithm}
Note that in all cases the computationally most expensive operations
are the multiplication with the dense $(p+1)^d\times (p+1)^d$ matrices
$A_{T_\rho,T_\sigma}$ and the solution of the dense system $D_{T_\rho}\vecn{\delta u}_{T_\rho}=\vecn{d}_{T_\rho}$ in each cell. As will be described in Section
\ref{sec:Implementation}, the defect can be calculated in a
matrix-free way by using the techniques from \cite{Muething2017}. The
standard approach to calculating the correction $\vecn{\delta u}_{T_\rho}$ is to
assemble the local matrices $D_{T_\rho}$, factorise them with an
LU-/Cholesky-approach and then use this to solve the linear
system. However, the costs for the factorisation- and back-substitution
steps grow like $(p+1)^{3d}$ and $(p+1)^{2d}$ respectively, and hence
this method becomes prohibitively expensive for high polynomial
degrees $p$. To avoid this, we solve the block-system
\begin{equation}
D_{T_\rho}\vecn{\delta u}_{T_\rho}=\vecn{d}_{T_\rho}\label{eqn:block_system}
\end{equation}
with a Krylov-subspace method such as the Conjugate
Gradient (CG) \cite{Hestenes1952} or Generalised Minimal Residual (GMRES) \cite{Saad1986}
iteration. This only requires operator applications which can be
implemented in a matrix-free way. If sum-factorisation techniques are
used, the overall cost of this inversion is $\order(n_{\operatorname{iter}}(\epsilon)\cdot d p^{d+1})$
where $n_{\operatorname{iter}}(\epsilon)$ is the number of iterations required
to solve the block-system.  An important observation, which
drastically effects the practical usage of the matrix-free smoothers,
is that it suffices to compute an inexact solution of
\eqref{eqn:block_system} with the Krylov
subspace method. For this we iterate until the two-norm of the
residual $||\vecn{d}_{T_\rho}-D_{T_\rho}\vecn{\delta u}_{T_\rho}||_2$ has been reduced by a fixed
factor $\epsilon\ll 1$. The resulting reduction in $n_{\operatorname{iter}}$
and preconditioner cost has to be balanced against the potential
degradation in the efficiency of the preconditioner. This dependence on
$\epsilon$ on the overall solution time will be studied below.

Using block smoothers, the method proposed here is not strictly robust with respect to
the polynomial degree $p$. This can be remedied by using overlapping block smoothers
\cite{Blatt2010,Bastian2012} where the blocks originate from an overlapping partitioning
of the mesh. However, this also increases the computational
cost substantially. Since in practical applications we consider moderate polynomial
degree anyway, we concentrate on simple non-overlapping block smoothers here
where the blocks originate from a single element.
\subsubsection{Intergrid operators}\label{sec:intergrid_operators}
Since we use nested function spaces with $\hat{V}_h\subset V_h^p$, injection
is the natural choice for the prolongation $P:\hat{V}_h\rightarrow V_h^p$.
Given a function $\hat{u}\in
\hat{V}_h$, we define
\begin{equation}
  P:\hat{u}\mapsto u\in V_h^p \quad\text{such that $u(\vec{x})=\hat{u}(\vec{x})$ for all $\vec{x}\in\Omega$}.
\end{equation}
Expanding in basis functions on element $T$ and using \eqref{eqn:ucupT} this
leads to
\begin{equation}
  \hat{u}_{\cap T} = \sum_{i\in \hat{I}_T} \hat{u}_i \hat{\psi}_i = \sum_{j\in I_T} u_j \psi^p_j = u_{\cap T}.
\end{equation}
Here $\hat{I}_T$ is the set of indices $i$ for which the coarse basis
function $\hat{\psi}_i(x)$ is non-zero on element $T$. For a general
basis, calculation of $u_j$ requires solution of a small system
involving the local $(p+1)^d\times (p+1)^d$ mass matrix $\int_T
\psi_i^p(\vec{x})\psi_j^p(\vec{x})\;dx$ in each cell $T$; this is expensive for higher
polynomial degrees $p$.  However, here we only consider the special
case that $\hat{V}_h$ is a subspace of $V_h^p$ and the basis is nodal.
In this case the unknowns in the higher order space can be found much
more efficiently by evaluating the low order basis functions at the
nodal points \(\vec{\zeta}_i\) of the higher order space,
\begin{equation}
  u_i = \sum_{j\in\hat{I}_T} \hat{u}_{j} \hat{\psi}_j(\vec{\zeta}_i) \qquad\text{for all $i\in I_T$}.
\end{equation}
In other words, the \(j\)-th column of the local prolongation matrix
$P|_T$ can be calculated by evaluating \(\hat{\psi}_j\) at all nodal
points,
\begin{equation}
  (P|_T)_{ij} = \hat{\psi}_j(\vec{\zeta}_i) \qquad\text{for $i\in I_T$ and $j\in \hat{I}_T$}.
\end{equation}
Since the DG basis functions form a partition of unity, for the $P_0$ subspace this simplifies even further to
\begin{equation}
  P|_T = (1,\hdots,1)^T.
\end{equation}
For the restriction operation $R:V_h^p\rightarrow\hat{V}_h$ we choose the
transpose of the prolongation, i.e. $R=P^T$. In the absence of
advection ($\vec{b}=0$) this choice preserves the symmetry of the system and
the WSIPG discretisation.
\subsubsection{Low order subspace matrix}
\label{sec:low-order-subspace}
The low order subspace matrix \(\hat{A}\) is needed to solve the
``coarse'' grid correction with the AMG preconditioner. In general
this matrix is given by the Galerkin product
\begin{equation}
  \hat{A} = R A P = P^T A P.\label{eqn:galerkin_product}
\end{equation}
However, calculating $\hat{A}$ according to this formula is expensive
for several reasons: first of all it requires two matrix-matrix
multiplications of large matrices with different sparsity patterns. In
addition, the matrix $A$ in the DG space has to be assembled, even
though in a matrix-free implementation it is never used in the DG
smoother. Furthermore, for our choice of coarse function spaces, the
matrix $\hat{A}$ has lots of entries which are formally zero. However,
if it is calculated according to \eqref{eqn:galerkin_product} in
an actual implementation those entries will not vanish exactly due to
rounding errors. This leads to drastically increased storage
requirements and unnecessary calculations whenever $\hat{A}$ is used
in the AMG preconditioner. All those issues can be avoided by directly
assembling the coarse grid matrix $\hat{A}$ as follows:

For the $P_0$ subspace it is shown in \ref{sec:low_order_matrix_P0}
that on an equidistant grid the matrix entries $\hat{A}_{ij}$ can be
calculated directly by discretising \eqref{eqn:linearproblem} with
a modified finite volume scheme in which (i) all flux terms are
multiplied by a factor $\alpha p(p+d-1)$ and (ii) the boundary flux is
scaled by an additional factor 2. For the conforming $Q_1$ subspace,
all contributions to $\hat{A}$ which arise from jump terms in
$A_{ij}=a_h(\psi^p_j,\psi^p_i)$ vanish. In fact, up to an additional
Dirichlet boundary term $\sum_{F\in \mathcal{F}_h^D}\gamma_F
\left(\hat{\psi}_i,\hat{\psi}_j\right)_F$, the matrix $\hat{A}_{ij}$ is identical to the
matrix obtained from a conforming $Q_1$ discretisation of
\eqref{eqn:linearproblem}, and it is this matrix that we use in
our implementation. In both cases the sparsity pattern of $\hat{A}$ is
the same as that obtained by a discretisation in the low-order
subspace. In summary this re-discretisation in the coarse level
subspace instead of using the Galerkin product in
\eqref{eqn:galerkin_product} results in a significantly reduced
sparsity pattern of the matrix $\hat{A}$ and avoids explicit assembly
of the large DG matrix $A$.
\section{Implementation in EXADUNE}\label{sec:Implementation}
The crucial ingredient for the development of the matrix-free block
smoothers described in Section \ref{sec:matr-free-prec} is the
efficient on-the-fly operator application of the full operator $A$ and
multiplication with the block-diagonal \(D_{T}\) for each element
\(T\in\mathcal{T}_h\). While the operator application is required in the
calculation of the local defect, the frequent multiplication with
$D_{T}$ is necessary in the iterative inversion of the diagonal
blocks. Our code is implemented in the EXADUNE \cite{Bastian2014} fork
of the DUNE framework \cite{Bastian2008a,Bastian2008b}. EXADUNE
contains several optimisations for adapting DUNE to modern computer
architectures.
\subsection{Matrix-free operator application}\label{sec:operator_application}
The matrix-free application of the operator $A$ itself has been
implemented and heavily optimised in \cite{Muething2017}, and here we
only summarise the most important ideas. Recall that the weak form
$a_h(u,v)$ in \eqref{eqn:bilinear_form} is split into volume-,
interface- and boundary terms defined in
\eqref{eqn:bilinear_form_split}. The local contributions to the system
matrix $A$ are evaluated in a \texttt{LocalOperator} class in
DUNE. More specifically, let $T^+$, $T^-$ be a pair of adjacent cells which
share a in internal face $F^i=\partial T^+(F^i)\cap \partial T^-(F^i)$ and
let $F^b\in \partial \tilde{T}(F^b) \cap \mathcal{F}_h^b$ be the face of a
cell $T$ adjacent to the boundary of the domain. For each cell $T$, internal
interface $F^i$ and boundary face $F^b$ the \texttt{LocalOperator} class
provides methods for evaluating the \textit{local} operator application
\begin{equation}
  \begin{aligned}
    \left(A^\textup{v}_{T,T}\vecn{u}_{T}\right)_{i_T} &:= a_h^{\textup{v}}(\vecn{u}_{\cap T},\psi_{i_T}^p)
\qquad\text{with}\quad  i_T \in I_{T}&\text{(volume)}\\
    \left(A^{\textup{if}}_{F^i;T',T''}\vecn{u}_{T''}\right)_{i_{T'}} &:= a_h^{\textup{if}}(\vecn{u}_{\cap T''},\psi_{i_{T'}}^p)
\qquad\text{with}\quad  i_{T'} \in I_{T'}&\text{(interface)}\\
    \left(A^{\textup{b}}_{F^b;\tilde{T}}\vecn{u}_{T}\right)_{i_{\tilde{T}}} &:= a_h^{b}(\vecn{u}_{\cap T},\psi_{i_{\tilde{T}}}^p)
\qquad\text{with}\quad  i_{\tilde{T}} \in I_{\tilde{T}}&\text{(boundary).}\\
  \end{aligned}
\end{equation}
For the interface term we calculate four contributions,
one for each of the combinations $(T',T'') = (T^-,T^-)$, $(T^-,T^+)$,
$(T^+,T^-)$ and $(T^+,T^+)$. 
By iterating over all elements and faces
of the grid and calling the corresponding local operator evaluations,
the full operator can be applied to any field $u\in V_h^p$.
This iteration over the grid is carried out by an instance of a
\texttt{GridOperator} class, which takes the \texttt{LocalOperator} as
template parameter.

In the following we briefly recapitulate the key ideas used to
optimise the matrix-free operator application with sum-factorisation
techniques and concentrate on application of the operator $A_{T,T}^{\text{v}}$.
In each cell $T$ this proceeds in three stages (see Algorithm 1 in \cite{Muething2017} for details)
\begin{enumerate}
\item Given the coefficient vector $\vecn{u}$, calculate the values of the field $u$ and its gradient at the quadrature points on the reference element $\hat{T}$.
\item At each quadrature point, evaluate the coefficients $K$, $\vec{b}$ and $c$ as well a geometric factors from the mapping between $T$ and $\hat{T}$. Using this and the function values from Stage 1, evaluate the weak form $a_h(\cdot,\cdot)$ at each quadrature point.
  \item Multiply by test functions (and their derivatives) to recover the entries of the coefficient vector $\vecn{v}=A\vecn{u}$.
\end{enumerate}
We now discuss the computational complexity of the individual stages.
\paragraph{Stage 1}
We assume that the basis functions in the reference element $\hat{T}$ have a tensor-product structure
\begin{equation}
  \hat{\phi}_{\vec{j}}(\hat{\vec{x}}) = \hat{\phi}_{(j_1,\dots,j_d)}(\hat{x}_1,\dots,\hat{x}_d) = \prod_{k=1}^d\hat{\theta}^{(k)}_{j_k}(\hat{x}_k).
\end{equation}
Those basis functions are enumerated by $d$-tuples $\vec{j}=(j_1,\dots,j_d)\in \vec{J}=J^{(1)}\times\dots\times J^{(d)}$, $J^{(k)}=\{1,\dots,n_k\}$.
For each cell $T$ we have an index map $g(T,\cdot)\rightarrow I$ which assigns a global index to each $d$-tuple $\vec{j}$ in the cell; this implies that in a particular cell $T$ any function $u\in V_h^p$ can be written as
\begin{equation}
  u(\vec{x}) = \hat{u}(\hat{\vec{x}}) = \sum_{\vec{j}\in\vec{J}} u_{g(T,\vec{j})}\hat{\phi}_{\vec{j}}(\hat{\vec{x}}) \qquad\text{for $\vec{x}=\mu_T(\vec{\hat{x}})\in T$.}\label{eqn:value_at_point}
\end{equation}
We use a tensor-product quadrature rule with points
\begin{equation}
  \vec{\xi}_{\vec{i}} = \left(\xi_{i_1}^{(1)},\dots,\xi_{i_d}^{(d)}\right)
\end{equation}
enumerated by $d$-tuples $\vec{i} = (i_1,\dots,i_d) \in I^{(1)}\times\dots\times I^{(d)}$, $I^{(k)} = \{1,\dots,m_k\}$. In the following we only consider $n_1=\dots=n_d=:n$ and $m_1=\dots=m_d=:m$, in this case that there are
$m^d$ quadrature points and $n^d$ basis functions per cell.
Following Eq. (13) in \cite{Muething2017}, the evaluation of the function $u$ at the quadrature point $\vec{\xi}_{\vec{i}}$ can be written as
\begin{equation}
  \hat{u}(\vec{\xi}_{\vec{i}}) = \sum_{\vec{j}\in\vec{J}} u_{g(T,\vec{j})}\hat{\phi}_{\vec{j}}(\vec{\xi}_{\vec{i}}) = \sum_{j_1=1}^n\dots\sum_{j_d=1}^n A_{i_1,j_1}^{(1)}\cdot\dots\cdot A_{i_d,j_d}^{(d)} u_{g(T,(j_1,\dots,j_d))}.\label{eqn:naive_evaluation}
\end{equation}
The entries of the $m\times n$ matrices $A^{(q)}$ are given by the one-dimensional shape functions evaluated at the quadrature points, $A_{\alpha,\beta}^{(k)} = \theta_\beta^{(k)}(\xi_\alpha^{(k)})$, $\alpha\in I^{(k)}$, $\beta\in J^{(k)}$.
\eqref{eqn:naive_evaluation} describes the multiplication of a vector of length $n^d$ by a $m^d\times n^d$ matrix $\hat{\Phi}$ with $\hat{\Phi}_{\vec{i},\vec{j}} := \hat{\phi}_{\vec{j}}(\vec{\xi}_{\vec{i}})=A_{i_1,j_1}^{(1)}\cdot\dots\cdot A_{i_d,j_d}^{(d)}$. Assuming that the ratio $\rho=n/m\le 1$ is constant, the naive cost of this operation is $\mathcal{O}((d+1)m^{2d})$ FLOPs.
A significant reduction in computational complexity is achieved by applying the small $m\times n$ matrices $A^{(k)}$ recursively in $d$ stages as
\begin{equation}
  \begin{aligned}
    u^{(0)}_{(j_1,\dots,j_d)} &= u_{g(T,(j_1,\dots,j_d))}\\
    u^{(s)}_{(i_1,\dots,i_s,j_{s+1},\dots,j_d)} &= \sum_{j_s=1}^n A_{i_s,j_s}^{(s)} u^{(s-1)}_{(i_1,\dots,i_{s-1},j_s,\dots,j_d)} \qquad\text{for $s=1,\dots,d$}\\
    \hat{u}(\vec{\xi}_{\vec{i}}) &= u^{(d)}_{(i_1,\dots,i_d)}.
  \end{aligned}
  \label{eqn:sumfact_evaluation}
\end{equation}
Each step in \eqref{eqn:sumfact_evaluation} requires the multiplication with an $m\times n$ matrix $A_{i_s,j_s}^{(s)}$ for $m^{s-1} n^{d-s}\le m^{d-1}$ vectors. Multiplication by an $m\times n$ matrix requires $2n\times m$ operations and hence the total cost is reduced to \mbox{$\mathcal{O}(d\cdot m^{d+1})$}. Since $m$ grows with the polynomial degree (typically $m=n=p+1$), for high polynomial degrees $p$ this sum-factorised approach is much faster than the naive algorithm with complexity \mbox{$\mathcal{O}((d+1)m^{2d})$}. Gradients of $u$ can be evaluated in the same way with slightly different matrices $A^{(s)}$.

The face integrals in the bilinear form require the evaluation on points of a quadrature rule of lower dimension $d-1$. Although this is complicated by the different embeddings of a face into the neighbouring cells, the same sum factorisation techniques can be applied. While for those terms the cost is only $\mathcal{O}(d\cdot m^d)$, for low polynomial degrees the absolute cost can still be larger than the cost of the volume terms (see Fig. 1 in \cite{Muething2017}). We conclude that - using sum factorisation - the cost of the first stage in the matrix-free operator application is $\mathcal{O}(d\cdot m^{d+1})$.

\paragraph{Stage 2} Since the operations can be carried out independently at the $m^d$ quadrature points, the cost of this stage is $\order(m^d)$.
\paragraph{Stage 3} Given the values $z_{(i_1,\dots,i_d)}:=\hat{u}(\vec{\xi}_{\vec{i}})$ at the quadrature points from Stage 2, the entries of the coefficient vector $\vec{v}=A\vec{u}$ can be calculated as
\begin{equation}
  v_{g(T,(j_1,\dots,j_d))} = \sum_{i_1=1}^m \dots\sum_{i_d=1}^m A_{j_1,i_1}^{(1)}\dots A_{j_d,i_d}^{(d)} z_{(i_1,\dots,i_d)} + \dots.\label{eqn:back_transform}
\end{equation}
Here ``$\dots$'' stands for terms of a similar structure which arise due to multiplication with the derivative of a test function in the bilinear form. For simplicity those terms are not written down here; an example can be found in Eq. (10) of \cite{Muething2017}. Note that \eqref{eqn:back_transform} is the same as \eqref{eqn:naive_evaluation}, but we now multiply by the transpose of the $m\times n$ matrices $A^{(k)}$. Sum factorisation can be applied in exactly the same way, resulting in a computational complexity of $\mathcal{O}(d\cdot m^{d+1})$.

Combining the total cost of all stages we find that the overall complexity of the sum-factorised operator application is $\order(d\cdot m^{d+1})=\order(d\cdot p^{d+1})$. Finally, observe that the operator application requires
reading the local dof-vector which has $n^d$ entries. A vector
of the same size has to be written back at the end of the
calculation. Overall, the total amount of memory moved is
$\order(n^d)$. The resulting arithmetic intensity is $\order(d\cdot m)=\order(d\cdot p)$ and the operation is clearly FLOP-bound for large $p$.

In EXADUNE \cite{Bastian2014} several optimisations are applied to
speed up the operator application and exploit modern FLOP-bound
architectures. The C++ vector class library \cite{Fog2016} is used for
the multiplication with the small dense matrices $A^{(k)}$. The code is vectorised by grouping the value of the function and its three spatial derivatives and evaluating the four values
$[\partial_{x_1}\hat{u}(\vec{\xi}),\partial_{x_2}\hat{u}(\vec{\xi}),\partial_{x_3}\hat{u}(\vec{\xi}),\hat{u}(\vec{\xi})]$ simultaneously in \eqref{eqn:naive_evaluation}. Efficient transfer of data from memory is
important. In the standard DUNE implementation, a copy of the local
dof-vector $u_T$ is created on each element. While this approach is
necessary for conforming finite element methods, for the DG method
used here the local dofs of neighbouring cells do not overlap and it
is therefore possible to directly operate on the global dof-vector,
thus avoiding expensive and unnecessary memory copies.  As reported in
\cite{Muething2017}, the operator application runs at close to $60\%$
of the peak floating point performance of an Intel Haswell CPU in some cases.
\subsection{Matrix-free application and inversion of block-diagonal}
\label{sec:block_inversion_implementation}
Since in our smoothers the block system \eqref{eqn:block_system} is solved iteratively, we need to be able to apply the block-diagonal part $D_T=A_{T,T}$ of the operator $A$ in each cell $T$ in an efficient way. For this, note that
\begin{equation}
  D_{T}\vecn{u}_T = A^{\textup{v}}_{T,T}\vecn{u}_T + \sum_{F^i\in \partial T\cap\mathcal{F}_h^i} A^{\textup{if}}_{F^i;T,T} \vecn{u}_T + \sum_{F^b\in \partial T\cap\mathcal{F}_h^b} A^{\textup{b}}_{F^b;T} \vecn{u}_T.\label{eqn:block_diag_apply}
\end{equation}
Based on the decomposition in
\eqref{eqn:block_diag_apply}, in each cell $T$ we first evaluate the
contribution $A^{\textup{v}}_{T,T}\vecn{u}_T$. We then loop
over all faces $F\in \partial T$ of $T$ and -- depending on whether the
face is interior or on the boundary -- add
$A^{\textup{if}}_{F^i;T,T} \vecn{u}_T$ or $A^{\textup{b}}_{F^b;T} \vecn{u}_T$.

To (approximately) solve the system $D_{T}\vecn{\delta u}_T=\vecn{d}_T$ in
line 3 of Algorithm \ref{alg:block_jacobi} and line 4 of Algorithm \ref{alg:block_SOR} we use a Krylov subspace method. This in
turn frequently applies the block-diagonal $D_T$ as described above.
To obtain the fastest possible element-local
matrix-free solver, additionally a preconditioner for \(D_{T}\) in the
Krylov method is required. For diffusion dominated problems we find
that dividing by the diagonal of $D_T$ is efficient. For strongly
convection-dominated problems which have large off-diagonal entries,
however, we use a tridiagonal preconditioner based on the
Thomas-algorithm \cite{Press2007}. While this requires explicit
assembly of the matrix $D_T$ at the beginning of the solve, we find
that the overhead of this is small. Since only the (tri-) diagonal of
$D_T$ is retained, it also does not increase memory requirements. In the
future the overhead could be reduced further by only assembling the (tri-)
diagonal entries of $D_T$.

Since ultimately the applications of the operators $D_T$ and $D_T^{-1}$
call the methods
$A^{\textup{v}}_{T,T}$, $A^{\textup{if}}_{F^i;T,T}$ and
$A^{\textup{b}}_{F^b;T}$ of the heavily optimised \texttt{LocalOperator} class
for $A$ (see Section \ref{sec:operator_application}) those operations will
run at a similar high performance as the matrix-free application of $A$ itself.
\section{Results}\label{sec:Results}
Unless stated otherwise, we work in a domain of size
$\Omega=[0,L_x]\times[0,L_y]\times[0,L_z]$ with $L_x=L_y=1$ and $L_z=2$. The problem
sizes for different polynomial degrees are shown in
Tab. \ref{tab:problemsizes}. For each degree, the number of grid cells
is chosen such that all cells have the same size (i.e. the grid is
isotropic) and the resulting matrix has a size of
$11-14\operatorname{GB}$ (see Section \ref{sec:memory} for a detailed discussion of memory requirements). All runs are carried out on 16 cores of an
Intel Xeon E5-2698v3 (Haswell, 2.30GHz) node of the ``donkey'' cluster at
Heidelberg University. For the problems in subsections
\ref{sec:results_diffusion} and \ref{sec:results_convection} the
processor layout is chosen as $2\times 2\times 4$, i.e. each core works on a
domain of the same size and there are no load-imbalances. For the
anisotropic SPE10 benchmark in subsection \ref{sec:results_spe10} the
processor layout is changed to $4\times 4 \times 1$ to avoid parallel
partitioning in the strongly coupled vertical direction. The code was
compiled with version 6.1 of the Gnu C++ compiler.
\begin{table}
\begin{center}
\begin{tabular}{lrrrrr}
\hline
degree & \#dofs & grid size& \#dofs & matrix & matrix vs.\\
& per cell & $n_x \times n_y \times n_z$ & & size & vector size\\
\hline\hline
1 & 8 & $128\times128\times256$ & $33.6\cdot 10^6$ & $14.0$ GB &$  56\times$\\
2 & 27 & $56\times56\times112$ & $9.5\cdot 10^6$ & $13.4$ GB &$ 189\times$\\
3 & 64 & $32\times32\times64$ & $4.2\cdot 10^6$ & $14.0$ GB &$ 448\times$\\
4 & 125 & $20\times20\times40$ & $2.0\cdot 10^6$ & $13.0$ GB &$ 875\times$\\
5 & 216 & $14\times14\times28$ & $1.2\cdot 10^6$ & $13.4$ GB &$1512\times$\\
6 & 343 & $10\times10\times20$ & $0.7\cdot 10^6$ & $12.3$ GB &$2401\times$\\
7 & 512 & $8\times8\times16$ & $0.5\cdot 10^6$ & $14.0$ GB &$3584\times$\\
8 & 729 & $6\times6\times12$ & $0.3\cdot 10^6$ & $12.0$ GB &$5103\times$\\
9 & 1000 & $5\times5\times10$ & $0.2\cdot 10^6$ & $13.0$ GB &$7000\times$\\
10 & 1331 & $4\times4\times8$ & $0.2\cdot 10^6$ & $11.8$ GB &$9317\times$\\
\hline
\end{tabular}
\caption{Problem sizes for different polynomial degrees. The last column lists the relative size of the matrix and a dof-vector. See Tabs. \ref{tab:memory_diffusion}, \ref{tab:memory_convection} and \ref{tab:memory_convection_small_nrestart} below for estimates of the memory requirements in the full implementation.}
\label{tab:problemsizes}
\end{center}
\end{table}
\subsection{Solver Implementations}\label{sec:solver_implementations}
To quantify gains in performance from the techniques described above we use the following three implementations of the hybrid multigrid algorithm in Section \ref{sec:hybr-mult-prec}:
\begin{enumerate}
\item \textbf{Matrix-free (MF)}. Both the DG operator and the DG smoother are applied in an entirely matrix-free way exploiting the sum-factorisation techniques discussed in Section \ref{sec:Implementation}. In the smoother the diagonal matrices $D_T$ in each cell are inverted iteratively to tolerance $\epsilon$ with a matrix-free iterative solver as described in Section \ref{sec:matr-free-prec}. In contrast to the following two implementations, storage requirement are drastically reduced in this matrix-free version of the code since it is not necessary to store the assembled DG matrix at all. The matrix $\hat{A}$ in the lower-order subspace is assembled directly as described in Section \ref{sec:low-order-subspace}, leading to further performance improvements in the setup phase (see results in \ref{sec:setup_costs}).
\item \textbf{Matrix-explicit (MX)}. As a reference implementation we follow the standard approach of assembling the DG matrix $A$ explicitly and storing the Cholesky-/LU- factors of the diagonal blocks. This matrix and the Cholesky-/LU-factors are then used for the application of operator and smoother in the DG space. To allow for a fair comparison we speed up the assembly process with sum-factorisation. In the setup phase we calculate the coarse grid matrix $\hat{A}$ in the lower order subspace via a Galerkin product $\hat{A}=P^TAP$.
\item \textbf{Partially matrix-free (PMF)}. To quantify any gains from the matrix-free inversion of the diagonal blocks, we also store the Cholesky-/LU- factorisation of those blocks in the assembled matrix but apply the DG operator itself in a matrix-free way as described in Section \ref{sec:operator_application}. In the smoother the blocks are inverted directly by using the pre-factorised block-matrices in each cell. As in the matrix-free implementation the matrix $\hat{A}$ is assembled directly. While the storage requirements are reduced by a factor $2d+1$ compared to the matrix-explicit code, it is still necessary to store the (factorised) diagonal blocks in each grid cell.
\end{enumerate}
\subsection{Diffusion equation}\label{sec:results_diffusion}
We first solve two purely diffusive problems on a regular grid in
$d=3$ dimensions. For the Poisson equation $-\Delta u =f$ the diffusion
tensor is the identity matrix and both the advection vector and the
zero order reaction term vanish. Homogeneous Dirichlet boundary
conditions are chosen on $\partial \Omega$.

While the Poisson equation is an important model problem, the fact
that the diffusion tensor is constant and diagonal simplifies the
on-the-fly operator evaluation and might give an unrealistic advantage
to the matrix-free implementation; in particular it is not necessary
to calculate an $3\times 3$ diffusion tensor at every quadrature point. To
avoid this we also include a computationally more challenging test
case of the form $-\nabla\cdot\left(K(\vec{x})\nabla u)\right)+c(\vec{x})u=f$ with a
spatially varying positive definite diffusion tensor and (small) zero
order term given by
\begin{xalignat}{2}
  K(\vec{x}) &= \sum_{k=1}^d P^{(k)}(x_k)\vec{n}_k \vec{n}^T_k,&
  c(\vec{x}) &= \gamma \sum_{k=1}^d x_k^2\qquad\text{with $\gamma =10^{-8}$}.
\end{xalignat}
In this expression the vectors $\vec{n}_k$ with $\vec{n}_k\cdot \vec{n}_\ell=\delta_{k\ell}$ form an
orthonormal basis of $\mathbb{R}^d$ and the polynomials $P^{(k)}(\zeta)$
are positive in the entire domain. More specifically, we choose
quadratic expressions
\begin{equation}
  P^{(k)}(\zeta) = p^{(k)}_0 + p^{(k)}_2 \zeta^2\qquad\text{with}\quad p^{(k)}_0, p^{(k)}_2 > 0.
\end{equation}
In contrast to the Poisson problem, the boundary conditions are also
more complex. On the face $x_0=1$ Neumann conditions with
$j(\vec{x})=\exp(-(\vec{x}-\vec{x}_{0}^{(N)})^2/(2\sigma_N^2))$ are applied, whereas on all other
faces we use Dirichlet boundary conditions with
$g(\vec{x})=\exp(-(\vec{x}-\vec{x}_{0}^{(D)})^2/(2\sigma_D^2))$. The Gaussians in those
expressions are centred around $\vec{x}_0^{(N)}=(1,0.5,0.5)$ and
$\vec{x}_0^{(D)}=(0,0.5,0.5)$ and have a width of $\sigma_N=\sigma_D=0.1$.

In both cases the right hand side is set to a Gaussian $f(\vec{x}) =
\exp\left(-(\vec{x}-\vec{x}_0)^2/(2\sigma^2)\right)$ with $\vec{x}_0=(0.75,0.5,0.3)$ and
$\sigma=0.1$. The resulting operator is symmetric positive definite and we
use a Conjugate Gradient (CG) iteration to solve the problem to a
relative tolerance $||\vec{r}||/||\vec{r}_0||=10^{-8}$ in the $L_2(\Omega)$ norm. Since
most of the solution time is spent in the preconditioner, we achieved
a significant speedup by only evaluating $K(\vec{x})$ and $c(\vec{x})$ in the
centre of each cell in the preconditioner (but the full spatially
varying expressions are used in the operator application of the
Krylov-solver).

\begin{figure}
\begin{minipage}{0.45\linewidth}
\includegraphics[width=\linewidth]{\figdir/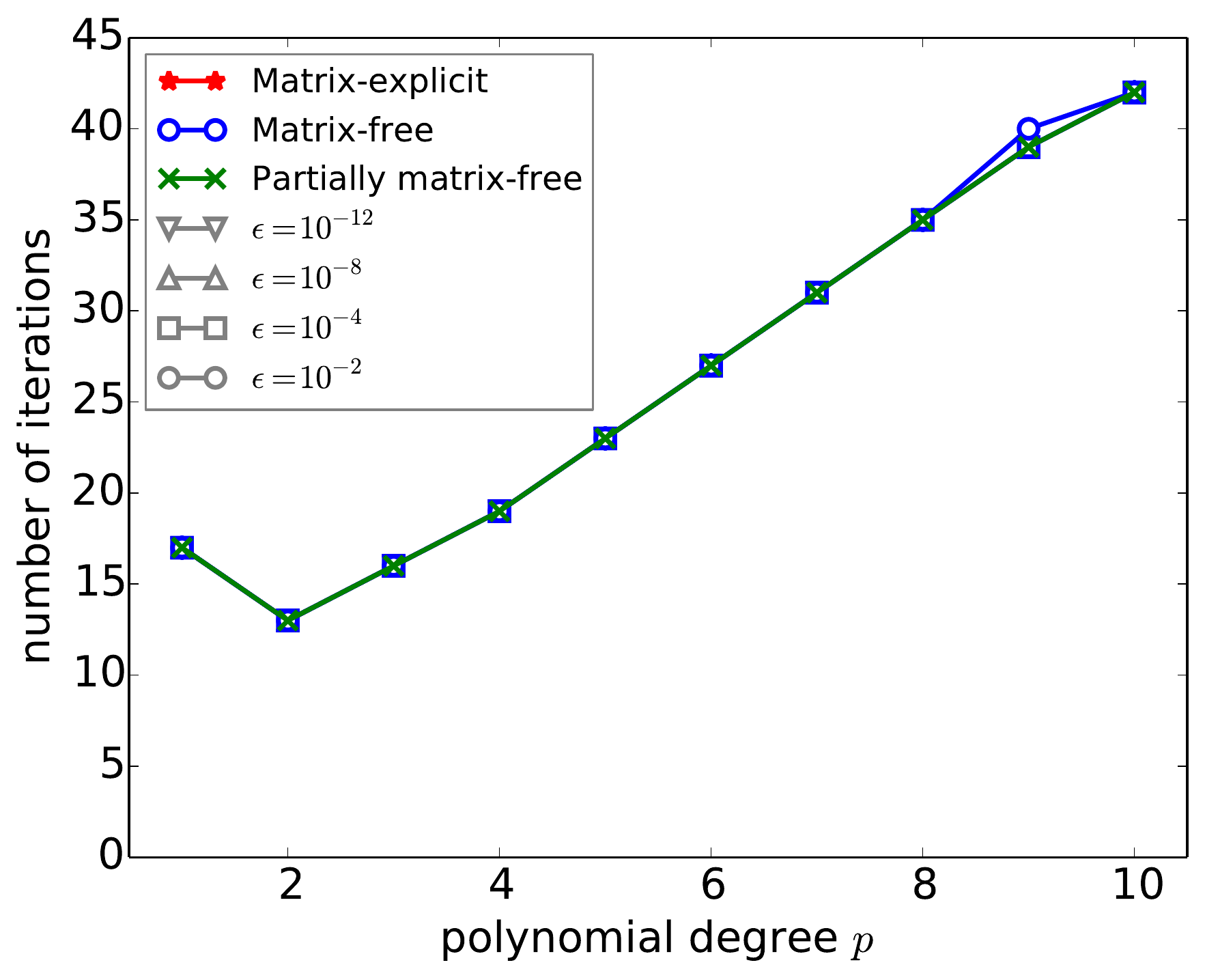}
\end{minipage}
\hfill
\begin{minipage}{0.45\linewidth}
\includegraphics[width=\linewidth]{\figdir/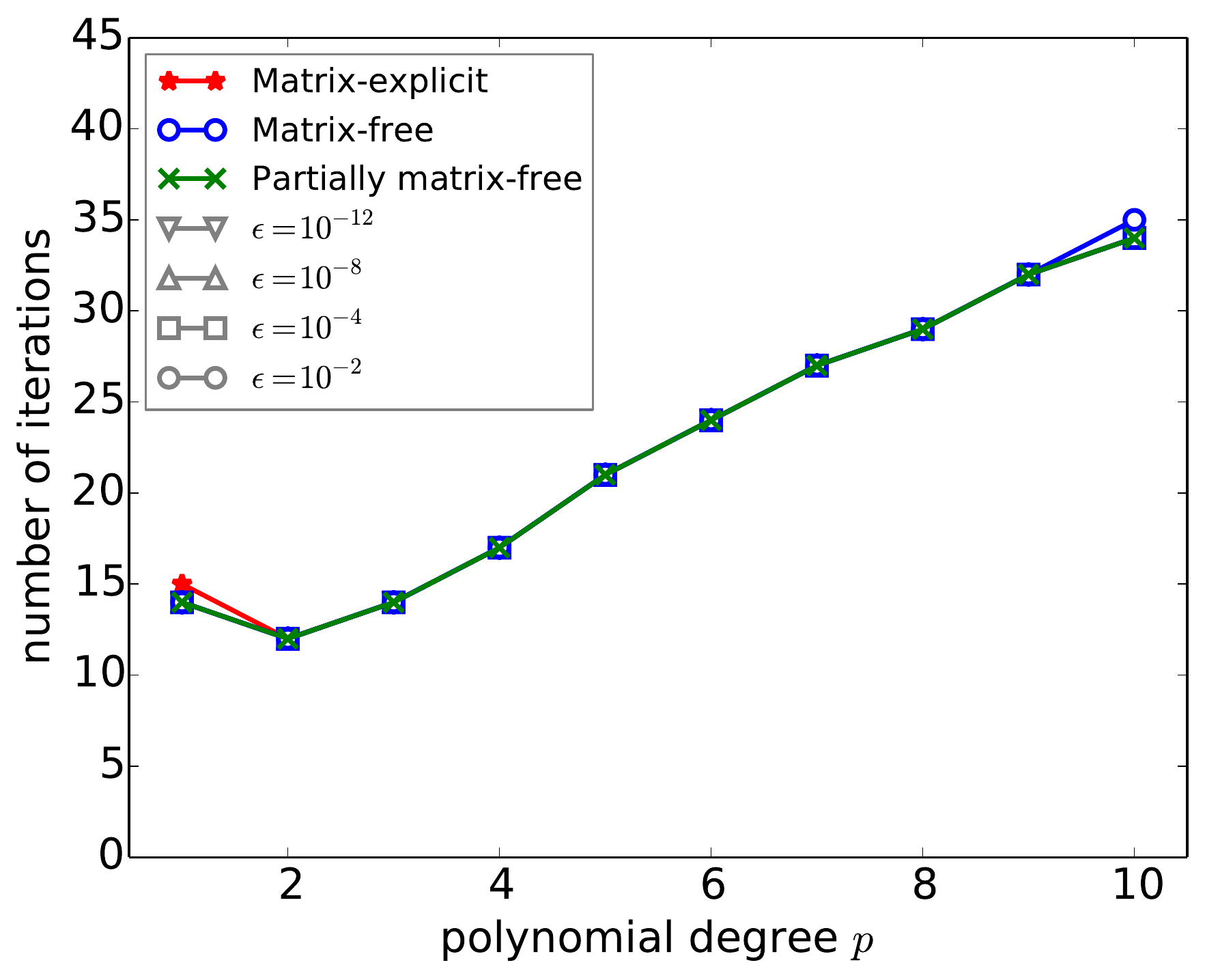}
\end{minipage}
\caption{Number of outer CG iterations for different
  implementations. Results are shown both for a constant-coefficient
  Poisson problem (left) and for the diffusion problem with varying
  coefficients (right). For the matrix-free implementation the
  block-solver tolerance $\epsilon$ is varied between $10^{-12}$ and
  $10^{-2}$.}
\label{fig:iterations_diffusion}
\end{figure}

\begin{figure}
\begin{minipage}{0.45\linewidth}
\includegraphics[width=\linewidth]{\figdir/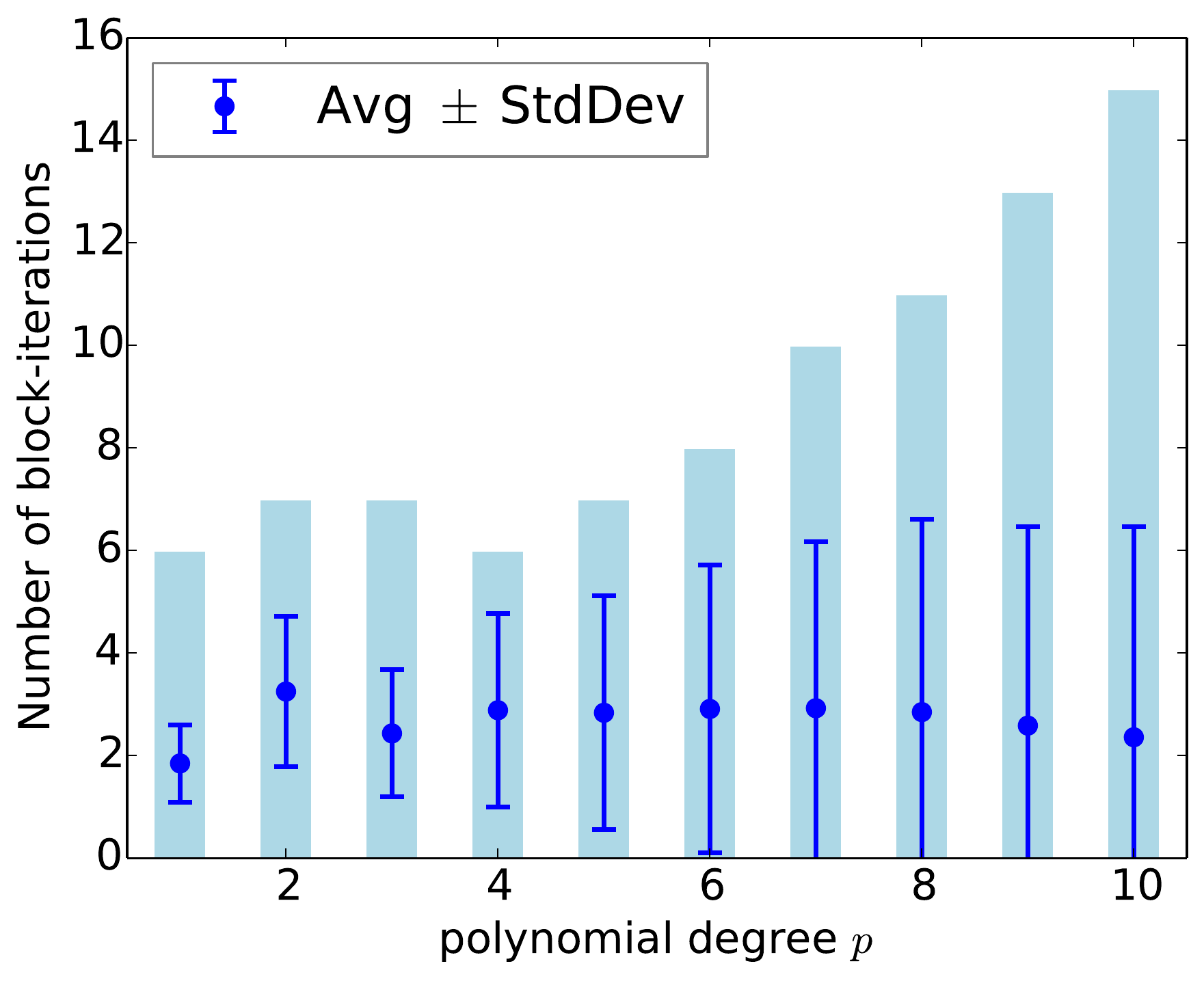}
\end{minipage}
\hfill
\begin{minipage}{0.45\linewidth}
\includegraphics[width=\linewidth]{\figdir/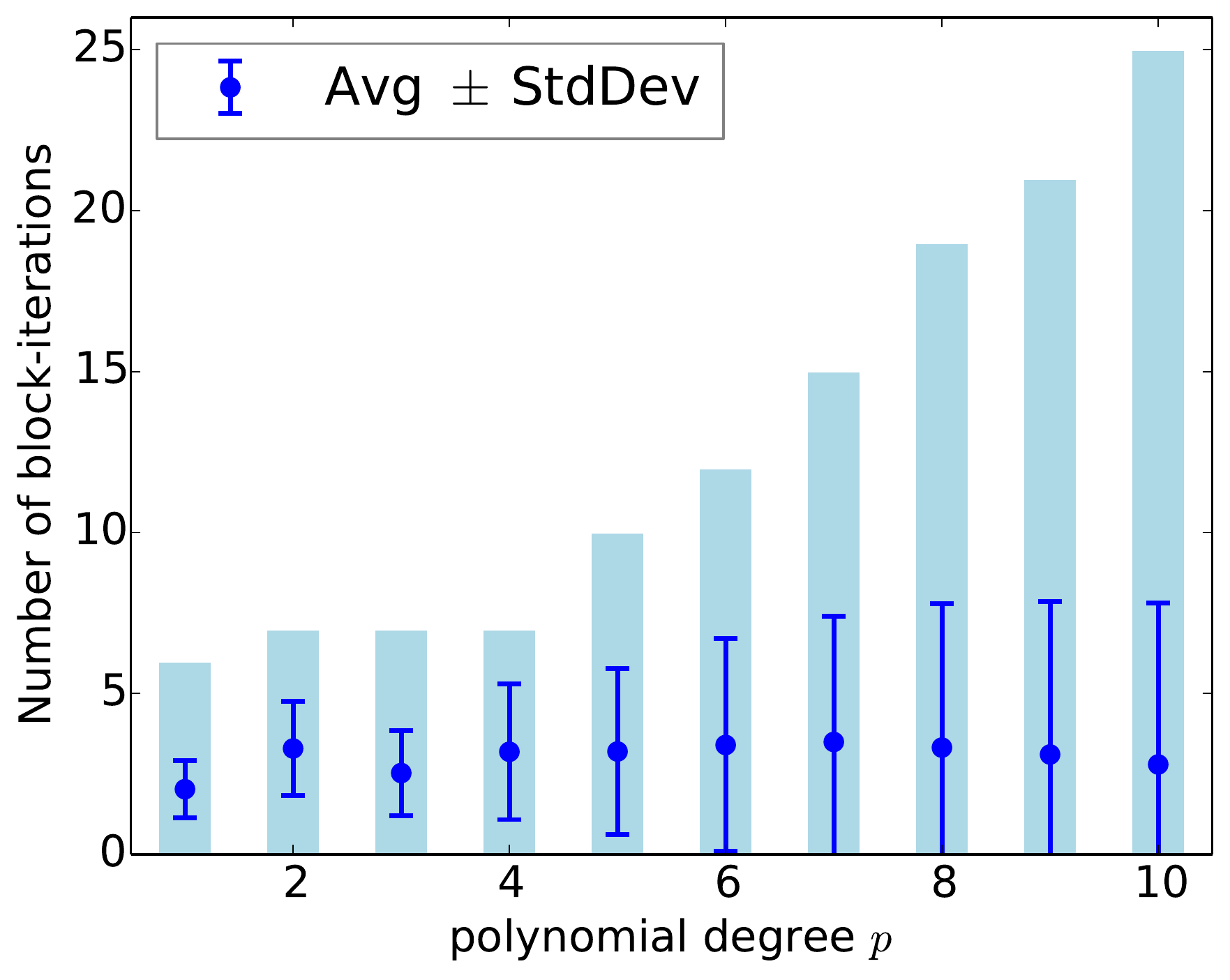}
\end{minipage}
\caption{Average number of block-solver iterations to solve to a
  tolerance of $\epsilon=10^{-2}$ for the Poisson problem (left) and the
  diffusion problem with spatially varying coefficients (right). The
  standard deviation and maximal number of solves is also shown (dark
  blue), the total range of iterations is indicated by the light blue
  bars in the background.}
\label{fig:iterations_blocksolver_diffusion}
\end{figure}

\begin{figure}
\begin{minipage}{0.45\linewidth}
\includegraphics[width=\linewidth]{\figdir/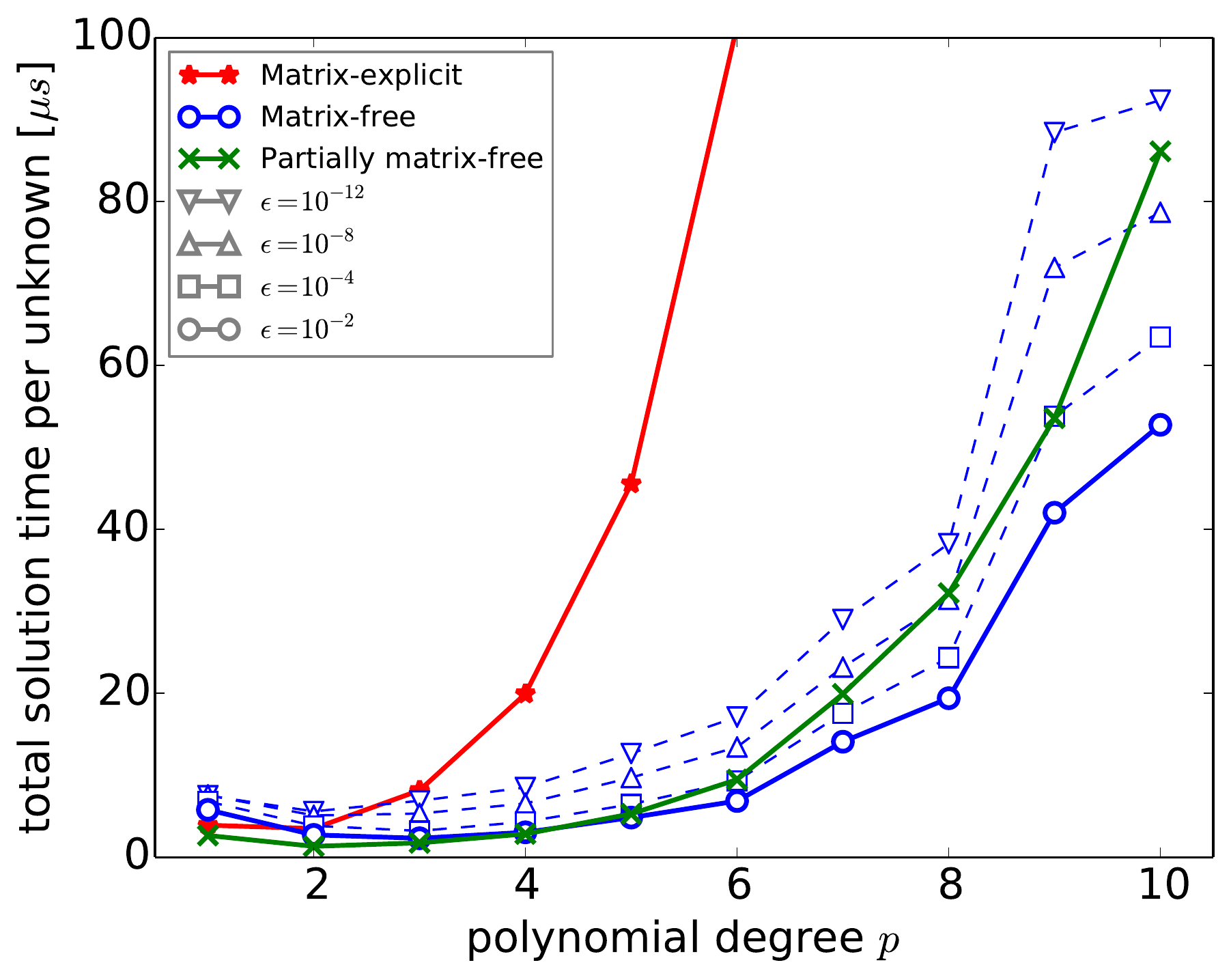}
\end{minipage}
\hfill
\begin{minipage}{0.45\linewidth}
\includegraphics[width=\linewidth]{\figdir/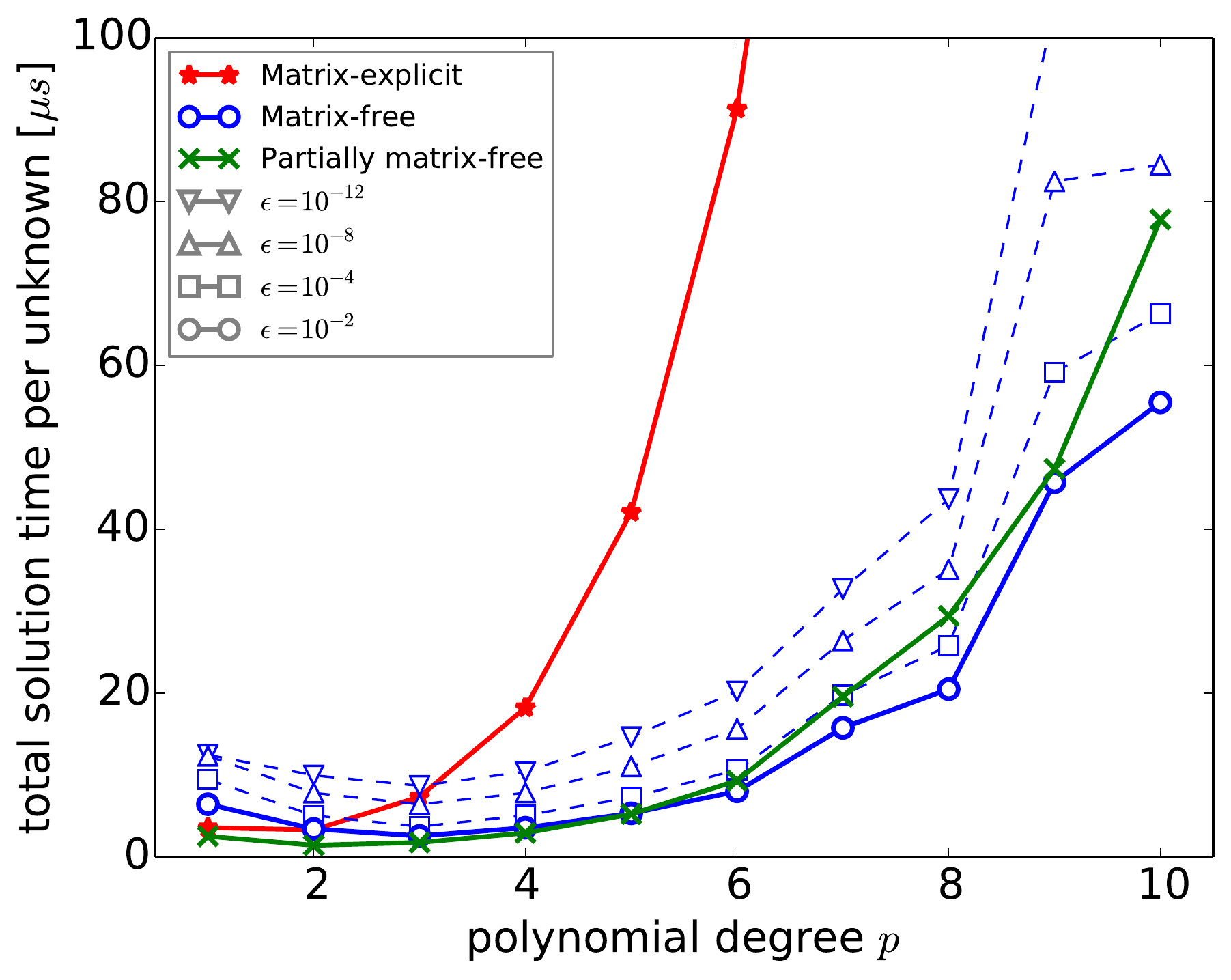}
\end{minipage}
\caption{Total solution time for different implementations and a range
  of block-solver tolerances $\epsilon$ for the Poisson problem (left) and
  the diffusion problem with spatially varying coefficients (right).}
\label{fig:solutiontime_diffusion}
\end{figure}

The block-Jacobi smoother is used in the DG space and the low order
subspace is spanned by conforming lowest order $Q_1$ elements. An
important parameter is the tolerance $\epsilon$ to which the diagonal blocks
are inverted: a relatively loose tolerance will require less
iterations of the block-solver, but might increase the total iteration
count of the outer Krylov-iteration. In
Fig. \ref{fig:iterations_diffusion} we show the number of CG
iterations for different tolerances $\epsilon$.
As can be seen from both figures, increasing the
block-solver tolerance to $\epsilon=10^{-2}$ has virtually no impact on the
number of outer CG iterations. The average number of inner
block-solver iterations which is required to reduce the residual by
two orders of magnitude is shown in
Fig. \ref{fig:iterations_blocksolver_diffusion}. In both cases the
average number of block-solver iterations is less than four for all
polynomial degrees and the maximal number of iterations never exceeds
15 (for the Poisson problem) and 25 (for the problem with spatially
varying coefficients).
The total solution time for a range of polynomial degrees is shown in
Fig. \ref{fig:solutiontime_diffusion} and (for the lowest polynomial
degrees $p$) in Tab. \ref{tab:solutiontime_diffusion}. For both
problems the partially matrix-free solver is superior to the
matrix-based version for all polynomial degrees. As expected from the
$\epsilon$-independence of the number of iterations shown in
Fig. \ref{fig:iterations_diffusion}, the matrix-free solver achieves
the best performance if the block-solver tolerance is chosen to be
very loose at $\epsilon=10^{-2}$. For higher polynomial degrees ($p\ge 5$) the
matrix-free solver gives the best overall performance and it always
beats the matrix-explicit solver except for the lowest degree
($p=1$). Although the partially matrix-free solver is fastest overall for low
polynomial degrees ($p\le 3$), Tab. \ref{tab:solutiontime_diffusion}
demonstrates that even in this case the fully matrix-free solver is only around
$2\times-3\times$
slower. It should also be kept in mind that using the matrix-free
solver also results in dramatically reduced memory requirements, as
will be discussed in Section \ref{sec:memory} below.

\begin{figure}
\begin{minipage}{0.45\linewidth}
\includegraphics[width=\linewidth]{\figdir/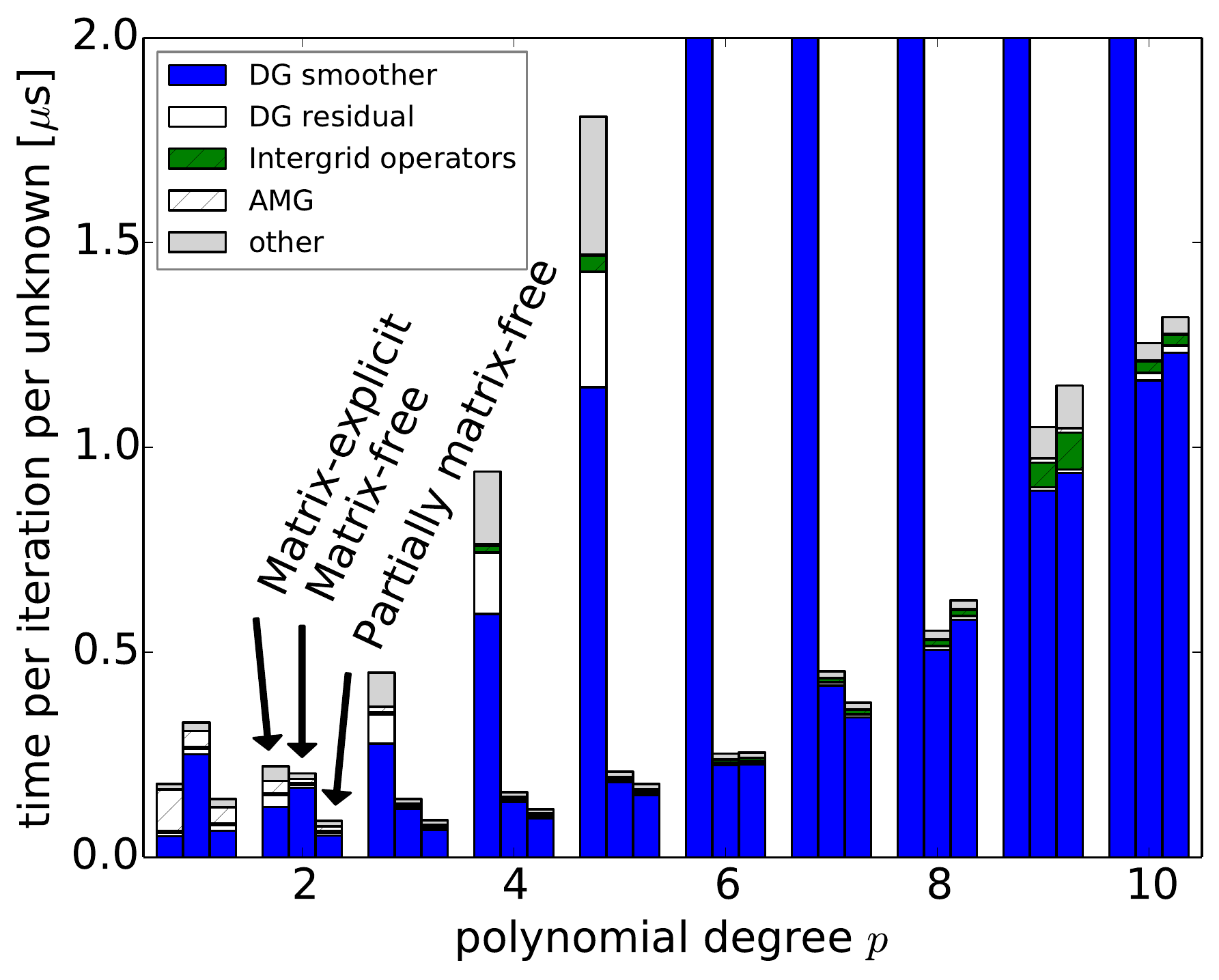}
\end{minipage}
\hfill
\begin{minipage}{0.45\linewidth}
\includegraphics[width=\linewidth]{\figdir/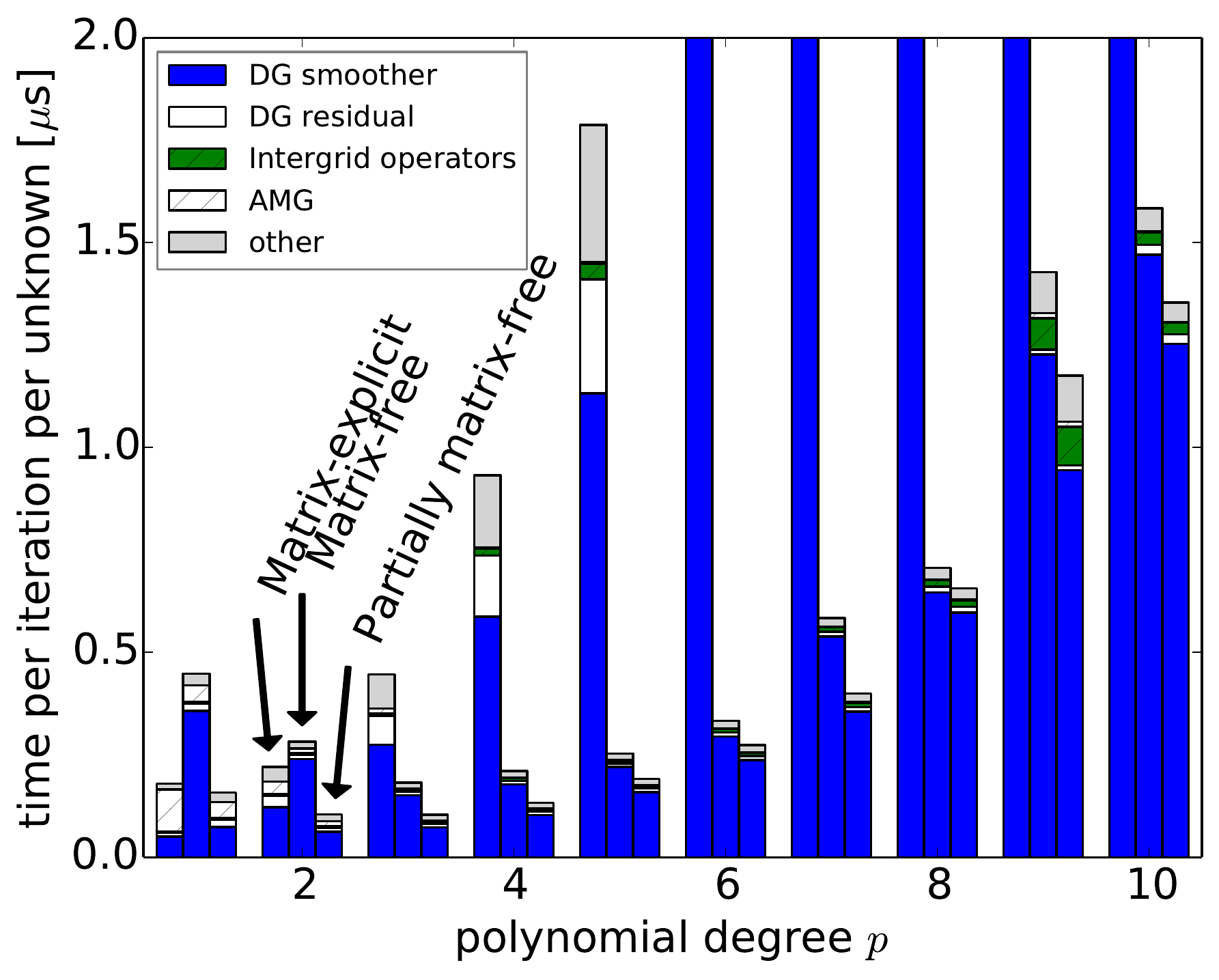}
\end{minipage}
\caption{Breakdown of the time per iteration for different
  implementations and a range of block-solver tolerances $\epsilon$ for the
  Poisson problem (left) and the diffusion problem with spatially
  varying coefficients (right).}
\label{fig:time_breakdown}
\end{figure}

\begin{table}
 \begin{center}
\begin{tabular}{ll|rrr|rrr}
\hline
& & \multicolumn{3}{|c|}{Poisson} & \multicolumn{3}{|c}{Inhomogeneous}\\
& degree $p$ & 1 & 2 & 3 & 1 & 2 & 3\\\hline\hline
Matrix-explicit &  & 3.91 & 3.50 & 8.16 & 3.60 & 3.34 & 7.37\\
\hline
 \multirow{4}{*}{Matrix-free} & $\epsilon=10^{-12}$ & 7.52 & 5.59 & 6.89 & 12.49 & 9.99 & 8.72\\
 & $\epsilon=10^{-8}$ & 7.59 & 5.10 & 5.33 & 12.35 & 7.86 & 6.45\\
 & $\epsilon=10^{-4}$ & 6.80 & 3.84 & 3.21 & 9.50 & 5.09 & 3.74\\
 & $\epsilon=10^{-2}$ & 5.79 & 2.73 & 2.31 & 6.47 & 3.45 & 2.60\\
\hline
Partially matrix-free &  & 2.66 & 1.33 & 1.75 & 2.56 & 1.46 & 1.80\\
\hline
\end{tabular}
\caption{Total solution time per unknown for the two diffusion problems and low polynomial degrees $p$. All times are given in $\mu s$.}
\label{tab:solutiontime_diffusion}
 \end{center}
\end{table}
Finally, a breakdown of the time per iteration is shown in
Fig. \ref{fig:time_breakdown}. Except for very low orders where the
AMG solve takes up a sizeable proportion of the runtime, this confirms
that the total solution time is completely dominated by the
application of the smoother and the calculation of the residual in the
DG space. A breakdown of setup costs can be found in
\ref{sec:setup_costs}.

\subsection{Convection dominated flow problem}\label{sec:results_convection}
While the block-Jacobi smoother is efficient for solving the diffusion
equations described in the previous section, it results in poor
convergence for larger values of the convection term. Here
block-SOR-type smoothers are particularly efficient and a suitable
preconditioner has to be chosen in the iterative inversion of the
block-matrices $D_T$. As a further non-trivial problem we therefore consider
a convection-dominated problem with constant coefficients which can be
written as $-\kappa \Delta u + \nabla\cdot(\vec{b} u)=f$.  The source term $f$ is adjusted such
that the exact solution is
$u(\vec{x})=\frac{x}{L_x}(1-\frac{x}{L_x})\frac{y}{L_y}(1-\frac{y}{L_y})\frac{z}{L_z}(1-\frac{z}{L_z})$. We
choose both an axis-parallel advection vector $\vec{b}=(1,0,0)$ and a more
general vector $\vec{b}=(1.0,0.5,0.3)$. The strength of the diffusion term
$\kappa$ is adjusted such that the grid Peclet-number
$\operatorname{Pe}=\max_i\{b_i\}h/\kappa$ is fixed at
$\operatorname{Pe}=2000$. This large grid-Peclet number makes the
problem heavily convection dominated and we precondition the outer Krylov-solver with two
iterations of a block-SSOR smoother instead of the hybrid multigrid
algorithm described in Section \ref{sec:hybr-mult-prec}. In the limit
$\operatorname{Pe}\rightarrow \infty$ and for a suitable ordering of the grid cells,
one block-SOR sweep over the grid would solve the equation
exactly. For the grid-Peclet number used here around 10-25 iterations
of the flexible GMRES \cite{Saad1993} solver are required to reduce
the residual by eight orders of magnitude, see
Fig. \ref{fig:convection_dominated_iterations}.  A standard GMRES
iteration (with restart 100) is used for inverting the block-diagonal
matrices $D_T$ in the matrix-free method.
\begin{figure}
  \begin{center}
    \includegraphics[width=0.8\linewidth]{\figdir/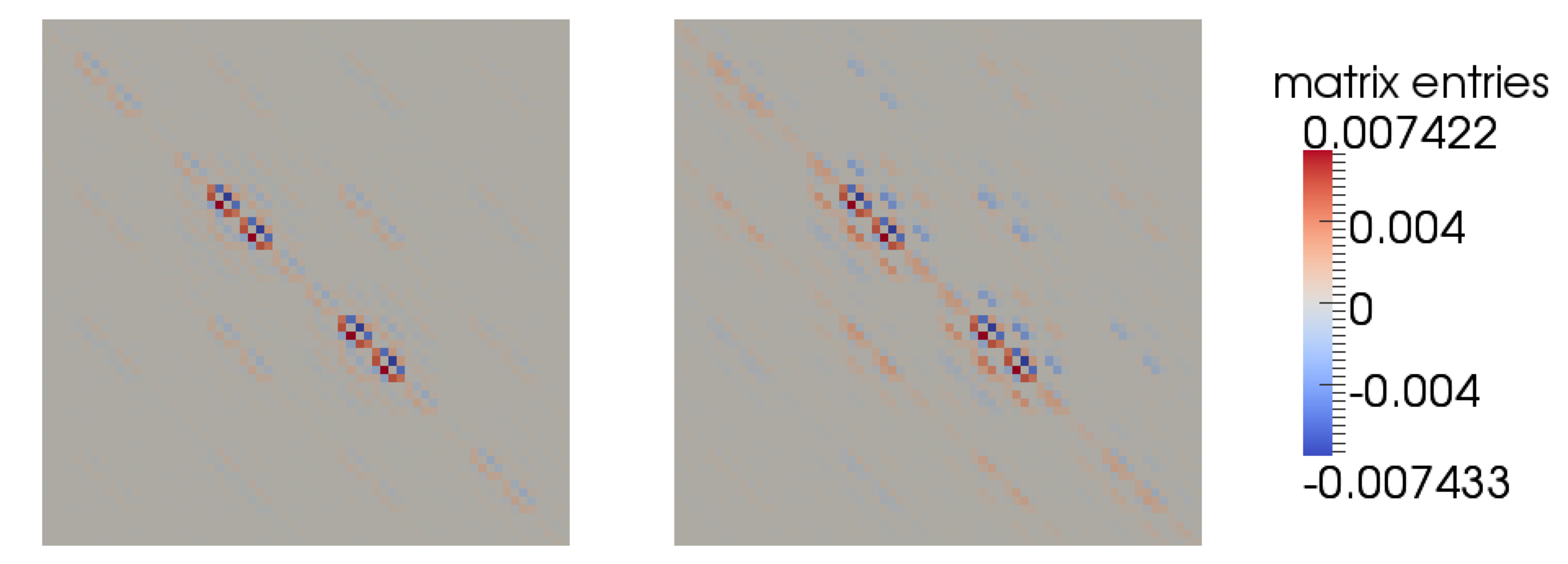}
    \caption{Block-diagonal $64\times 64$ matrix at degree $p=3$ for the
      convection dominated problem with $b=(1,0,0)$ (left) and
      $b=(1.0,0.5,0.3)$ (right). The Peclet number 2000 in both
      cases.}
    \label{fig:block_matrix}
  \end{center}
\end{figure}

The fast iterative inversion of the block-diagonal matrices $D_T$ requires a
suitable preconditioner. A typical $64\times 64$ block-matrix for degree
$n=3$ is shown in Fig. \ref{fig:block_matrix}. This matrix is clearly
not diagonally dominant, in fact is has very large off-diagonal
entries.  The simple diagonal preconditioner, which was successfully
used for all previous experiments on the diffusion problem, will
therefore not be very efficient; in fact we find that the
block-iteration will not converge at all. For the particular choice of
an x-axis aligned advection vector $\vec{b}=(1,0,0)$ and a suitable ordering
of the degrees of freedom, the first sub- and super-diagonal are
particularly large. To account for the structure of the matrix we
replaced the diagonal preconditioner by a tridiagonal solve with the
Thomas-algorithm \cite{Press2007}. While this is more expensive, the
cost of one preconditioner application is still proportional to number
of unknowns per cell.

The structure of the matrix can be explained heuristically as follows:
consider the weak form of the one-dimensional advection term
\begin{equation}
  A_{ij} = a(\tilde{\psi}^p_i,\tilde{\psi}^p_j) \propto \int \tilde{\psi}^p_i(x) \partial_x \tilde{\psi}^p_j(x) \;dx.
\end{equation}
For high polynomial degree and the nodal Gauss-Lobatto basis chosen
here, each of the basis functions $\tilde{\psi}^p_i(x)$ can be
approximated by a Gaussian centred at the $i$-th Gauss-Lobatto point
and width which half the distance to the neighbouring nodal
points. Simple symmetry arguments then imply that the diagonal entry
$A_{ii}$ is zero and the off-diagonal elements $A_{i,i+1}$ and
$A_{i,i-1}$ are large and of approximately opposite
size. Qualitatively this explains the structure seen in
Fig. \ref{fig:block_matrix}. For the non-axis-aligned cases additional
off-diagonal bands appear at a distance of $p+1$ and $(p+1)^2$ from
the main diagonal.

\begin{figure}
\begin{minipage}{0.45\linewidth}
\begin{center}
 \includegraphics[width=1.0\linewidth]{\figdir/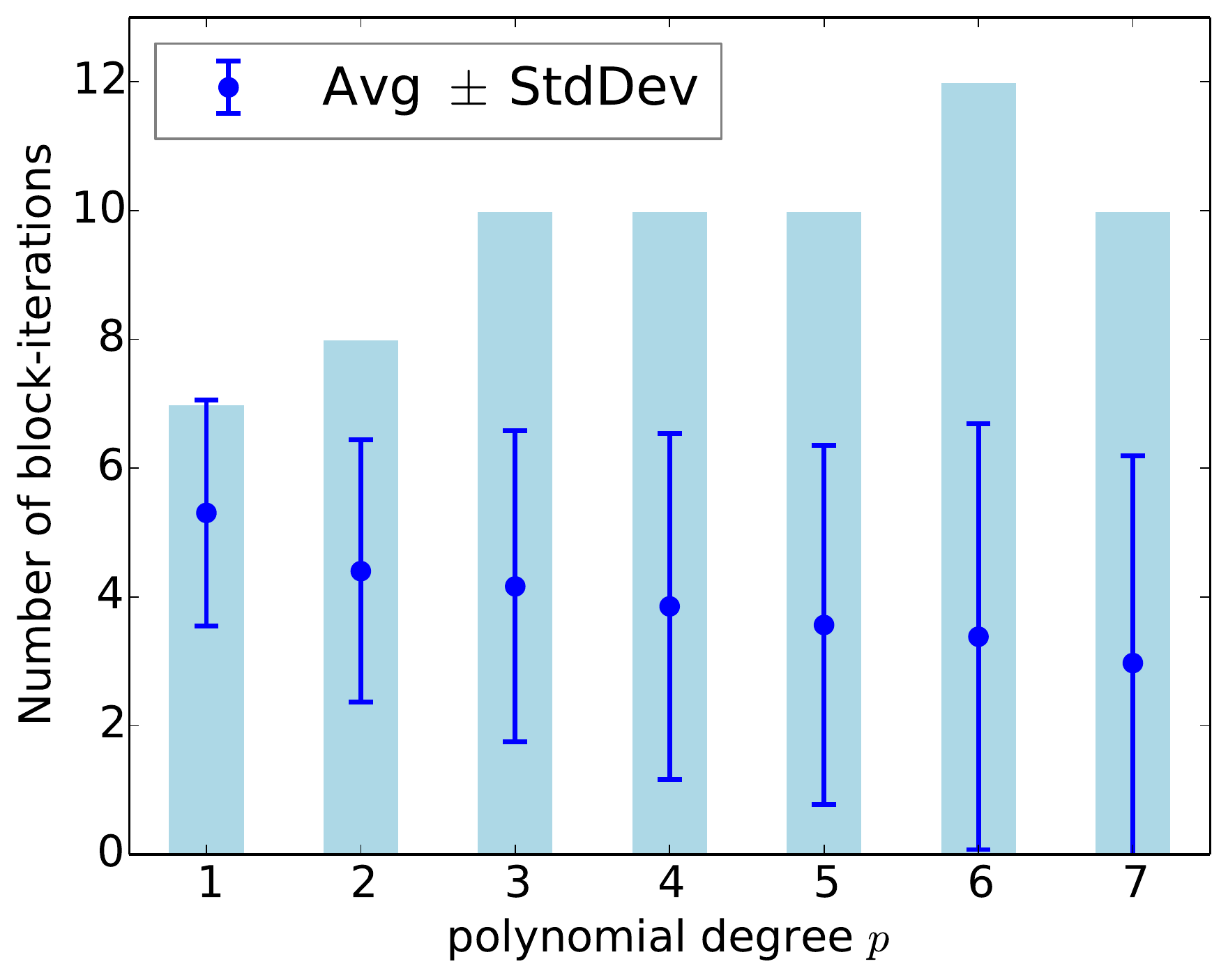}
\end{center}
\end{minipage}
\hfill
\begin{minipage}{0.45\linewidth}
\begin{center}
 \includegraphics[width=1.0\linewidth]{\figdir/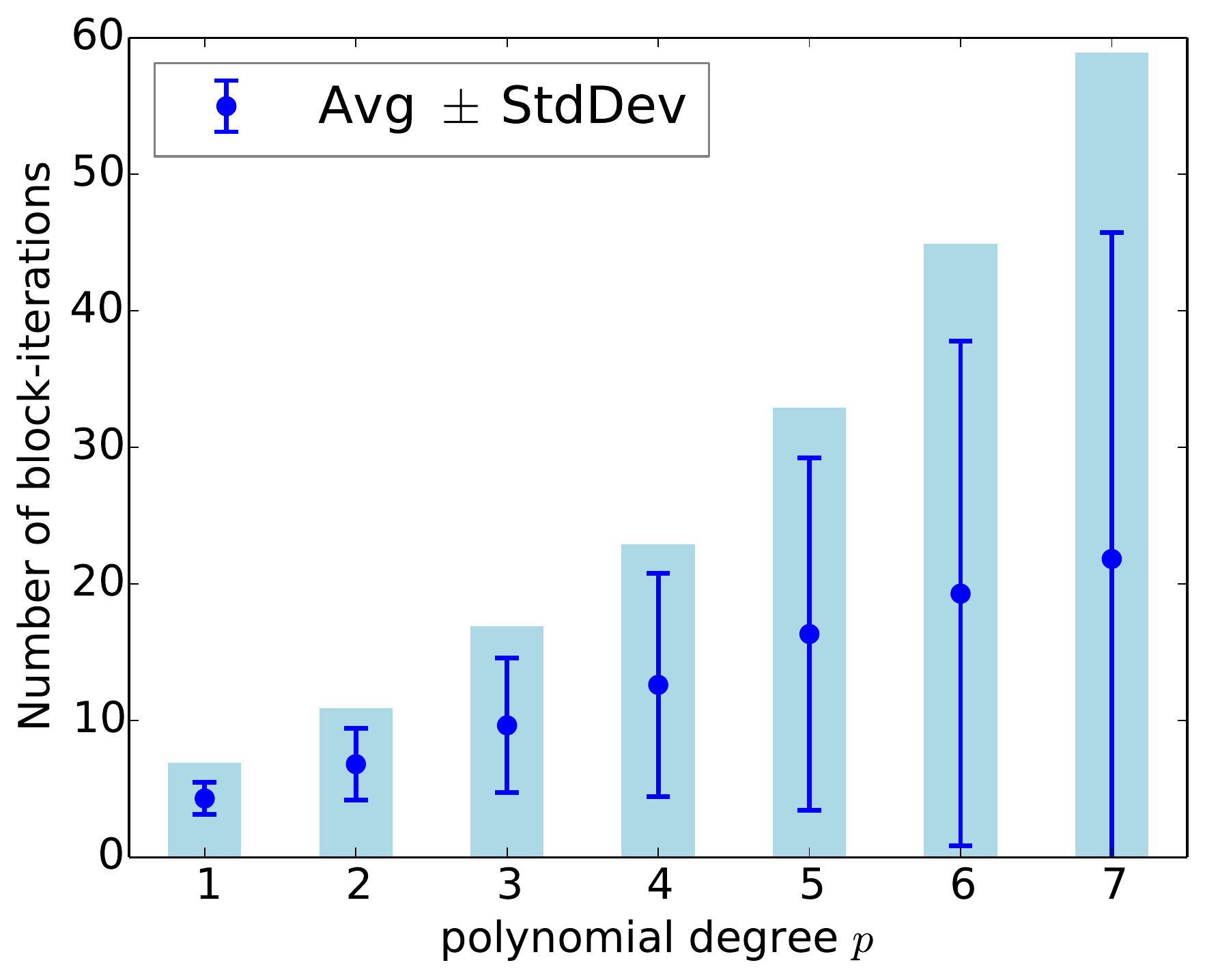}
\end{center}
\end{minipage}
 \caption{Number of GMRES iterations for the block solver with tri-diagonal preconditioner and tolerance $\epsilon=10^{-2}$. The convection-dominated problem is solved with $\vec{b}=(1,0,0)$ (left) and $\vec{b}=(1.0,0.5,0.3)$ (right).}
 \label{fig:convection_dominated_blockiter}
\end{figure}

\begin{figure}
\begin{minipage}{0.45\linewidth}
\begin{center}
 \includegraphics[width=1.0\linewidth]{\figdir/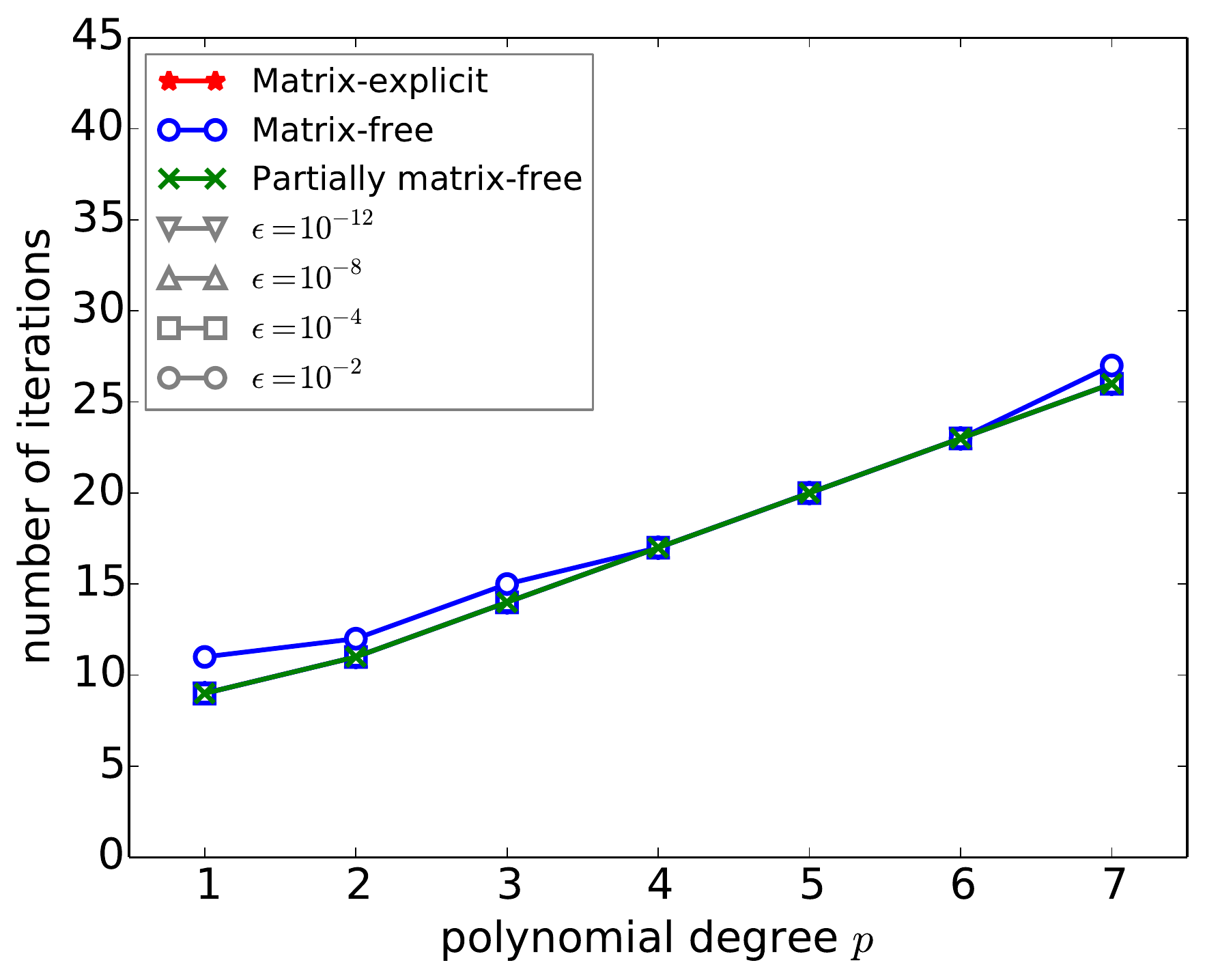}
\end{center}
\end{minipage}
\hfill
\begin{minipage}{0.45\linewidth}
\begin{center}
 \includegraphics[width=1.0\linewidth]{\figdir/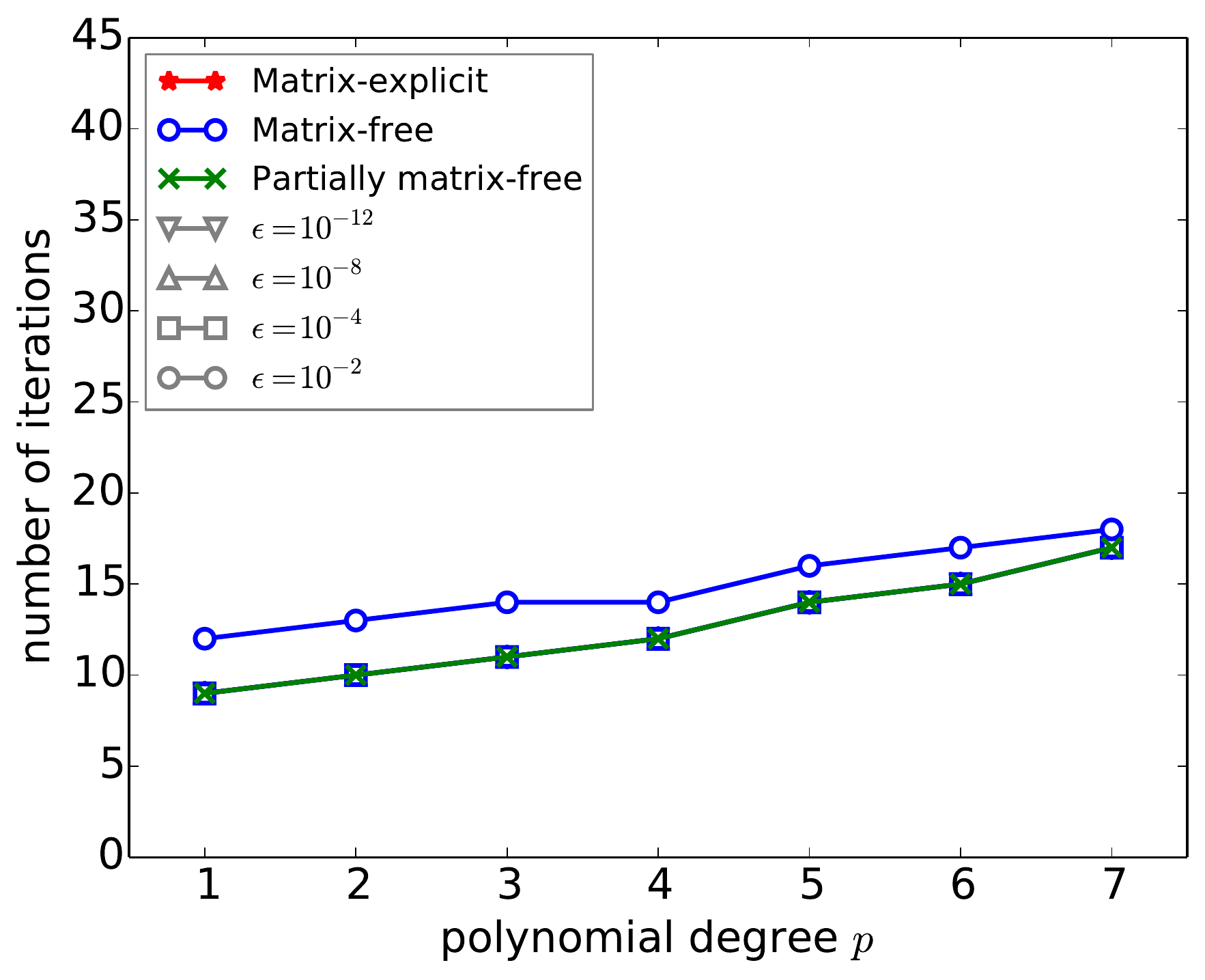}
\end{center}
\end{minipage}
 \caption{Number of outer flexible GMRES iterations. The convection-dominated problem is solved with $b=(1,0,0)$ (left) and $b=(1.0,0.5,0.3)$ (right).}
 \label{fig:convection_dominated_iterations}
\end{figure}

\begin{figure}
\begin{minipage}{0.45\linewidth}
\begin{center}
 \includegraphics[width=1.0\linewidth]{\figdir/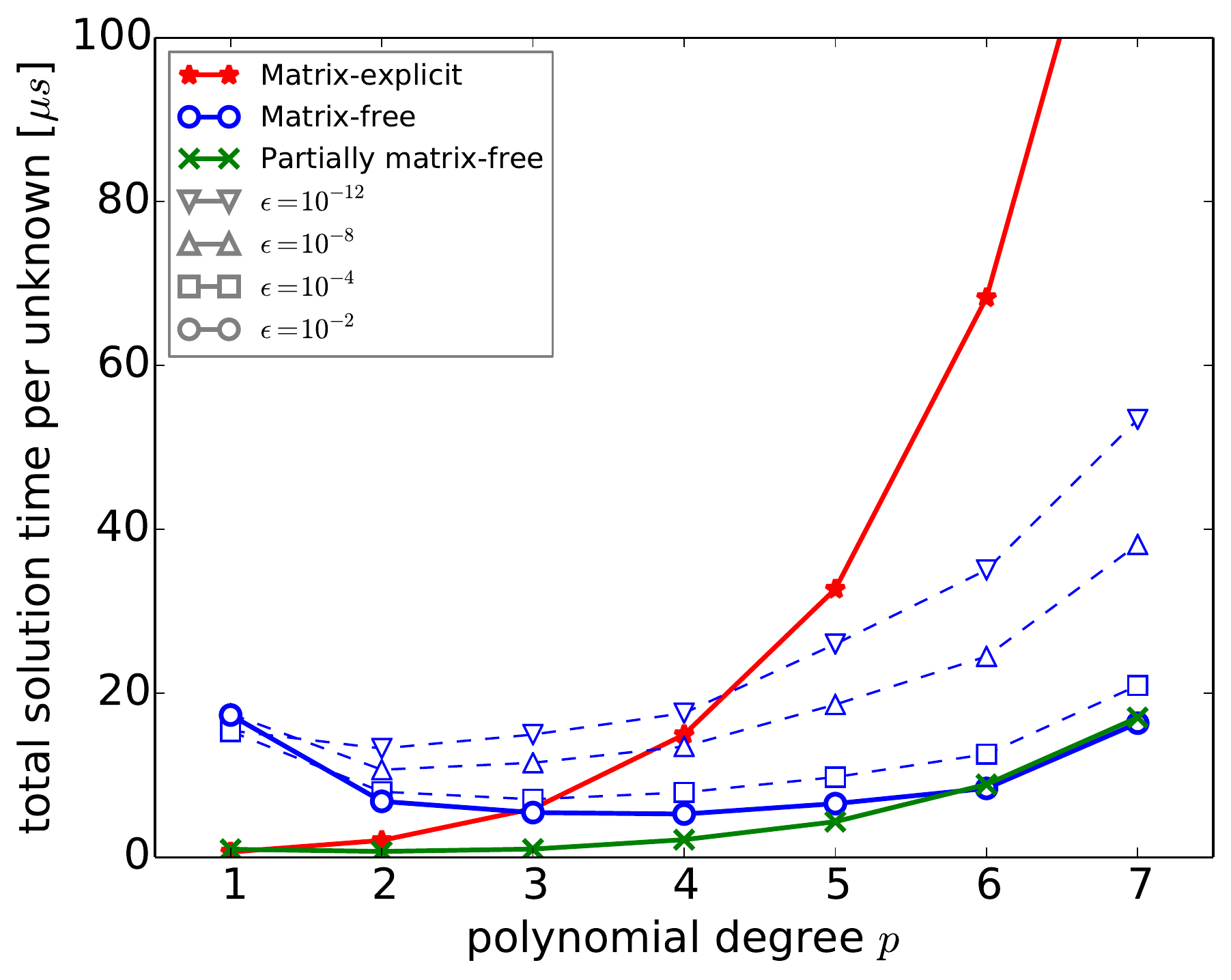}
\end{center}
\end{minipage}
\hfill
\begin{minipage}{0.45\linewidth}
\begin{center}
 \includegraphics[width=1.0\linewidth]{\figdir/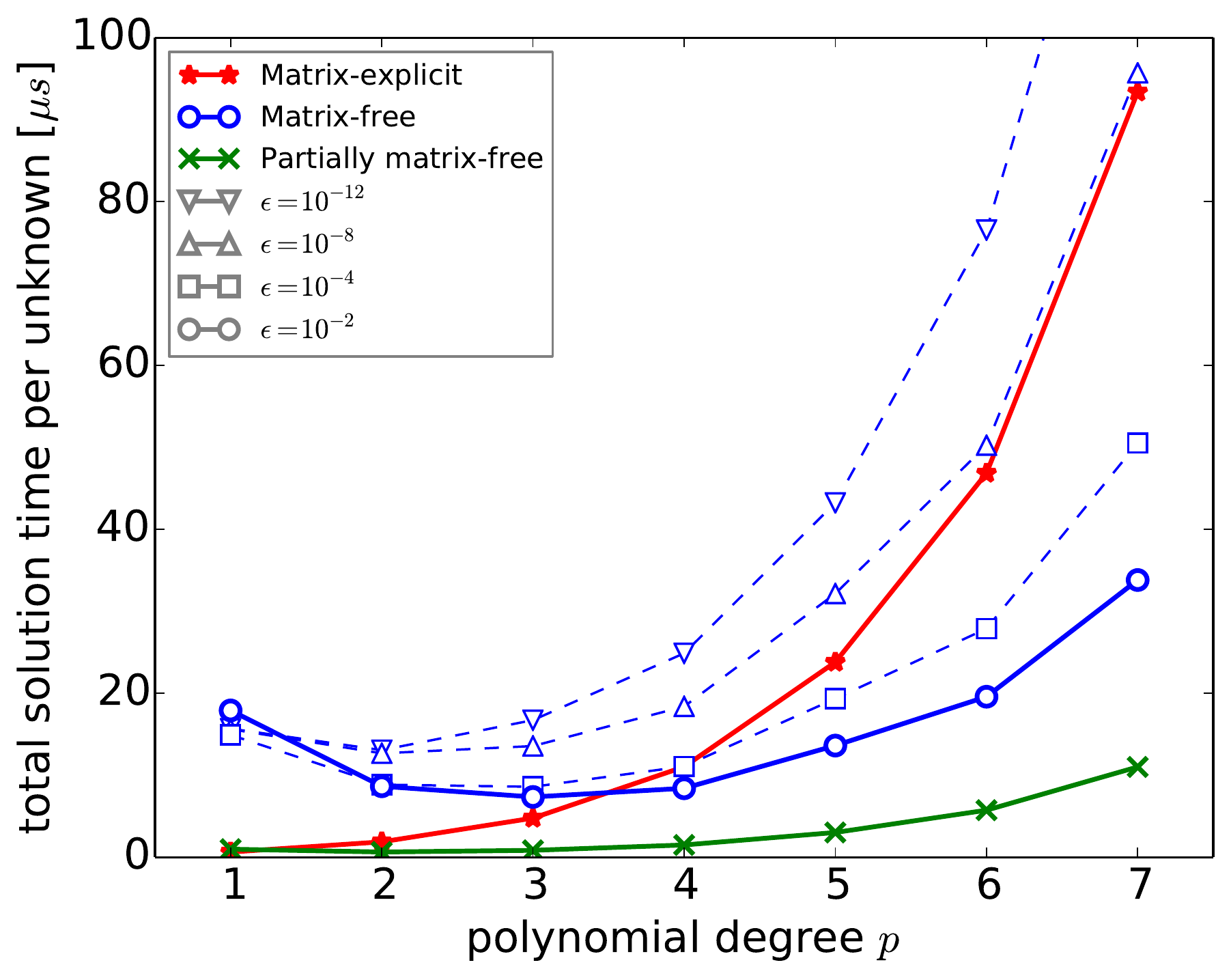}
\end{center}
\end{minipage}
 \caption{Total solution time per unknown. The convection-dominated problem is solved with $\vec{b}=(1,0,0)$ (left) and $\vec{b}=(1.0,0.5,0.3)$ (right).}
 \label{fig:convection_dominated_solutiontime}
\end{figure}

The number of block-solver iterations for the convection dominated
problem is shown in Fig. \ref{fig:convection_dominated_blockiter}
which should be compared to
Fig. \ref{fig:iterations_blocksolver_diffusion}. On average around
four iterations are required to reduce the residual by two orders of
magnitude for the axis-aligned advection vector $\vec{b}=(1,0,0)$ and there
is no strong dependence on the polynomial degree. For the more general
advection vector $\vec{b}=(1.0,0.5,0.3)$ the number of iterations grows
since the tridiagonal preconditioner becomes less efficient as
expected due to the additional sub-diagonals in the matrix.
The number of outer flexible GMRES iterations to solve the full system
to a tolerance of $10^{-8}$ is shown in
Fig. \ref{fig:convection_dominated_iterations}.  For a given
polynomial degree the number of iterations is virtually identical for
all solvers and tolerances $\epsilon$ for the block-diagonal solver. The only
exception is the matrix-free solver with tolerance $\epsilon=10^{-2}$ which
leads to a $\approx 20\%$ increase in the iteration count.
Finally, the total solution time is shown in
Fig. \ref{fig:convection_dominated_solutiontime} and
Tab. \ref{tab:solutiontime_convection}. In all cases the partially
matrix-free solver clearly beats the matrix-based implementation. As
expected from Fig. \ref{fig:convection_dominated_iterations} the
performance of the matrix-free code mainly depends on the block-solver
tolerance $\epsilon$.  For the loosest tolerance $\epsilon=10^{-2}$ it is faster
than the matrix-explicit solver for moderately high polynomial degrees
($p\ge 3$ for $\vec{b}=(1,0,0)$ and $p\ge 4$ for $\vec{b}=(1.0,0.5,0.3)$). For high
polynomial degree, the matrix-free solve time is in the same ballpark
as the partially matrix-free time for the axis aligned advection
vector $\vec{b}=(1,0,0)$. For $\vec{b}=(1.0,0.5,0.3)$ however, the matrix-free
solver is substantially less efficient than the partially matrix-free
version. This is can be traced back to the growth in the number of iterations in the matrix-free inversion of $D_T$ shown in Fig. \ref{fig:convection_dominated_blockiter} (right). Again it should be kept in mind that the reduced memory
requirements of the matrix-free solver allow the solution of larger
problems.

\begin{table}
 \begin{center}
\begin{tabular}{ll|rrr|rrr}
\hline
& & \multicolumn{3}{|c|}{$\vec{b}=(1,0,0)$} & \multicolumn{3}{|c}{$\vec{b}=(1.0,0.5,0.3)$}\\
& degree $p$ & 2 & 3 & 4 & 2 & 3 & 4\\\hline\hline
Matrix-explicit &  & 2.06 & 5.92 & 14.97 & 1.88 & 4.77 & 11.04\\
\hline
 \multirow{4}{*}{Matrix-free} & $\epsilon=10^{-12}$ & 13.27 & 14.97 & 17.59 & 13.07 & 16.73 & 24.89\\
 & $\epsilon=10^{-8}$ & 10.68 & 11.53 & 13.51 & 12.68 & 13.57 & 18.39\\
 & $\epsilon=10^{-4}$ & 7.99 & 7.07 & 7.91 & 8.86 & 8.60 & 11.08\\
 & $\epsilon=10^{-2}$ & 6.82 & 5.45 & 5.27 & 8.66 & 7.36 & 8.42\\
\hline
Partially matrix-free &  & 0.71 & 1.00 & 2.14 & 0.64 & 0.84 & 1.50\\
\hline
\end{tabular}
\caption{Total solution time per unknown for convection dominated problem. All times are given in $\mu s$.}
\label{tab:solutiontime_convection}
 \end{center}
\end{table}
\subsection{Memory savings in matrix-free implementation}\label{sec:memory}
\newcommand{\Nvector}{M_\text{vector}}
\newcommand{\Nlocalvector}{M_\text{local vector}}
\newcommand{\Nmatrix}{M_\text{matrix}}
\newcommand{\Nblockdiag}{M_\text{blockdiag}}
\newcommand{\Ncell}{N_\text{cell}}
As argued above, a key advantage of the completely matrix-free solver 
is its significantly lower demand on storage space. To estimate the memory
requirements for all implementations carefully, first observe that 
on a $d$-dimensional grid with $\Ncell$ cells an entire
vector in DG space with polynomial degree $p$ consists
of $\Nvector=\Ncell (p+1)^d$ double precision numbers. The corresponding
sizes for the full DG matrix and its blockdiagonal (used in the partially
matrix-free implementation) are $\Nmatrix=7\Ncell(p+1)^{2d}$ and
$\Nblockdiag=\Ncell(p+1)^{2d}$. In the following estimates we will take into account vectors and matrices on the first coarse multigrid level (i.e. in the $Q_1$ and $P_0$ subspace); we can safely neglect all data stored on the
next-coarser levels. Since all grids used in the numerical experiments above
contain at least 100 cells, we also ignore temporary storage for cell-local
vectors.

First consider the setup for the diffusion problem in Section \ref{sec:results_diffusion}. The solution
and right hand side, both of which are of size $\Nvector$, have to be stored
in every implementation. The outer CG iteration requires two additional
temporaries of size $\Nvector$ each. Likewise, the hybrid multigrid stores
three temporaries of size $\Nvector$ (two in the completely matrix-based
implementation) and the corresponding vector storage in the $Q_1$ subspace is
$3\Ncell$ (three vectors with one unknown per cell). While the DG system
matrix of size
$\Nmatrix$ is only stored in the matrix-based implementation, the coarse
level matrix of size $27\Ncell$, which describes the $Q_1$ discretisation of
the problem, is required by the AMG solver in all cases.
In the partially matrix-free case we store the entries of the local
block-diagonal of total size $\Nblockdiag$; in the matrix-free case we store
the diagonal only (for preconditioning the cell-local iterative solver),
which is equivalent to storing another vector of size $\Nvector$.
Overall, this results in the following storage requirements for the
matrix-explicit (MF), partially matrix-free (PMF) and fully matrix-free (MF)
cases:
\begin{equation}
\begin{aligned}
M^{(\text{diffusion})}_{\text{total}} &= 
\begin{cases}
\Nmatrix + 6\Nvector + 30\Ncell & \\
\Nblockdiag + 7\Nvector + 30\Ncell & \\
8\Nvector + 30\Ncell &
\end{cases}
= \begin{cases}
\left(7(p+1)^{2d}+6(p+1)^d+30\right)\Ncell & \text{[MX]}\\
\left((p+1)^{2d}+7(p+1)^d+30\right)\Ncell & \text{[PMF]}\\
\left(8(p+1)^d+30\right)\Ncell & \text{[MF]}
\end{cases}
\end{aligned}
\label{eqn:memory_diffusion}
\end{equation}
Numerical values for all considered problem sizes (see Tab. \ref{tab:problemsizes}) in $d=3$ dimensions are shown in Tab. \ref{tab:memory_diffusion}. The table shows that even for moderate polynomial degrees the relative saving of the partially matrix-free implementation approaches a constant factor of $\approx 7\times$, since only the block-diagonal has to be stored. The relative saving for the fully matrix-free method is significantly larger, and of two- to three orders of magnitiude, even for moderately large degrees. Note that asymptotically the savings grow with the third power of $p+1$ as $\approx \frac{7}{8}(p+1)^d$.
\begin{table}
\begin{center}
\begin{tabular}{rr|r|rr|rr}
\hline
degree $p$ & \# dofs & \multicolumn{5}{c}{memory requirements $M_{\text{total}}^{(\text{diffusion})}$} \\
           &         & \multicolumn{1}{|c|}{MX} & \multicolumn{2}{|c|}{PMF} & \multicolumn{2}{|c}{MF} \\
\hline\hline
$ 1$ & $ 33.6 \cdot 10^{6}$ & $ 17.6$ GB & $  5.0$ GB & $[  3.5 \times]$ & $3.154$ GB &  $[   5.6 \times]$ \\
$ 2$ & $  9.5 \cdot 10^{6}$ & $ 14.9$ GB & $  2.7$ GB & $[  5.6 \times]$ & $0.691$ GB &  $[  21.5 \times]$ \\
$ 3$ & $  4.2 \cdot 10^{6}$ & $ 15.2$ GB & $  2.4$ GB & $[  6.4 \times]$ & $0.284$ GB &  $[  53.7 \times]$ \\
$ 4$ & $  2.0 \cdot 10^{6}$ & $ 14.1$ GB & $  2.1$ GB & $[  6.7 \times]$ & $0.132$ GB &  $[ 106.9 \times]$ \\
$ 5$ & $  1.2 \cdot 10^{6}$ & $ 14.4$ GB & $  2.1$ GB & $[  6.8 \times]$ & $0.077$ GB &  $[ 186.5 \times]$ \\
$ 6$ & $  0.7 \cdot 10^{6}$ & $ 13.2$ GB & $  1.9$ GB & $[  6.9 \times]$ & $0.044$ GB &  $[ 297.6 \times]$ \\
$ 7$ & $  0.5 \cdot 10^{6}$ & $ 15.1$ GB & $  2.2$ GB & $[  6.9 \times]$ & $0.034$ GB &  $[ 445.5 \times]$ \\
$ 8$ & $  0.3 \cdot 10^{6}$ & $ 12.9$ GB & $  1.9$ GB & $[  6.9 \times]$ & $0.020$ GB &  $[ 635.4 \times]$ \\
$ 9$ & $  0.2 \cdot 10^{6}$ & $ 14.0$ GB & $  2.0$ GB & $[  7.0 \times]$ & $0.016$ GB &  $[ 872.5 \times]$ \\
$10$ & $  0.2 \cdot 10^{6}$ & $ 12.7$ GB & $  1.8$ GB & $[  7.0 \times]$ & $0.011$ GB &  $[1162.1 \times]$ \\
\hline
\end{tabular}
\caption{Estimated memory requirements as defined in Eq. \eqref{eqn:memory_diffusion} for different implementations for the solution of the diffusion problem in Section \ref{sec:results_diffusion}. The savings relative to the matrix-explicit implementation are shown in square brackets.}
\label{tab:memory_diffusion}
\end{center}
\end{table}

The convection-dominated problem in Section \ref{sec:results_convection} uses a slightly different setup. The number of temporary vectors in the flexible GMRES solver with a restart value of $n_{\text{restart}}^{(\text{outer})}=100$ is $2n_{\text{restart}}^{(\text{outer})}+2$ for the outer solve. In contrast to the diffusion problem in Section \ref{sec:results_diffusion}, a block-SSOR smoother is used instead of the hybrid multigrid algorithm. This reduces the storage requirements since the SSOR smoother operates directly on the solution vector, and does not need the temporary storage for additional vectors and coarse level matrices/vectors which are required in the multigrid algorithm. However, since the local blocks are inverted iteratively with a GMRES method which uses $n_{\text{restart}}^{(\text{block})}+1$ temporary local vectors of size $\Nlocalvector=(p+1)^d$ and the restart value $n_{\text{restart}}^{(\text{block})}=100$ in our numerical experiments is (unnecessarily) large, we choose to count those vectors in the fully matrix-free method. Furthermore, since the iterative block-diagonal solvers are preconditioned with the tridiagonal Thomas algorithm, we need two additional global vectors to store the bands of off-diagonal entries in the matrix-free algorithm.
Overall, this results in the following storage requirements:
\begin{equation}
\begin{aligned}
M^{(\text{convection})}_{\text{total}} &= 
\begin{cases}
\Nmatrix  + (2n_{\text{restart}}^{(\text{outer})}+4)\Nvector & \text{[MX]}\\
\Nblockdiag + (2n_{\text{restart}}^{(\text{outer})}+4)\Nvector & \text{[PMF]} \\
(2n_{\text{restart}}^{(\text{outer})}+7)\Nvector + (n_{\text{restart}}^{(\text{block})}+1)\Nlocalvector & \text{[MF]}
\end{cases}\\[1ex]
&= \begin{cases}
\left(7(p+1)^{2d}+(2n_{\text{restart}}^{(\text{outer})}+4)(p+1)^d \right)\Ncell & \text{[MX]}\\
\left((p+1)^{2d}+(2n_{\text{restart}}^{(\text{outer})}+4)(p+1)^d \right)\Ncell & \text{[PMF]}\\
(2n_{\text{restart}}^{(\text{outer})}+7)(p+1)^d \Ncell + (n_{\text{restart}}^{(\text{block})}+1)(p+1)^d & \text{[MF]}
\end{cases}
\end{aligned}
\label{eqn:memory_convection}
\end{equation}
\begin{table}
\begin{center}
\begin{tabular}{rr|r|rr|rr}
\hline
degree $p$ & \# dofs & \multicolumn{5}{c}{memory requirements $M_{\text{total}}^{(\text{convection})}$} \\
           &         & \multicolumn{1}{|c|}{MX} & \multicolumn{2}{|c|}{PMF} & \multicolumn{2}{|c}{MF} \\
\hline\hline
$ 1$ & $ 33.6 \cdot 10^{6}$ & $ 65.0$ GB & $ 53.0$ GB & $[  1.2 \times]$ & $51.750$ GB &  $[   1.3 \times]$ \\
$ 2$ & $  9.5 \cdot 10^{6}$ & $ 27.8$ GB & $ 16.3$ GB & $[  1.7 \times]$ & $14.626$ GB &  $[   1.9 \times]$ \\
$ 3$ & $  4.2 \cdot 10^{6}$ & $ 20.4$ GB & $  8.4$ GB & $[  2.4 \times]$ & $6.469$ GB &  $[   3.1 \times]$ \\
$ 4$ & $  2.0 \cdot 10^{6}$ & $ 16.1$ GB & $  4.9$ GB & $[  3.3 \times]$ & $3.085$ GB &  $[   5.2 \times]$ \\
$ 5$ & $  1.2 \cdot 10^{6}$ & $ 15.2$ GB & $  3.7$ GB & $[  4.1 \times]$ & $1.828$ GB &  $[   8.3 \times]$ \\
$ 6$ & $  0.7 \cdot 10^{6}$ & $ 13.3$ GB & $  2.8$ GB & $[  4.8 \times]$ & $1.058$ GB &  $[  12.6 \times]$ \\
$ 7$ & $  0.5 \cdot 10^{6}$ & $ 14.8$ GB & $  2.8$ GB & $[  5.3 \times]$ & $0.809$ GB &  $[  18.3 \times]$ \\
$ 8$ & $  0.3 \cdot 10^{6}$ & $ 12.5$ GB & $  2.2$ GB & $[  5.7 \times]$ & $0.486$ GB &  $[  25.6 \times]$ \\
$ 9$ & $  0.2 \cdot 10^{6}$ & $ 13.4$ GB & $  2.2$ GB & $[  6.0 \times]$ & $0.386$ GB &  $[  34.7 \times]$ \\
$10$ & $  0.2 \cdot 10^{6}$ & $ 12.1$ GB & $  1.9$ GB & $[  6.2 \times]$ & $0.264$ GB &  $[  45.8 \times]$ \\
\hline
\end{tabular}
  \caption{Estimated memory requirements as defined in Eq. \eqref{eqn:memory_convection} for different implementations for the solution of the convection dominated problem in Section \ref{sec:results_convection} with $n_\text{restart}^{(\text{outer})}=n_{\text{restart}}^{(\text{block})}=100$. The savings relative to the matrix-explicit implementation are shown in square brackets.}
\label{tab:memory_convection}
\end{center}
\end{table}
The memory requirements and savings in the (partially-) matrix-free method for our setup (see Tab. \ref{tab:problemsizes}) in $d=3$ dimensions are shown in Tab. \ref{tab:memory_convection}. As in Tab. \ref{tab:memory_diffusion}, the relative saving in memory consumption of the partially matrix-free method is asymptotically $\approx 7\times$. Due to the large number of temporary vectors in the GMRES method, the saving of the fully matrix-free method is less pronounced, but still nearly two orders of magnitude for the highest polynomial degrees. In hindsight and looking at Figs. \ref{fig:convection_dominated_iterations} and \ref{fig:convection_dominated_blockiter}, the chosen restart values were probably too conservative and memory consumption could have been reduced without impact on results (except for the $\vec{b}=(1.0,0.5,0.3)$ case, where the matrix-free method is less efficient) by choosing smaller values of $n_{\text{restart}}^{(\text{outer})}=15$ and $n_{\text{restart}}^{(\text{block})}=12$. The corresponding numbers in this case are shown in Tab. \ref{tab:memory_convection_small_nrestart} and show potential savings for the matrix-free method which are only a factor of $4\times$ smaller than those observed for the diffusion problem in Tab. \ref{tab:memory_diffusion}.
\begin{table}
\begin{center}
\begin{tabular}{rr|r|rr|rr}
\hline
degree $p$ & \# dofs & \multicolumn{5}{c}{memory requirements $M_{\text{total}}^{(\text{convection})}$} \\
           &         & \multicolumn{1}{|c|}{MX} & \multicolumn{2}{|c|}{PMF} & \multicolumn{2}{|c}{MF} \\
\hline\hline
$ 1$ & $ 33.6 \cdot 10^{6}$ & $ 22.5$ GB & $ 10.5$ GB & $[  2.1 \times]$ & $9.250$ GB &  $[   2.4 \times]$ \\
$ 2$ & $  9.5 \cdot 10^{6}$ & $ 15.8$ GB & $  4.3$ GB & $[  3.7 \times]$ & $2.614$ GB &  $[   6.0 \times]$ \\
$ 3$ & $  4.2 \cdot 10^{6}$ & $ 15.1$ GB & $  3.1$ GB & $[  4.9 \times]$ & $1.156$ GB &  $[  13.0 \times]$ \\
$ 4$ & $  2.0 \cdot 10^{6}$ & $ 13.5$ GB & $  2.4$ GB & $[  5.7 \times]$ & $0.551$ GB &  $[  24.6 \times]$ \\
$ 5$ & $  1.2 \cdot 10^{6}$ & $ 13.7$ GB & $  2.2$ GB & $[  6.2 \times]$ & $0.327$ GB &  $[  41.8 \times]$ \\
$ 6$ & $  0.7 \cdot 10^{6}$ & $ 12.4$ GB & $  1.9$ GB & $[  6.5 \times]$ & $0.189$ GB &  $[  65.8 \times]$ \\
$ 7$ & $  0.5 \cdot 10^{6}$ & $ 14.1$ GB & $  2.1$ GB & $[  6.6 \times]$ & $0.145$ GB &  $[  97.8 \times]$ \\
$ 8$ & $  0.3 \cdot 10^{6}$ & $ 12.1$ GB & $  1.8$ GB & $[  6.7 \times]$ & $0.087$ GB &  $[ 138.7 \times]$ \\
$ 9$ & $  0.2 \cdot 10^{6}$ & $ 13.1$ GB & $  1.9$ GB & $[  6.8 \times]$ & $0.069$ GB &  $[ 189.8 \times]$ \\
$10$ & $  0.2 \cdot 10^{6}$ & $ 11.9$ GB & $  1.7$ GB & $[  6.9 \times]$ & $0.047$ GB &  $[ 252.0 \times]$ \\
\hline
\end{tabular}
\caption{Estimated potential memory requirements as defined in Eq. \eqref{eqn:memory_convection} for different implementations for the solution of the convection dominated problem in Section \ref{sec:results_convection} with $n_\text{restart}^{(\text{outer})}=15$, $n_{\text{restart}}^{(\text{block})}=12$. The savings relative to the matrix-explicit implementation are shown in square brackets.}
\label{tab:memory_convection_small_nrestart}
\end{center}
\end{table}
\subsection{SPE10 test case}\label{sec:results_spe10}
We finally consider a challenging problem from a real-world
application. The SPE10 test case was used to compare different
reservoir modelling packages \cite{Christie2001}. Due to the large jumps in
the permeability field the solution of the
SPE10 problem is not smooth. At first sight, using a high-order method seems
to be questionable in this case since it it will not lead to the best
performance for a given error in the solution. However, in most
applications the solution is smooth almost everywhere in the domain and any
irregularities are confined to small, lower-dimensional subdomains. A
typical example are sharp drops in pressure across fronts in
atmospheric weather prediction. The main purpose of the following numerical
experiments on the SPE10 benchmark is therefore to demonstrate the
efficiency of our methods even for very challenging setups (in
particular large, high contrast jumps in the coefficients).

\begin{figure}
  \begin{center}
  \includegraphics[width=0.8\linewidth]{\figdir/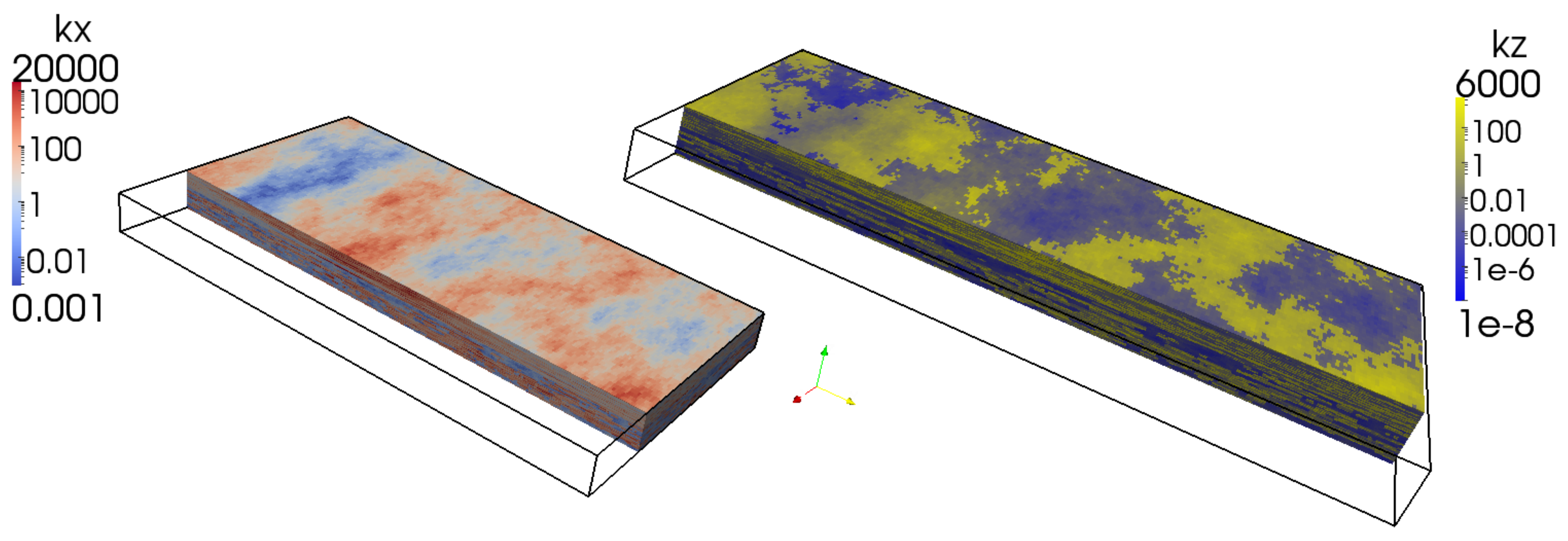}
  \caption{Permeability fields $K_x$ and $K_z$ from the SPE10 testcase}
  \label{fig:SPE10kx}
  \end{center}
\end{figure}

The second SPE10 dataset describes the permeability tensor
$K=\operatorname{diag}(K_x,K_y,K_z)$
of a three dimensional reservoir of size $L_x\times L_y\times L_z =
1,200\operatorname{ft}\times2,200\operatorname{ft}\times170 \operatorname{ft}$
which is divided into $60\times220\times85=1.122\cdot10^6$ cells of size
$20\operatorname{ft}\times10\operatorname{ft}\times2
\operatorname{ft}$. Permeability is assumed to be constant in each
cell, but can vary strongly over the domain (see
Fig. \ref{fig:SPE10kx} and Tab. \ref{tab:SPE10variation}).
\begin{table}
 \begin{center}
  \begin{tabular}{lrrr}
    \hline
    Field & minimum & mean & maximum\\
    \hline\hline
    $K_x$ & $6.65\cdot10^{-4}$ & $3.55\cdot 10^2$ & $2.0\cdot10^4$\\
    $K_y$ & $6.65\cdot10^{-4}$ & $3.55\cdot 10^2$ & $2.0\cdot10^4$\\
    $K_z$ & $6.65\cdot10^{-8}$ & $5.84\cdot 10^1$ & $6.0\cdot10^3$\\
    \hline
  \end{tabular}
  \caption{minimum, mean and maximum value of the three diagonal entries of permeability field for the SPE10 testcase}
  \label{tab:SPE10variation}
 \end{center}
\end{table}
Transport of a fluid with viscosity $\mu$ through the reservoir is
described by Darcy's law which relates pressure $u$ and flux $\vec{q}$ according to
\begin{equation}
  K \nabla u = -\mu \vec{q}.\label{eqn:Darcy}
\end{equation}
For an incompressible fluid the stationary continuity equation in the
absence of source terms is given by $\nabla\cdot \vec{q}=0$. Inserting this into
\eqref{eqn:Darcy} and assuming constant viscosity results in the
elliptic problem
\begin{equation}
  -\nabla\cdot (K\nabla u) = 0.
\end{equation}
This diffusion equation is solved with homogeneous Neumann boundary
conditions $\vec{\nu}\cdot \nabla u=0$ at $z=0$ and $z=L_z$ and Dirichlet boundary
conditions $u=-y$ on all other surfaces.  The problem is challenging
both due to the rapid variations in the diffusion tensor $K$ and the
strong anisotropy introduced by the large aspect ratio of the cells.

As for the diffusion equations in Section \ref{sec:results_diffusion},
the hybrid multigrid algorithm described in Section
\ref{sec:hybr-mult-prec} is used as a preconditioner. However, in
contrast to the conforming $Q_1$ space used there, for the SPE10
problem the low-order subspace is the piecewise constant $P_0$
space. The relative norm of the residual is not guaranteed to decrease
monotonically in the CG iteration, and (particularly for the
ill-conditioned matrix in this problem) this norm can be a poor
indicator of convergence. We therefore use the estimator for the
energy norm $||\vec{e}||_A=\sqrt{\vec{e}^T A \vec{e}}$ of the error
$\vecn{e}=\vecn{u}-\vecn{u}_{\text{exact}}$ given in \cite[Section~3.1]{Strakos2005}. This
does not add any significant overhead since only a small number of
additional scalar products are required in every iteration. The
quality of this estimate depends on the ``error-bandwidth'', i.e. the
number of terms that are used in the sum (3.5) in
\cite{Strakos2005}. For all results obtained here we used an
error-bandwidth of 4; increasing this further to 8 does not have any
discernible impact on the results.

\begin{figure}
  \begin{minipage}{0.5\linewidth}
    \includegraphics[width=\linewidth]{\figdir/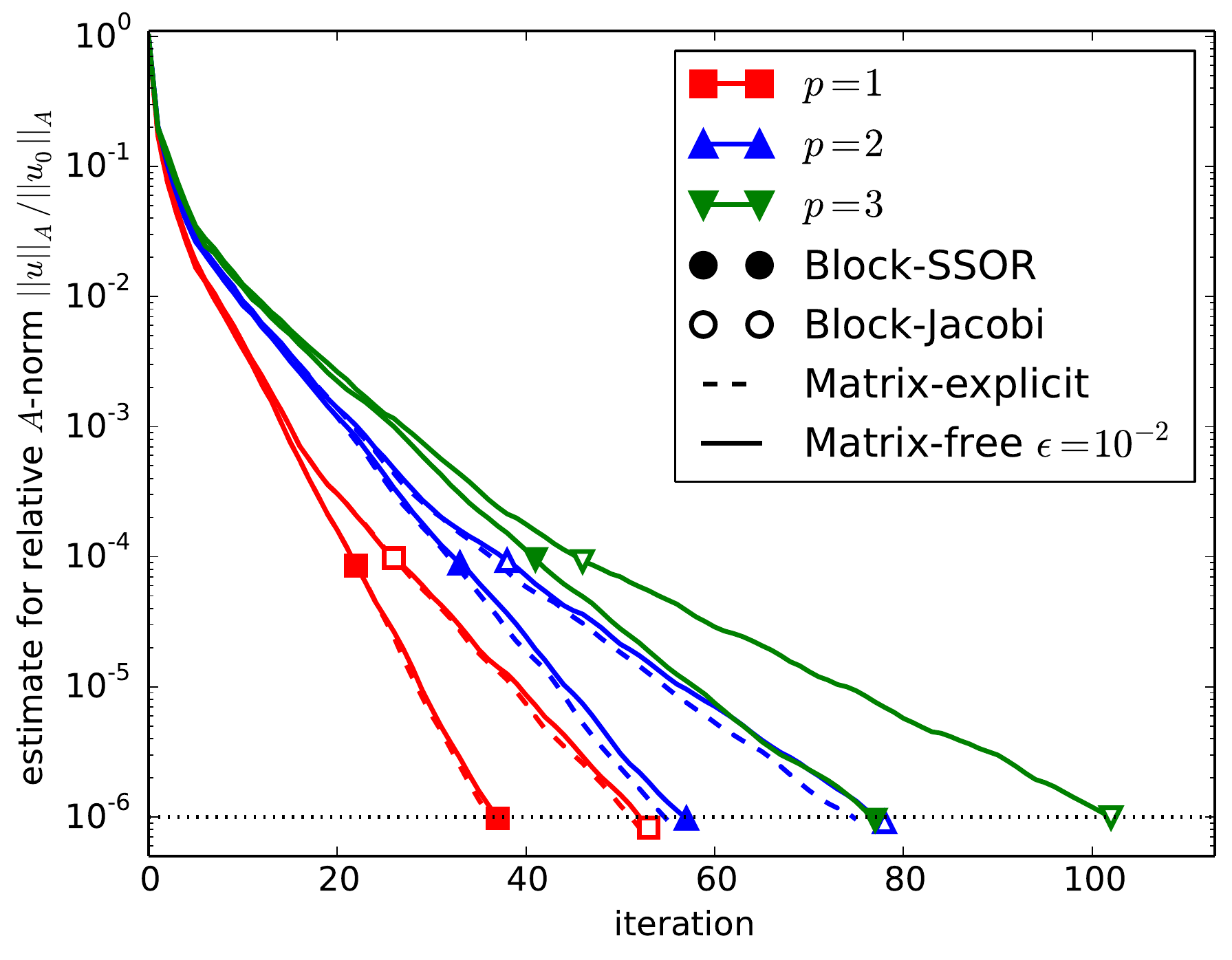}
  \end{minipage}
  \hfill
  \begin{minipage}{0.5\linewidth}
    \includegraphics[width=\linewidth]{\figdir/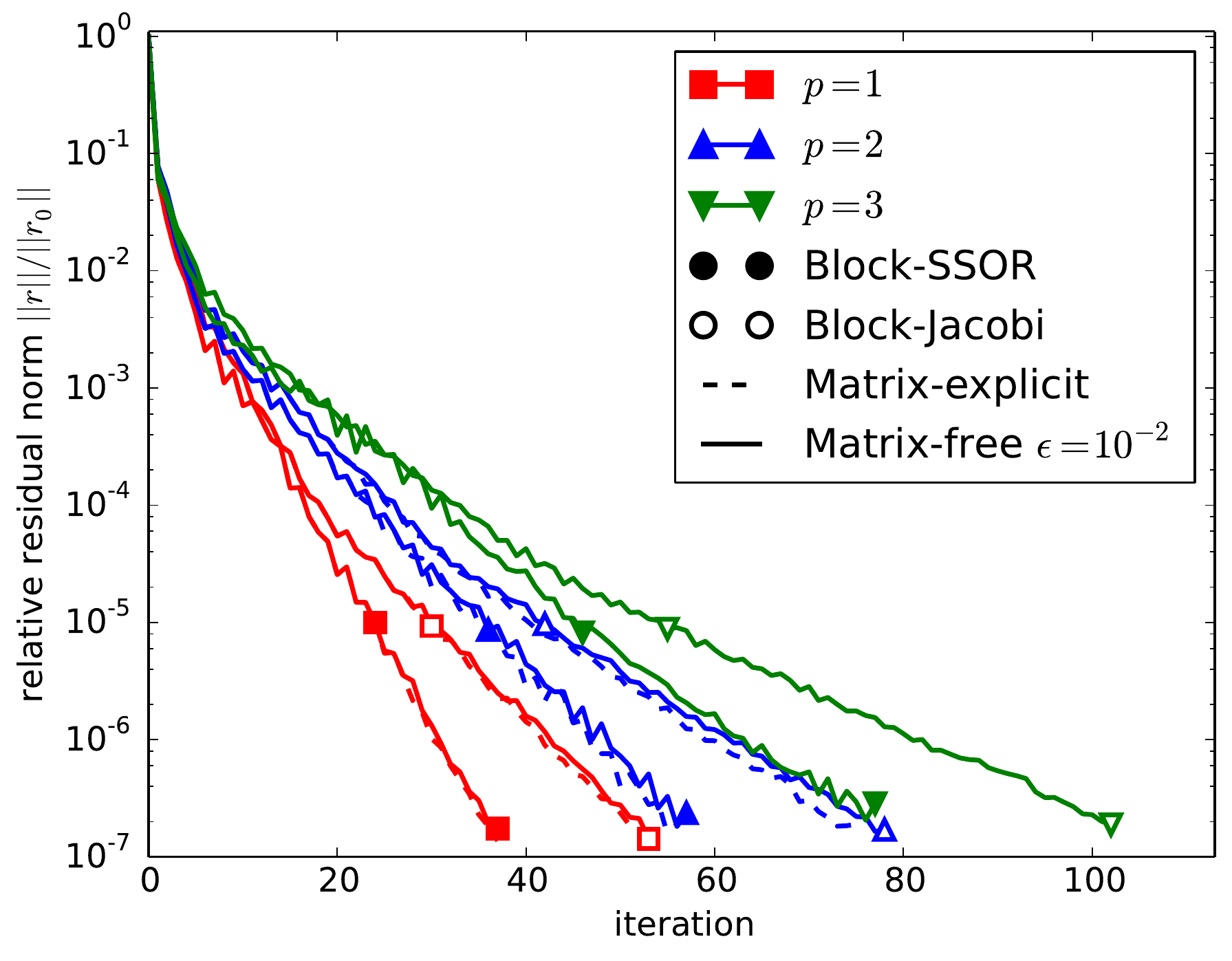}
  \end{minipage}
  \caption{Convergence history for SPE10 benchmark. Both the relative
    energy norm (left) and the relative residual norm (right) are
    shown for polynomial degrees 1 (red squares), 2 (blue upward
    triangles) and 3 (green downward triangles). Results for the
    block-SSOR smoother are marked by filled symbols and results for
    the block-Jacobi smoother by empty symbols.}
  \label{fig:convergence_spe10}
\end{figure}

The SPE10 problem is solved to a relative tolerance of $10^{-6}$ in
the energy norm. As Fig. \ref{fig:convergence_spe10} demonstrates, the
error is reduced monotonically (left), whereas there are significant
jumps in the relative residual (right). Solving to tighter tolerances
is not of practical interest since other sources of error dominate the
physical problem. The number of iterations and solution times for
different polynomial degrees on one node of the ``donkey'' cluster are
shown in Tab. \ref{tab:SPE10results_donkey} (run-ahead iterations for
the initial estimate of the error in the energy norm are not counted
but included in the reported runtimes). In the matrix-free solver the
diagonal blocks were inverted to a relative tolerance of $\epsilon=10^{-2}$.
\begin{table}
  \begin{center}
\begin{tabular}{lc|rrr|rrr}
\hline
 & degree & \multicolumn{3}{|c|}{\# iterations} & \multicolumn{3}{|c}{time per unknown $[\mu s]$}\\
smoother & & MX & MF & PMF & MX & MF & PMF\\\hline\hline
\multirow{3}{*}{jacobi} & 1  &  52 &  53 &  52 & 5.67 & 38.00 & 8.58\\
& 2  &  75 &  78 &  74 & 23.80 & 24.15 & 6.99\\
& 3  & --- & 102 & --- & --- & 22.25 & ---\\
\hline
\multirow{3}{*}{ssor} & 1  &  37 &  37 &  37 & 7.44 & 57.67 & 15.84\\
& 2  &  55 &  57 &  55 & 34.00 & 36.09 & 11.56\\
& 3  & --- &  77 & --- & --- & 34.14 & ---\\
\hline
\end{tabular}
  \caption{SPE10 benchmark on one node of ``donkey''. The number of
    iterations and total solution time per unknown is shown for
    different polynomial degrees. For the highest polynomial degree
    only matrix-free solver results are reported due to memory
    limitations.}\label{tab:SPE10results_donkey}
  \end{center}
\end{table}
We find that in all cases a block-Jacobi smoother is more efficient
than block-SSOR. For the lowest polynomial degree ($p=1$) the
matrix-explicit solver results in the shortest solution time, although
it is only around $1.5\times$ faster than the partially matrix-free
solver. However, for $p=2$ the partially matrix-free solver is $3\times$
faster than the matrix-explicit version which has about the same
performance as the matrix-free implementation.  For the highest
polynomial degree ($p=3$) neither the full matrix nor the
block-diagonal fit into memory and only results for the matrix-free
solver are reported.
\section{Conclusion}\label{sec:Conclusion}
We described the efficient implementation of matrix-free multigrid block
smoothers for higher order DG discretisations of the
convection-diffusion equation. Since the implementation is FLOP-bound,
this leads to significantly better utilisation of modern manycore
CPUs. Sum-factorisation techniques reduce the computational complexity
from $\mathcal{O}(p^{2d})$ to $\mathcal{O}(d\cdot p^{d+1})$ and make the
method particularly competitive for high polynomial degrees $p$. Our
numerical experiments confirm that this leads to a significant
reduction of the solution time even for relatively low polynomial
degrees. Since the matrix does not have to be assembled and stored
explicitly, our approach also leads to significantly reduced storage
requirements. For moderately high polynomial degrees the totally
matrix-free implementation typically gives the best overall
performance. For low polynomial degrees we find that a partially
matrix-free implementation, in which only the operator application is
matrix-free, can reduce the runtime even further. The construction of
a suitable preconditioner for the matrix-free block-inversion in
convection-dominated problems is challenging if the advection field is
not axis aligned, and therefore in this setup the partially
matrix-free implementation (which avoids this problem) is particularly
efficient.

With the SPE10 benchmark we have included a challenging, real world
application. To quantify the performance for other realistic setups we
will also apply the matrix-free solvers to the pressure correction equation
in atmospheric models in the future.

Since the local blocks are only inverted in the preconditioner,
additional approximations could improve performance further. For
example, the use of lower-order Gauss-Lobatto quadrature rules for
calculation of the volume integral can reduce the number of function
evaluations since the quadrature points coincide with those of the
surface integral.

The implementation of the highly-optimised matrix-free operators described
in this work required a significant amount of low-level optimisation, leading
to sophisticated code which is hard to maintain. In the future, we will
explore the use of automatic code generation techniques which are currently
integrated into the DUNE framework \cite{Kempf2017} and have been shown to
give promising results.
\section*{Acknowledgments}
We would like to thank Lawrence Mitchell (Imperial College London) for
useful discussions about alternative preconditioning methods based on
low-order discretisation on a nodal subgrid. Similarly, we are
grateful to Rob Scheichl (Bath) for helpful comments throughout the
project. The authors acknowledge support by the state of
Baden-W\"{u}rttemberg through bwHPC. Part of the work was carried out during
extended visits of one of the authors to Heidelberg University, funded
by MATCH (MAThematics Center Heidelberg) and a startup grant from the
University of Bath. This project was supported in part by Deutsche
Forschungsgemeinschaft under grant Ba 1498/10-2 within the SPPEXA programme.

\appendix
\section{Discretisation in the lower- order subspaces}\label{sec:lower_order_discretisation}
In this appendix we show how the coarse grid matrix $\hat{A}=P^TAP$ can be calculated efficiently by re-discretisation in the subspace $\hat{V}_h$. As discussed in Section \ref{sec:low-order-subspace} this avoids the expensive calculation via the explicit Galerkin product.
\subsection{Piecewise constant subspace $P_0$}\label{sec:low_order_matrix_P0}
The $P_0$ basis functions have a vanishing gradient on each element.
This simplifies the evaluation of \(a_h(\hat{\psi}_j, \hat{\psi}_i)\) since
the diffusion and convection term do not contribute to the volume
integral $a_h^{\textup{v}}(\hat{\psi}_j,\hat{\psi}_j)$ and the consistency
and symmetrisation term vanish in the face integral
$a_h^{\textup{if}}(\hat{\psi}_j,\hat{\psi}_j)$.  The entries of the Galerkin
product with $P_0$ as a subspace are
\begin{equation}
\begin{split}
  \left(P^T A P\right)_{ij} &= \sum_{T\in\mathcal{T}_h} \int_T c
  \hat{\psi}_j \hat{\psi}_i \, dx
  +\sum_{F\in\mathcal{F}_h^{i}} \int_F
  \Phi(\hat{\psi}_j^-,\hat{\psi}_i^+,\vec{\nu}_F\cdot \vec{b})\llbracket
  \hat{\psi}_i \rrbracket \, ds
  +\sum_{F\in\mathcal{F}_h^{D}}
  \int_F  \Phi(\hat{\psi}_j,0,\vec{\nu}_F\cdot \vec{b})\llbracket \hat{\psi}_i
  \rrbracket \, ds \\
  &+ \alpha p (p+d-1) \left(\sum_{F\in\mathcal{F}_h^{i}} \int_F
  \frac{\langle \delta_{K\nu_F}^- , \delta_{K\nu_F}^+ \rangle |F|}{\min(|T^-(F)|,|T^+(F)|)}
  \llbracket \hat{\psi}_j \rrbracket \llbracket \hat{\psi}_i
  \rrbracket  \, ds
  +
  \sum_{F\in\mathcal{F}_h^{D}} \int_F
  \frac{\delta_{K\nu_F}^-|F|}{|T^-(F)|}
  \hat{\psi}_j \hat{\psi}_i\, ds \right).
\end{split}
\label{eq:p0_galerkin1}
\end{equation}
This should be compared to the finite volume discretisation of
\eqref{eqn:linearproblem} and \eqref{eqn:boundary_conditions}
which can be obtained as follows: let \(\vec{x}(T)\) be the center of cell
\(T\), \(\vec{x}(F)\) the center of face \(F\) and \(\vec{x}^-(F)\) the center of
inside cell, \(\vec{x}^+(F)\) the center of outside cell with respect to
\(F\). The evaluation of the convection and reaction term at \(\vec{x}(F)\)
and \(\vec{x}(T)\), respectively, is denoted in shorthand by \(\hat{\vec{b}} =
\vec{b}(\vec{x}(F)\)), \(\hat{c} = c(\vec{x}(T))\).  We approximate the diffusive normal
flux with a central difference approximation and use the harmonic
average of the diffusion tensor on \(F\in\mathcal{F}_h^i\). On boundary
faces \(F\in\mathcal{F}_h^b\) the diffusive normal flux is approximated
by a one-sided difference quotient. This results in the finite volume
operator
\begin{equation}
\begin{split}
  \hat{A}_{ij}^{\textup{FV}} &= \sum_{T\in\mathcal{T}_h}
  \hat{c}\,\hat{\psi}_j\hat{\psi}_i |T|
  +\sum_{F\in\mathcal{F}_h^{i}}
  \Phi(\hat{\psi}_j^-,\hat{\psi}_i^+,\vec{\nu}_F\cdot \hat{\vec{b}})
  \llbracket \hat{\psi}_i \rrbracket |F|
  +\sum_{F\in\mathcal{F}_h^{D}}
  \Phi(\hat{\psi}_j,0,\vec{\nu}_F\cdot \hat{\vec{b}})
  \llbracket \hat{\psi}_i \rrbracket |F| \\
  &+ \sum_{F\in\mathcal{F}_h^{i}}
  \frac{\langle \delta_{K\nu_F}^- , \delta_{K\nu_F}^+ \rangle}
  {\left\| \vec{x}^-(F) - \vec{x}^+(F)\right\|_2}
  \llbracket \hat{\psi}_j \rrbracket
  \llbracket \hat{\psi}_i \rrbracket |F|
  + \sum_{F\in\mathcal{F}_h^{D}}
  \frac{\delta_{K\nu_F}}{\left\|\vec{x}^-(F) - \vec{x}(F)\right\|_2}
  \hat{\psi}_j \hat{\psi}_i |F| .
\end{split}
\label{eq:p0_matrix}
\end{equation}
If in $P^T A P$ the upwind flux and $\int_T c\;ds$ are approximated by
the box-rule (or $c$ and $\vec{b}$ are piecewise constant on the grid
cells), the first three terms in \eqref{eq:p0_galerkin1} and
\eqref{eq:p0_matrix} are identical. This is an admissible
approximation since it is used in a subspace correction. To relate the
remaining two integrals over interior- and boundary-faces, observe
that on an equidistant mesh \(|T^-(F)| = |T^+(F)| \;\forall
F\in\mathcal{F}_h^i\). If the grid is also axi-parallel, the inverse
distances in \eqref{eq:p0_matrix} can be written as
\begin{equation}
\frac{1}{\left\| \vec{x}^-(F) - \vec{x}^+(F)\right\|_2} =
  \frac{|F|}{\min(|T^-(F)|,|T^+(F)|)} \quad \forall F\in\mathcal{F}_h^i
\end{equation}
on the interior faces and as
\begin{equation}
  \frac{1}{\left\|\vec{x}^-(F) - \vec{x}(F)\right\|_2} =
  \frac{|F|}{2|T^-(F)|}
  \quad\forall F\in\mathcal{F}_h^b
\end{equation}
on the boundary faces. Provided the diffusivity $K$ is also piecewise
constant on the grid cells (or, on each face is approximated by its
mean value), the final two terms in \eqref{eq:p0_galerkin1} and
\eqref{eq:p0_matrix} differ by multiplicative factors of $\alpha p(p+d-1)$
and $2\alpha p(p+d-1)$ respectively. Hence, instead of constructing
$\hat{A}$ via the Galerkin product $\hat{A}=P^T AP$, it can be
obtained from a modified finite-volume discretisation, if the last two
terms in \eqref{eq:p0_matrix} are scaled appropriately.
\subsection{Conforming piecewise linear subspace $Q_1$}\label{sec:low_order_matrix_Q1}
The \(Q_1\) conforming piecewise linear basis functions are
continuous over the interior edges. Exploiting this fact for
\(a_h(\hat{\psi}_j, \hat{\psi}_i)\) implies that all terms containing
the jump operator \(\llbracket \cdot \rrbracket\) that are summed over
\(\mathcal{F}_h^i\) vanish. The entries of the Galerkin product with \(
Q_1\) as a subspace read
\begin{equation}
\begin{split}
  \left(P^T A P\right)_{ij} &= \sum_{T\in\mathcal{T}_h} \int_T \left((K\nabla\hat{\psi}_j
  - \vec{b}\hat{\psi}_j)\cdot \nabla \hat{\psi}_i
  +\hat{\psi}_j \hat{\psi}_i \right) dx\\
  &+\sum_{F\in\mathcal{F}_h^{D}} \int_F\left(
\Phi(\hat{\psi}_j,0,\vec{\nu}_F\cdot \vec{b}) \hat{\psi}_i
  -\vec{\nu}_F\cdot (\omega^- K\nabla \hat{\psi}_j) \,\hat{\psi}_i
  - \hat{\psi}_j \vec{\nu}_F\cdot (\omega^- K \nabla \hat{\psi}_i)
  + \gamma_{F} \hat{\psi}_j \hat{\psi}_i \right) ds .
\end{split}
\end{equation}
This is identical to the conforming piecewise linear Finite Element
operator for the model problem \eqref{eqn:elliptic_eqn}. The only
differences are the integrals over the Dirichlet boundary which is a
remainder from the DG bilinear form. In contrast to the \(P_0\)
subspace \(P^T A P\) is rediscretised by calculating \(a_h(\hat{\psi}_j,
\hat{\psi}_i) \, \forall i,j\) with the optimisation that the integral for all
\(F\in\mathcal{F}_h^i\) is skipped. Note that although $\hat{A}_{ij}$ is
calculated by assembling the DG operator $a_h(\hat{\psi}_i,\hat{\psi}_j)$ on
a low-dimensional subspace, the polynomial degree which is used in the
calculation of $\gamma_F$ is the degree $p$ of the DG space $V_h^p$, and
not the degree of $Q_1$.
\section{Setup costs}\label{sec:setup_costs}
The setup costs for the diffusion problems described in Section
\ref{sec:results_diffusion} are shown for a selection of polynomial
degrees in in Tables \ref{tab:setup_poisson} and
\ref{tab:setup_inhomogeneous}. Note that the (partially) matrix-free
implementation has a significantly lower setup cost since it does not
require assembly of the DG matrix and construction of the coarse grid
matrix $\hat{A}$ via the Galerkin product.
\begin{table}
\begin{center}
\begin{tabular}{|l|rrr|rrr|rrr|}
\hline
 & \multicolumn{3}{|c|}{$p=2$} & \multicolumn{3}{|c|}{$p=3$} & \multicolumn{3}{|c|}{$p=6$}\\
 & MX & MF & PMF & MX & MF & PMF & MX & MF & PMF\\
\hline\hline
Assemble prolongation  &  0.09 &  0.07 &  0.07 &  0.05 &  0.04 &  0.04 &  0.02 &  0.01 &  0.01\\
DG Matrix setup+assembly  &  2.36 & --- & --- &  2.50 & --- & --- &  3.85 & --- & ---\\
DG smoother setup (block-factor)  &  0.25 & --- &  1.10 &  0.26 & --- &  1.09 &  0.43 & --- &  1.64\\
Diagonal matrix assembly  & --- &  0.91 & --- & --- &  0.91 & --- & --- &  1.43 & ---\\
Galerkin product $\hat{A}=P^TAP$  &  2.01 & --- & --- &  0.99 & --- & --- &  2.89 & --- & ---\\
Coarse matrix assembly  & --- &  0.18 &  0.18 & --- &  0.04 &  0.04 & --- &  0.00 &  0.00\\
AMG setup  &  1.03 &  0.37 &  0.37 &  0.22 &  0.10 &  0.10 &  0.02 &  0.01 &  0.01\\
\hline
total setup  &  5.73 &  1.54 &  1.73 &  4.03 &  1.09 &  1.27 &  7.21 &  1.46 &  1.67\\
\hline
total solve time & 33.13 & 25.85 & 12.64 & 34.25 &  9.71 &  7.35 & 69.61 &  4.71 &  6.41\\
\hline
\end{tabular}
\caption{Breakdown of setup time for the Poisson problem for different polynomial degrees $p$. To put the setup costs into perspective, the total solve time (which includes the setup time) is also given}
\label{tab:setup_poisson}
\end{center}
\end{table}

\begin{table}
\begin{center}
\begin{tabular}{|l|rrr|rrr|rrr|}
\hline
 & \multicolumn{3}{|c|}{$p=2$} & \multicolumn{3}{|c|}{$p=3$} & \multicolumn{3}{|c|}{$p=6$}\\
 & MX & MF & PMF & MX & MF & PMF & MX & MF & PMF\\
\hline\hline
Assemble prolongation  &  0.09 &  0.07 &  0.07 &  0.05 &  0.04 &  0.04 &  0.02 &  0.01 &  0.01\\
DG Matrix setup+assembly  &  3.24 & --- & --- &  3.25 & --- & --- &  4.52 & --- & ---\\
DG smoother setup (block-factor)  &  0.25 & --- &  1.28 &  0.26 & --- &  1.23 &  0.43 & --- &  1.84\\
Diagonal matrix assembly  & --- &  1.07 & --- & --- &  1.03 & --- & --- &  1.67 & ---\\
Galerkin product $\hat{A}=P^TAP$  &  1.97 & --- & --- &  0.98 & --- & --- &  2.86 & --- & ---\\
Coarse matrix assembly  & --- &  0.20 &  0.19 & --- &  0.04 &  0.04 & --- &  0.00 &  0.00\\
AMG setup  &  1.05 &  0.39 &  0.39 &  0.22 &  0.10 &  0.10 &  0.02 &  0.01 &  0.01\\
\hline
total setup  &  6.60 &  1.73 &  1.93 &  4.77 &  1.21 &  1.42 &  7.85 &  1.70 &  1.87\\
\hline
total solve time & 31.70 & 32.74 & 13.82 & 30.95 & 10.89 &  7.51 & 62.58 &  5.52 &  6.37\\
\hline
\end{tabular}
\caption{Breakdown of setup time for the inhomogeneous diffusion
  problem for different polynomial degrees $p$. To put the setup costs
  into perspective, the total solve time (which includes the setup
  time) is also given}
\label{tab:setup_inhomogeneous}
\end{center}
\end{table}
\clearpage
\bibliographystyle{elsarticle-num}

\begin{thebibliography}{10}
\expandafter\ifx\csname url\endcsname\relax
  \def\url#1{\texttt{#1}}\fi
\expandafter\ifx\csname urlprefix\endcsname\relax\def\urlprefix{URL }\fi
\expandafter\ifx\csname href\endcsname\relax
  \def\href#1#2{#2} \def\path#1{#1}\fi

\bibitem{Wheeler1978}
M.~F. Wheeler, {An elliptic collocation-finite element method with interior
  penalties}, SIAM Journal on Numerical Analysis 15~(1) (1978) 152--161.
\newblock \href {http://dx.doi.org/10.1137/0715010}
  {\path{doi:10.1137/0715010}}.

\bibitem{Arnold1982}
D.~N. Arnold, {An interior penalty finite element method with discontinuous
  elements}, SIAM journal on numerical analysis 19~(4) (1982) 742--760.
\newblock \href {http://dx.doi.org/10.1137/0719052}
  {\path{doi:10.1137/0719052}}.

\bibitem{Oden1998}
J.~Oden, I.~Babu\^{s}ka, C.~E. Baumann, A discontinuous $hp-$ finite element
  method for diffusion problems, Journal of Computational Physics 146~(2)
  (1998) 491 -- 519.
\newblock \href {http://dx.doi.org/10.1006/jcph.1998.6032}
  {\path{doi:10.1006/jcph.1998.6032}}.

\bibitem{Baumann1999}
C.~E. Baumann, J.~T. Oden, {A discontinuous $hp-$ finite element method for
  convection diffusion problems}, Computer Methods in Applied Mechanics and
  Engineering 175~(3) (1999) 311--341.
\newblock \href {http://dx.doi.org/10.1016/S0045-7825(98)00359-4}
  {\path{doi:10.1016/S0045-7825(98)00359-4}}.

\bibitem{Riviere1999}
B.~Rivi{\`e}re, M.~F. Wheeler, V.~Girault, {Improved energy estimates for
  interior penalty, constrained and discontinuous Galerkin methods for elliptic
  problems. Part I}, Computational Geosciences 3~(3-4) (1999) 337--360.
\newblock \href {http://dx.doi.org/10.1023/A:1011591328604}
  {\path{doi:10.1023/A:1011591328604}}.

\bibitem{Arnold2002}
D.~N. Arnold, F.~Brezzi, B.~Cockburn, L.~D. Marini, {Unified analysis of
  discontinuous Galerkin methods for elliptic problems}, SIAM journal on
  numerical analysis 39~(5) (2002) 1749--1779.
\newblock \href {http://dx.doi.org/10.1137/S0036142901384162}
  {\path{doi:10.1137/S0036142901384162}}.

\bibitem{Riviere2008}
B.~Riviere, {Discontinuous Galerkin Methods for Solving Elliptic and Parabolic
  Equations: Theory and Implementation}, Frontiers in Applied Mathematics,
  Society for Industrial and Applied Mathematics (SIAM, 3600 Market Street,
  Floor 6, Philadelphia, PA 19104), 2008.
\newblock \href {http://dx.doi.org/10.1137/1.9780898717440}
  {\path{doi:10.1137/1.9780898717440}}.

\bibitem{DiPietro2011}
D.~Di~Pietro, A.~Ern, {Mathematical Aspects of Discontinuous Galerkin Methods},
  Math{\'e}matiques et Applications, Springer Berlin Heidelberg, 2011.
\newblock \href {http://dx.doi.org/10.1007/978-3-642-22980-0}
  {\path{doi:10.1007/978-3-642-22980-0}}.

\bibitem{Cockburn2002}
B.~Cockburn, C.~Dawson, Approximation of the velocity by coupling discontinuous
  galerkin and mixed finite element methods for flow problems, Computational
  Geosciences 6~(3) (2002) 505--522.
\newblock \href {http://dx.doi.org/10.1023/A:1021203618109}
  {\path{doi:10.1023/A:1021203618109}}.

\bibitem{Bastian2003}
P.~Bastian, Higher Order Discontinuous Galerkin Methods for Flow and Transport
  in Porous Media, Springer Berlin Heidelberg, Berlin, Heidelberg, 2003, pp.
  1--22.
\newblock \href {http://dx.doi.org/10.1007/978-3-642-19014-8\_1}
  {\path{doi:10.1007/978-3-642-19014-8\_1}}.

\bibitem{Bastian2014II}
P.~Bastian, {A fully-coupled discontinuous Galerkin method for two-phase flow
  in porous media with discontinuous capillary pressure}, Computational
  Geosciences 18~(5) (2014) 779--796.
\newblock \href {http://dx.doi.org/10.1007/s10596-014-9426-y}
  {\path{doi:10.1007/s10596-014-9426-y}}.

\bibitem{Siefert2014}
C.~Siefert, R.~Tuminaro, A.~Gerstenberger, G.~Scovazzi, S.~S. Collis,
  \href{https://doi.org/10.1007/s10596-014-9419-x}{Algebraic multigrid
  techniques for discontinuous galerkin methods with varying polynomial order},
  Computational Geosciences 18~(5) (2014) 597--612.
\newblock \href {http://dx.doi.org/10.1007/s10596-014-9419-x}
  {\path{doi:10.1007/s10596-014-9419-x}}.
\newline\urlprefix\url{https://doi.org/10.1007/s10596-014-9419-x}

\bibitem{Chorin1968}
A.~J. Chorin, {Numerical solution of the Navier-Stokes equations}, Mathematics
  of computation 22~(104) (1968) 745--762.
\newblock \href {http://dx.doi.org/10.1090/S0025-5718-1968-0242392-2}
  {\path{doi:10.1090/S0025-5718-1968-0242392-2}}.

\bibitem{Temam1969}
R.~T{\'e}mam, Sur l'approximation de la solution des {\'e}quations de
  navier-stokes par la m{\'e}thode des pas fractionnaires (ii), Archive for
  Rational Mechanics and Analysis 33~(5) (1969) 377--385.
\newblock \href {http://dx.doi.org/10.1007/BF00247696}
  {\path{doi:10.1007/BF00247696}}.

\bibitem{2016arXiv161200657P}
M.~Piatkowski, S.~M{\"u}thing, P.~Bastian,
  \href{https://www.sciencedirect.com/science/article/pii/S0021999117308732}{A
  stable and high-order accurate discontinuous galerkin based splitting method
  for the incompressible navier–stokes equations}, Journal of Computational
  Physics 356 (2018) 220 -- 239.
\newblock \href {http://dx.doi.org/https://doi.org/10.1016/j.jcp.2017.11.035}
  {\path{doi:https://doi.org/10.1016/j.jcp.2017.11.035}}.
\newline\urlprefix\url{https://www.sciencedirect.com/science/article/pii/S0021999117308732}

\bibitem{Restelli2009}
M.~Restelli, F.~X. Giraldo, {A conservative discontinuous Galerkin
  semi-implicit formulation for the Navier--Stokes equations in nonhydrostatic
  mesoscale modeling}, SIAM Journal on Scientific Computing 31~(3) (2009)
  2231--2257.
\newblock \href {http://dx.doi.org/10.1137/070708470}
  {\path{doi:10.1137/070708470}}.

\bibitem{Dedner2015}
A.~Dedner, R.~Kl{\"o}fkorn, {On Efficient Time Stepping using the Discontinuous
  Galerkin Method for Numerical Weather Prediction}, in: Parallel Computing: On
  the Road to Exascale, Advances in Parallel Computing, 2015, pp. 627--636.
\newblock \href {http://dx.doi.org/10.3233/978-1-61499-621-7-627}
  {\path{doi:10.3233/978-1-61499-621-7-627}}.

\bibitem{Trottenberg2000}
U.~Trottenberg, C.~W. Oosterlee, A.~Schuller, Multigrid, Academic Press,
  Cambridge, Massachusetts, 2000.

\bibitem{McCalpin1995}
J.~D. McCalpin, {Memory Bandwidth and Machine Balance in Current High
  Performance Computers}, IEEE Computer Society Technical Committee on Computer
  Architecture (TCCA) Newsletter (1995) 19--25.

\bibitem{Orszag1980}
S.~A. Orszag, {Spectral methods for problems in complex geometries}, Journal of
  Computational Physics 37~(1) (1980) 70--92.
\newblock \href {http://dx.doi.org/10.1016/0021-9991(80)90005-4}
  {\path{doi:10.1016/0021-9991(80)90005-4}}.

\bibitem{Karniadakis2005}
G.~Karniadakis, S.~Sherwin, Spectral/hp Element Methods for {CFD}, Oxford
  University Press, 2005.

\bibitem{Vos2010}
P.~E. Vos, S.~J. Sherwin, R.~M. Kirby, From h to p efficiently: Implementing
  finite and spectral/hp element methods to achieve optimal performance for
  low-and high-order discretisations, Journal of Computational Physics 229~(13)
  (2010) 5161--5181.
\newblock \href {http://dx.doi.org/10.1016/j.jcp.2010.03.031}
  {\path{doi:10.1016/j.jcp.2010.03.031}}.

\bibitem{Kronbichler2012}
M.~Kronbichler, K.~Kormann, A generic interface for parallel cell-based finite
  element operator application, Computers \& Fluids 63 (2012) 135--147.
\newblock \href {http://dx.doi.org/10.1016/j.compfluid.2012.04.012}
  {\path{doi:10.1016/j.compfluid.2012.04.012}}.

\bibitem{Muething2017}
S.~M\"{u}thing, M.~Piatkowski, P.~Bastian, {High-performance Implementation of
  Matrix-free High-order Discontinuous Galerkin Methods}, Submitted to SIAM
  SISC\href {http://arxiv.org/abs/1711.10885} {\path{arXiv:1711.10885}}.

\bibitem{Brandt1984}
A.~Brandt, S.~F. McCormick, J.~W. Ruge, {Algebraic multigrid {(AMG)} for sparse
  matrix equations}, in: D.~J. Evans (Ed.), {Sparsity And Its Applications},
  1984, pp. 258--283.

\bibitem{Stueben2001}
K.~St\"{u}ben, A review of algebraic multigrid, Journal of Computational and
  Applied Mathematics 128~(1) (2001) 281 -- 309, numerical Analysis 2000. Vol.
  VII: Partial Differential Equations.
\newblock \href {http://dx.doi.org/10.1016/S0377-0427(00)00516-1}
  {\path{doi:10.1016/S0377-0427(00)00516-1}}.

\bibitem{Saad2003}
Y.~Saad, {Iterative Methods for Sparse Linear Systems}, 2nd Edition, Society
  for Industrial and Applied Mathematics, 2003.
\newblock \href {http://dx.doi.org/10.1137/1.9780898718003}
  {\path{doi:10.1137/1.9780898718003}}.

\bibitem{Bastian2012}
P.~Bastian, M.~Blatt, R.~Scheichl, {Algebraic multigrid for discontinuous
  Galerkin discretizations of heterogeneous elliptic problems}, Numerical
  Linear Algebra with Applications 19~(2) (2012) 367--388.
\newblock \href {http://dx.doi.org/10.1002/nla.1816}
  {\path{doi:10.1002/nla.1816}}.

\bibitem{Kronbichler2016}
M.~Kronbichler, W.~A. Wall, {A performance comparison of continuous and
  discontinuous Galerkin methods with multigrid solvers, including a new
  multigrid scheme for HDG }\href {http://arxiv.org/abs/1611.03029}
  {\path{arXiv:1611.03029}}.

\bibitem{dealII84}
W.~Bangerth, D.~Davydov, T.~Heister, L.~Heltai, G.~Kanschat, M.~Kronbichler,
  M.~Maier, B.~Turcksin, D.~Wells, {The \texttt{deal.II} Library, Version 8.4},
  Journal of Numerical Mathematics 24.
\newblock \href {http://dx.doi.org/10.1515/jnma-2016-1045}
  {\path{doi:10.1515/jnma-2016-1045}}.

\bibitem{Fischer1997}
P.~F. Fischer, An overlapping schwarz method for spectral element solution of
  the incompressible navier-stokes equations, Journal of Computational Physics
  133~(1) (1997) 84 -- 101.
\newblock \href {http://dx.doi.org/10.1006/jcph.1997.5651}
  {\path{doi:10.1006/jcph.1997.5651}}.

\bibitem{Brown2010}
J.~Brown, {Efficient Nonlinear Solvers for Nodal High-Order Finite Elements in
  {3D}}, Journal of Scientific Computing 45~(1-3) (2010) 48--63.
\newblock \href {http://dx.doi.org/10.1007/s10915-010-9396-8}
  {\path{doi:10.1007/s10915-010-9396-8}}.

\bibitem{Austin2012}
T.~M. Austin, M.~Brezina, B.~Jamroz, C.~Jhurani, T.~A. Manteuffel, J.~Ruge,
  {Semi-automatic sparse preconditioners for high-order finite element methods
  on non-uniform meshes}, Journal of Computational Physics 231~(14) (2012)
  4694--4708.
\newblock \href {http://dx.doi.org/10.1016/j.jcp.2012.03.013}
  {\path{doi:10.1016/j.jcp.2012.03.013}}.

\bibitem{Bastian2008b}
P.~Bastian, M.~Blatt, A.~Dedner, C.~Engwer, R.~Kl\"{o}fkorn, M.~Ohlberger,
  O.~Sander, {A generic grid interface for parallel and adaptive scientific
  computing. {Part I}: abstract framework}, Computing 82~(2-3) (2008) 103--119.
\newblock \href {http://dx.doi.org/10.1007/s00607-008-0003-x}
  {\path{doi:10.1007/s00607-008-0003-x}}.

\bibitem{Bastian2014}
P.~Bastian, C.~Engwer, D.~G{\"o}ddeke, O.~Iliev, O.~Ippisch, M.~Ohlberger,
  S.~Turek, J.~Fahlke, S.~Kaulmann, S.~M{\"u}thing, D.~Ribbrock, EXA-DUNE:
  Flexible PDE Solvers, Numerical Methods and Applications, Springer
  International Publishing, Cham, 2014, pp. 530--541.
\newblock \href {http://dx.doi.org/10.1007/978-3-319-14313-2\_45}
  {\path{doi:10.1007/978-3-319-14313-2\_45}}.

\bibitem{Christie2001}
M.~Christie, M.~Blunt, et~al., {Tenth SPE comparative solution project: A
  comparison of upscaling techniques}, in: SPE Reservoir Simulation Symposium,
  Society of Petroleum Engineers, 2001.
\newblock \href {http://dx.doi.org/10.2118/72469-PA}
  {\path{doi:10.2118/72469-PA}}.

\bibitem{Epshteyn2007}
Y.~Epshteyn, B.~Rivi{\`e}re, {Estimation of penalty parameters for symmetric
  interior penalty Galerkin methods}, Journal of Computational and Applied
  Mathematics 206~(2) (2007) 843--872.
\newblock \href {http://dx.doi.org/10.1016/j.cam.2006.08.029}
  {\path{doi:10.1016/j.cam.2006.08.029}}.

\bibitem{Hartmann2008}
R.~Hartmann, P.~Houston, {An optimal order interior penalty discontinuous
  Galerkin discretization of the compressible Navier--Stokes equations},
  Journal of Computational Physics 227~(22) (2008) 9670--9685.
\newblock \href {http://dx.doi.org/10.1016/j.jcp.2008.07.015}
  {\path{doi:10.1016/j.jcp.2008.07.015}}.

\bibitem{Ern2009}
A.~Ern, A.~F. Stephansen, P.~Zunino, {A discontinuous Galerkin method with
  weighted averages for advection–diffusion equations with locally small and
  anisotropic diffusivity}, IMA Journal of Numerical Analysis 29~(2) (2009)
  235--256.
\newblock \href
  {http://arxiv.org/abs/http://imajna.oxfordjournals.org/content/29/2/235.full.pdf+html}
  {\path{arXiv:http://imajna.oxfordjournals.org/content/29/2/235.full.pdf+html}},
  \href {http://dx.doi.org/10.1093/imanum/drm050}
  {\path{doi:10.1093/imanum/drm050}}.

\bibitem{Kirby2011}
R.~C. Kirby, {Fast simplicial finite element algorithms using Bernstein
  polynomials}, Numerische Mathematik 117~(4) (2011) 631--652.
\newblock \href {http://dx.doi.org/10.1007/s00211-010-0327-2}
  {\path{doi:10.1007/s00211-010-0327-2}}.

\bibitem{Kirby2012}
R.~C. Kirby, K.~T. Thinh, {Fast simplicial quadrature-based finite element
  operators using Bernstein polynomials}, Numerische Mathematik 121~(2) (2012)
  261--279.
\newblock \href {http://dx.doi.org/10.1007/s00211-011-0431-y}
  {\path{doi:10.1007/s00211-011-0431-y}}.

\bibitem{Blatt2010}
M.~Blatt, {A parallel algebraic multigrid method for elliptic problems with
  highly discontinuous coefficients}, Ph.D. thesis, Heidelberg, Univ., Diss.,
  2010 (2010).

\bibitem{Blatt2012}
M.~Blatt, O.~Ippisch, P.~Bastian, {A Massively Parallel Algebraic Multigrid
  Preconditioner based on Aggregation for Elliptic Problems with Heterogeneous
  Coefficients }\href {http://arxiv.org/abs/1209.0960}
  {\path{arXiv:1209.0960}}.

\bibitem{Hestenes1952}
M.~R. Hestenes, E.~Stiefel, {Methods of Conjugate Gradients for Solving Linear
  Systems}, Journal of Research of the National Bureau of Standards 49~(6).

\bibitem{Saad1986}
Y.~Saad, M.~H. Schultz, {GMRES: A Generalized Minimal Residual Algorithm for
  Solving Nonsymmetric Linear Systems}, SIAM Journal on Scientific and
  Statistical Computing 7~(3) (1986) 856--869.
\newblock \href {http://dx.doi.org/10.1137/0907058}
  {\path{doi:10.1137/0907058}}.

\bibitem{Bastian2008a}
P.~Bastian, M.~Blatt, A.~Dedner, C.~Engwer, R.~Kl\"{o}fkorn, R.~Kornhuber,
  M.~Ohlberger, O.~Sander, {A generic grid interface for parallel and adaptive
  scientific computing. {Part II}: implementation and tests in {DUNE}},
  Computing 82~(2-3) (2008) 121--138.
\newblock \href {http://dx.doi.org/10.1007/s00607-008-0004-9}
  {\path{doi:10.1007/s00607-008-0004-9}}.

\bibitem{Fog2016}
{VCL} {C++ Vector Class Library},
  \url{http://www.agner.org/optimize/\#vectorclass}.

\bibitem{Press2007}
W.~H. Press, S.~A. Teukolsky, W.~T. Vetterling, B.~P. Flannery, {Numerical
  Recipes 3rd Edition: The Art of Scientific Computing}, Cambridge University
  Press, 2007.

\bibitem{Saad1993}
Y.~Saad, {A flexible inner-outer preconditioned GMRES algorithm}, SIAM Journal
  on Scientific Computing 14~(2) (1993) 461--469.
\newblock \href {http://dx.doi.org/10.1137/0914028}
  {\path{doi:10.1137/0914028}}.

\bibitem{Strakos2005}
Z.~Strako{\v{s}}, P.~Tich{\'y}, {Error Estimation in Preconditioned Conjugate
  Gradients}, BIT Numerical Mathematics 45~(4) (2005) 789--817.
\newblock \href {http://dx.doi.org/10.1007/s10543-005-0032-1}
  {\path{doi:10.1007/s10543-005-0032-1}}.

\bibitem{Kempf2017}
D.~Kempf, {Generating performance-optimized finite element assembly kernels for
  dune-pdelab, Talk at the DUNE User Meeting 2017} (2017).

\end{thebibliography}

\end{document}